\documentclass[english]{article}
\usepackage[T1]{fontenc}
\usepackage[latin9]{inputenc}
\usepackage{verbatim}
\usepackage{amsmath}
\usepackage{amsthm}
\usepackage{amssymb}
\usepackage{esint}
\usepackage{varwidth}
\usepackage{url}
\usepackage{constants}
\usepackage{graphicx}
\usepackage{float}
\usepackage{hyperref}
\usepackage{mathabx} 
\usepackage{subcaption}
\usepackage{commath}
\usepackage{tikz}
\usetikzlibrary{decorations.pathreplacing}
\usepackage{authblk}
\usepackage{lipsum}

\def\RR{{\mathbb R}}

\hypersetup{
   colorlinks=true,%
   citecolor=black,%
   filecolor=black,%
   linkcolor=black,%
   urlcolor=blue
}

\makeatletter
\theoremstyle{plain}
\newtheorem{thm}{\protect\theoremname}
  \theoremstyle{plain}
  \newtheorem{lem}[thm]{\protect\lemmaname}
  \theoremstyle{plain}
  \newtheorem{corr}[thm]{\protect\corrolaryname}
  \theoremstyle{plain}
  \newtheorem{defn}[thm]{\protect\definitionname}
  \theoremstyle{plain}
  \newtheorem{prop}[thm]{\protect\propositionname}
  \theoremstyle{definition}

  \theoremstyle{plain}
  \newtheorem{conditions}[thm]{\protect\conditionsname}

\def \l {\left}
\def \r {\right}
\def \interior {\operatorname{int}} 

\newcommand{\wt}[1]{\widetilde{#1}}
\DeclareMathOperator*{\argmin}{\arg\!\min}
\DeclareMathOperator*{\argmax}{\arg\!\max}

\def \pert {\widetilde{Z}} 

\makeatother

\usepackage{babel}
  \providecommand{\definitionname}{Definition}
  \providecommand{\examplename}{Example}
  \providecommand{\lemmaname}{Lemma}
  \providecommand{\corrolaryname}{Corollary}
  \providecommand{\propositionname}{Proposition}
  \providecommand{\conditionsname}{Conditions}
\providecommand{\theoremname}{Theorem}

\begin{document}

\title{Sparse non-negative super-resolution --- simplified and stabilised}
\author[1]{Armin Eftekhari}
\author[2,3]{Jared Tanner}
\author[4]{Andrew Thompson}
\author[2]{Bogdan Toader}
\author[5]{Hemant Tyagi}
\affil[1]{EPFL, Switzerland}
\affil[2]{University of Oxford, UK}
\affil[3]{Alan Turing Institute, UK}
\affil[4]{National Physical Laboratory, UK}
\affil[5]{INRIA Lille - Nord Europe, France}

\maketitle

\begin{center}
  \large
  \textit{To David L. Donoho, a uniquely positive gentleman,\\
in celebration of his 60th birthday and with 
thanks for his inspiration and support.}
\end{center}

\begin{abstract}
The convolution of a discrete measure,
$x=\sum_{i=1}^ka_i\delta_{t_i}$, with a local window function, $\phi(s-t)$,
is a common model for a measurement device whose resolution is
substantially lower than that of the objects being observed.
Super-resolution concerns localising the point sources
$\{a_i,t_i\}_{i=1}^k$ with an accuracy beyond the essential
support of $\phi(s-t)$, typically from $m$ samples
$y(s_j)=\sum_{i=1}^k a_i\phi(s_j-t_i)+\delta_j$, where $\delta_j$ indicates
an inexactness in the sample value.   We consider the setting of
$x$ being non-negative and seek to characterise all non-negative
measures approximately consistent with the samples.  We first show
that $x$ is the unique non-negative measure consistent with the
samples provided the samples are exact, i.e. $\delta_j=0$, $m\ge 2k+1$
samples are available, and $\phi(s-t)$ generates a Chebyshev system.
This is independent of how close the sample locations are and {\em does
not rely on any regulariser beyond non-negativity}; as such, it extends
and clarifies the work 
by Schiebinger et al. in \cite{schiebinger2015superresolution}
and De Castro et al. in \cite{de2012exact}, who achieve the same
results but require a total variation regulariser, which we
show is unnecessary.

Moreover, we characterise non-negative solutions $\hat{x}$ consistent with the
samples within the bound $\sum_{j=1}^m \delta_j^2 \le \delta^2$.  Any such
non-negative measure is within ${\mathcal O}(\delta^{1/7})$ of the
discrete measure $x$ generating the samples in the generalised
Wasserstein distance.  Similarly, we show using somewhat different
techniques that the integrals of $\hat{x}$ and $x$ over
$(t_i-\epsilon,t_i+\epsilon)$ are similarly close, converging to one
another as $\epsilon$ and $\delta$ approach zero.  We also show how to
make these general results, for windows that form a Chebyshev system,
precise for the case of $\phi(s-t)$ being a Gaussian window.
The main innovation of these results is that non-negativity alone is
sufficient to localise point sources beyond the essential sensor
resolution and that, while regularisers such as total variation might
be particularly effective, they are not required in the non-negative setting.
\end{abstract}


\section{Introduction}

Super-resolution concerns recovering a resolution beyond the essential
size of the point spread function of a sensor.  For instance, a particularly
stylised example concerns multiple point sources which, because of the
finite \emph{resolution} or \emph{bandwidth} of the sensor, may not be  
visually distinguishable.  Various instances of this problem exist in
applications such as astronomy
\cite{puschmann2005super},  imaging in chemistry, medicine and
neuroscience
\cite{betzig2006imaging,hess2006ultra,rust2006sub,mccutchen1967superresolution,greenspan2009super,ekanadham2011neural,hell2009primer,tur2011innovation}, 
spectral estimation \cite{thomson1982spectrum,tang2015near},
geophysics \cite{khaidukov2004diffraction}, and system
identification \cite{shah2012linear}.  Often in these application much
is known about the point spread function of the sensor, or can be
estimated and, given such model information, it is possible to identify
point source locations with accuracy substantially below the essential
width of the sensor point spread function.  Recently there has been substantial interest from the
mathematical community in posing algorithms and proving
super-resolution guarantees in this setting, 
see for instance \cite{Candes2014,Tang2013,Fyhn2013,
Demanet2013,Duval2015,Denoyelle2017,Azais2015,bendory2016robust}. 
Typically these approaches borrow notions from compressed sensing
\cite{donoho2006compressed,candes2006robust,Candes2005}.  In
particular, the aforementioned contributions to super-resolution
consider what is known as the Total Variation norm minimisation over measures
which are consistent with the samples.  In this manuscript we
show first that, for suitable point spread 
functions, such as the Gaussian, any discrete non-negative
measure composed of $k$ point sources is uniquely defined from $2k+1$
of its samples, and moreover that this uniqueness is
independent of the separation between the point sources.  We then show
that by simply imposing non-negativity, which is typical in many
applications, any non-negative measure suitably consistent with
the samples is similarly close to the discrete non-negative measure
which would generate the noise free samples.  These results
substantially simply results by
\cite{schiebinger2015superresolution,de2012exact} and show that, while
regularisers such as Total Variation may be particularly effective,
in the setting of non-negative point sources such regularisers are not
necessary to achieve stability.

\subsection{Problem setup}

Throughout this manuscript we consider non-negative measures in relation to
discrete measures.  To be concrete, let $x$ be a $k$-discrete 
\emph{non-negative} Borel measure supported on the 
interval $I=[0,1] \subset \mathbb{R}$, given by
\begin{equation}
  x=\sum_{i=1}^k a_{i}\cdot \delta_{t_{i}}\quad
  \mbox{with}\;\;a_i>0 \quad\mbox{ and } t_i\in int(I)\quad\mbox{ for
    all } i.
  \label{eq:def of x}
\end{equation}
Consider also real-valued and {continuous} functions $\{\phi_{j}\}_{j=1}^m$
and let $\{y_{j}\}_{j=1}^m$ be the possibly noisy measurements
collected from $x$ by convolving against sampling functions $\phi_j(t)$: 
\begin{equation}
  y_j = \int_I \phi_j(t) x(\dif t) + \delta_j
  = \sum_{i=1}^{k} a_i \phi_j(t_i) + \delta_j,
  \label{eq:def_of_y}
\end{equation}
where $\delta_j$ with $\sum_{j=1}^m \delta_j^2 \le \delta^2$
can represent additive noise.  
Organising the $m$ samples from
\eqref{eq:def_of_y} in matrix notation by letting 
\begin{equation}\label{eq:Phi}
y:= [ y_1\;\cdots\; y_m]^T\in\RR^m,\quad\quad
\Phi(t):=[\phi_1(t)\;\cdots\;\phi_m(t)]^T\in\RR^m
\end{equation}
allows us to state the program we investigate:
\begin{equation}
  \text{find }z\ge 0\text{ subject to }
\left\|y-\int_I\Phi(t)z(\dif t)\right\|_2\le \delta',
  \label{eq:mai}
\end{equation}
with $\delta \leq \delta'$.
Herein we characterise non-negative measures consistent with
measurements \eqref{eq:def_of_y}  in relation to the discrete measure
\eqref{eq:def of x}.  That is, we consider any non-negative 
 Borel measure $z$ from the Program \eqref{eq:mai} 
\footnote{An equivalent formulation of Program \eqref{eq:mai} 
  minimises $\|y-\int_I \Phi(t)z(\dif t)\|_2$ over all 
  non-negative measures on $I$ (without any constraints). 
  In this context, however, we find it somewhat more intuitive 
  to work with Program \eqref{eq:mai}, particularly 
  considering the importance of the case $\delta=0$.
}
and show that any such $z$ is close to $x$ given by
\eqref{eq:def of x} in an appropriate metric, see Theorems
\ref{thm:main Gaussian}, \ref{thm:main Gaussian lambda},
\ref{thm:main noisy}, \ref{thm:noisy_expl} and  \ref{thm:grouping_sources}.
Note that when solving Program \eqref{eq:mai}, the measure $x$ is 
completely unkown, so neither the source locations $\{t_i\}_{i=1}^k$ 
and weights $\{a_i\}_{i=1}^k$ nor their number $k$ is known in advance.
Moreover, we show that the $x$ from \eqref{eq:def of x} is the unique
solution to Program \eqref{eq:mai} when $\delta'=0$; e.g. in the
setting of exact samples, $\eta_i=0$ for all $i$.
Program \eqref{eq:mai} is particularly notable in that there is no
regulariser of $z$ beyond imposing non-negativity and, rather than
specify an algorithm to select a $z$ which satisfies Program
\eqref{eq:mai}, we consider all admissible solutions.  
The aim of doing this is to highlight that non-negativity is the main
regulariser, especially in the noise-free setting. 
However, a practitioner solving this problem in the context of 
sparse measures would be advised to include additional regularisers 
such as the TV norm or a sparsity constraint in the context of 
non-convex methods to encourage sparisty, specifically in the noisy setting.

The admissible solutions of Program \eqref{eq:mai} are determined by the
source and sample locations, which we denote as 
\begin{equation}\label{eq:def_T_S}
  T=\{t_i\}_{i=1}^k\subset \interior(I)
  \quad\text{ and }\quad
  S=\{s_{j}\}_{j=1}^m\subseteq I
\end{equation}
respectively, as well as the particular functions $\phi_j(t)$ used to sample the
$k$-sparse non-negative measure $x$ from \eqref{eq:def of x}.  
Lastly,  we introduce the notions of minimum separation and sample
proximity, which we use to characterise solutions of Program \eqref{eq:mai}.
\begin{defn}
  \textbf{(Minimum separation and sample proximity) }\emph{\label{def:min_sep}}
  For finite $\tilde{T} = T \cup \{0, 1\} \subset I$, let $\Delta(T) > 0$
  be the minimum separation between the points in $T$ along with the endpoints of $I$, namely 
  \begin{equation}
    \Delta(T) = \min_{T_i, T_j \in \tilde{T}, i \ne j } |T_i - T_j|.
    \label{eq:def of sep}
  \end{equation}
  We define the sample proximity to be the number $\lambda \in (0,\frac12)$
  such that, for each source location $t_i$, there exists a closest sample location
  $s_{l(i)} \in S$ to $t_i$ with
  \begin{equation}
    |t_i - s_{l(i)}| \leq \lambda\Delta(T).
    \label{eq:closest_sample_thm}
  \end{equation}
\end{defn}

We describe the nearness of solutions to Program \eqref{eq:mai} in
terms of an additional parameter $\epsilon$ associated with intervals
around the sources $T$; that is we let $\epsilon\le\Delta(T)/2$ and
define intervals as:
\begin{equation}\label{eq:T_epsilon}
  T_{i,\epsilon}:=\left\{ t:\left|t-t_{i}\right|\le\epsilon\right\} \cap I,\qquad i\in[k],
  \quad
  T_{\epsilon}:=\bigcup_{i=1}^k T_{i,\epsilon},
\end{equation}
where $[k] = {1,2,\ldots,k}$, and set $T_{i,\epsilon}^{C}$ and $T_{\epsilon}^{C}$ to be the complements
of these sets with respect to $I$.
In order to make the most general result of Theorems \ref{thm:main
  noisy} and \ref{thm:noisy_expl}  more interpretable, we turn to
presenting them in Section \ref{sec:example} for the case of
$\phi_j(t)$ being shifted Gaussians.

  \subsection{Main results simplified to Gaussian window
    \label{sec:example}}

In this section we consider $\phi_j(t)$ to be shifted Gaussians with centres 
at the source locations $s_j$, specifically
\begin{equation}
  \phi_{j}(t)=g(t-s_{j})= e^{-\frac{(t-s_{j})^{2}}{\sigma^2}}.
  \label{eq:gaussian window}
\end{equation}
We might interpret \eqref{eq:gaussian window} as the ``point spread
function'' of the sensing mechanism being a Gaussian window and $s_j$
the sample locations in the sense that  
\begin{equation}\label{eq:conv}
  \int_I \phi_j(t) x(\dif t) = \int_I g(t-s_j) x(\dif t) = (g\star x)(s_j),
  \qquad \forall j\in [m],
\end{equation}
evaluates the ``filtered'' copy of $x$ at locations $s_j$
where $\star$ denotes convolution.

\begin{figure}[h]
  \centering
  \includegraphics[width=\textwidth]{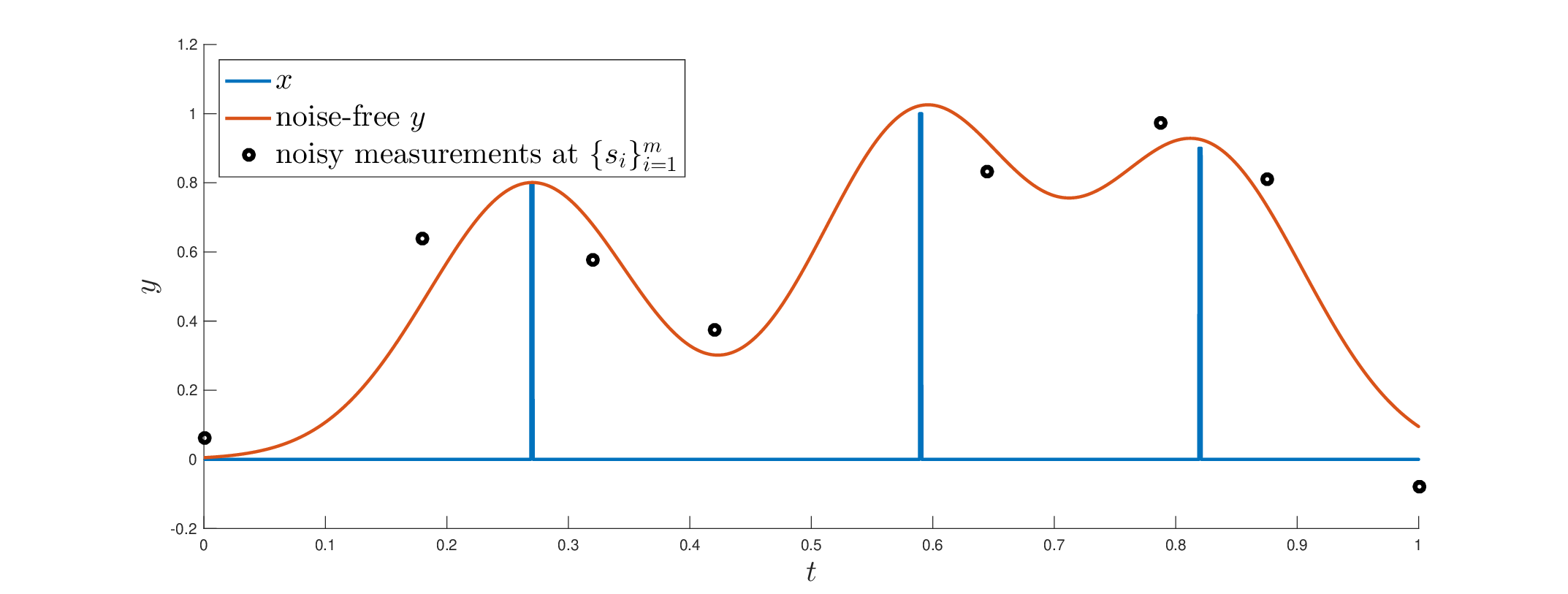}
  \caption{Example of discrete measure $x$ and measurements where $\phi_j(t) = \phi(t - s_j)$
  for $s_j\in S$ and the Gaussian kernel $\phi(t) = e^{-\frac{t^2}{\sigma^2}}$.}
  \label{fig:example_sig}
\end{figure}

As an illustration, Figure \ref{fig:example_sig} shows
the discrete measure $x$ in blue for $k=3$, the continuous function
$y(s)=(g\star x)(x)$ in red, and the noisy samples $y(s_j)$ at the
sample locations $S$ represented as the black circles. 

The conditions we impose to ensure stability of Program \eqref{eq:mai}
for $\phi(t)$ Gaussian as in \eqref{eq:gaussian window} are as follows:
\begin{conditions}
  \textbf{(Gaussian window conditions)} 
  When the window function is a Gaussian $\phi(t) = e^{-\frac{t^2}{\sigma^2}}$, 
  we require its width $\sigma$ and the source and sampling locations
  from \eqref{eq:def_T_S} to satisfy the following conditions:
  \begin{enumerate}
    \item Samples define the interval boundaries: $s_1=0$ and $s_m=1$, 
    \item Samples near sources: 
      for every $i\in[k]$, there exists a 
      pair of samples $s,s'\subset S$, one on each side 
      of $t_i$, such that $|s-t_i| \le \eta$ 
      and $s'-s \in [C_1\eta, C_2\eta]$ for some $C_1 \in (0,1]$
      and $C_2 \in [1,2)$ and $\eta \leq \sigma^2$ small enough; 
      which is quantified in Lemma \ref{lem:bounds on b for Gaussian}. 
    \item Sources away from the boundary: 
      $\sigma \sqrt{\log(1/\eta^3)} \le t_i,s_j \le 1-\sigma \sqrt{\log(1/\eta^3)}$ 
      for every $i\in[k]$ and $j\in [2:m-1]$, 
    \item Minimum separation of sources: $\sigma \leq \sqrt{2}$ 
      and $\Delta(T) > \sigma \sqrt{\log{(3 + \frac{4}{\sigma^2})}}$,
      where the minimum separation $\Delta(T)$ of the sources is defined
      in Definition \ref{def:min_sep}.
  \end{enumerate}
  \label{cond:conds thms 1 2}
\end{conditions}

The four properties in Conditions \ref{cond:conds thms 1 2} can be
interpreted as follows: Property 1 imposes that the sources are within
the interval defined by the minimum and maximum sample; Property 2
ensures that there is a pair of samples near each source which
translates into a sampling density condition in relation to the
minimum separation between sources and in particular requires the
number of samples $m\ge 2k+2$; 
Property 3 constrains the width of the Gaussian $\sigma$
through the sampling density $\eta$ (in particular 
$\sigma \ge \epsilon/\sqrt{\log(1/\eta^3)}$, where $\epsilon$ is the 
minimum distance between a source and the sampling boundary, 
which implies that having samples far away from the sources 
requires a wider point spread function. Conversely, a smaller distance
of the sources to the boundary $\epsilon$ allows a narrower point
spread function);
Properties 3 and 4 are technical conditions used in the proof to bound the 
eigenvalues of an associated stability matrix, and we expect may be improved,
though a condition constraining the size of $\sigma$ as compared to $\eta$ 
is likely necessary in the noisy setting as otherwise the samples in 
\eqref{eq:def_of_y} may have little dependence on the source locations $t_i$.

We can now present our main results on the robustness of Program
\eqref{eq:mai} as they apply to the Gaussian window; these are 
Theorem \ref{thm:main Gaussian}, which follows from 
Theorem \ref{thm:main noisy}, and 
Theorem \ref{thm:main Gaussian lambda}, which follows from 
Theorem \ref{thm:noisy_expl}.
However, before stating the stability results, it is important to note
that, in the setting of exact samples, $\eta_i=0$, the solution of
Program \eqref{eq:mai} is unique when $\delta'=0$.

\begin{prop}
\label{prop:unique_gaussian} \emph{\bf (Uniqueness of exactly
  sampled sparse non-negative measures for $\phi(t)$ Gaussian)} Let $x$ be a non-negative
  $k$-sparse discrete measure supported on $I$, see (\ref{eq:def of x}).  
  If $\delta=0$, $m\ge 2k+1$ and $\{\phi_{j}\}_{j=1}^m$
  are shifted Gaussians as in \eqref{eq:gaussian window}, 
  then $x$ is the unique solution of
  Program (\ref{eq:mai}) with $\delta'=0$. 
\end{prop} 

Proposition \ref{prop:unique_gaussian} states that Program
\eqref{eq:mai} successfully localises the $k$ 
impulses present in $x$ given only $2k+1$ measurements when
$\phi_j(t)$ are shifted Gaussians whose centres are in $I$. 
Theorems \ref{thm:main Gaussian} and \ref{thm:main Gaussian lambda}
extend this uniqueness condition to 
show that any solution to Program \eqref{eq:mai} with $\delta'>0$ is
proportionally close to the unique solution when $\delta'=0$.

\begin{thm} \emph{\textbf{(Wasserstein stability of Program
      \eqref{eq:mai} for $\phi(t)$ Gaussian)}} \label{thm:main Gaussian}
  Let $I=[0,1]$ and consider a $k$-sparse non-negative measure $x$ supported on $T\subset \interior(I)$. 
  Consider also an arbitrary increasing sequence $\{s_{j}\}_{j=1}^m\subset\mathbb{R}$ and, 
  for positive $\sigma$, let $\{\phi_{j}(t)\}_{j=1}^{m}$  
  be defined in \eqref{eq:gaussian window},
  which form $\Phi$ according to (\ref{eq:Phi}). 
  If $m\ge 2k+2$ and 
  Conditions \ref{cond:conds thms 1 2} hold, 
  then Program (\ref{eq:mai}) with $\delta'=\delta$ 
  is stable in the sense that 
  \begin{equation}\label{eq:gw_main}
    d_{GW}(x,\widehat{x}) 
    \le 
    F_1 \cdot \delta + \|x\|_{TV} \cdot \epsilon
  \end{equation}
  for all $\epsilon\le \Delta(T)/2$ where $d_{GW}$ is the generalised
  Wasserstein distance as defined in \eqref{eq:def of gen EMD} and
    the exact expression of $F_1 = 
  F_1(k,\Delta(T),\frac{1}{\sigma},\frac{1}{\epsilon},\eta)$  
  is given in the proof (see \eqref{eq:final F1}
  in Section \ref{sec:proof of main Gaussian}).
  In particular, for $\sigma < \frac{1}{\sqrt{3}}$ and
  $ \Delta(T) > \sigma \sqrt{\log{\frac{5}{\sigma^2}}}$,
  we have:
  \begin{equation}
    F_1(k,\Delta(T), \frac{1}{\sigma}, \frac{1}{\epsilon},\eta)
    <
    \frac{
      c_1 k C_1(\frac{1}{\epsilon})
    }{
      \eta \sigma^2 
    }
    \left[
      \frac{c_2}{\sigma^6 (1-3\sigma^2)^2}
    \right]^k,
    \label{eq:gaussian in particular}
  \end{equation}
  if
  \begin{equation}
    \eta \leq \min \left\{
        \frac{
          c_3 \sigma^6(1-3\sigma^2)
        }{
          (k+1)^{\frac32}
        },
        \frac{
          c_4 \bar{C}^{\frac16} \sigma^{\frac23}
        }{
          (k+1)^{\frac13}
        }
      \right\},
    \label{eq:cond eta partic}
  \end{equation}
  where $c_1,c_2,c_3,c_4$ are universal constants 
  and $C_1\left(\frac{1}{\epsilon}\right)$ is given
  by \eqref{eq:c1 eps def} in Section \ref{sec:bounds coeffs}
\end{thm} 

The central feature of Theorem \ref{thm:main Gaussian} is that the
proportionality to $\delta$ and $\epsilon$ of the Wasserstein distance between any solution to
Program \eqref{eq:mai} and the unique solution for $\delta'=0$ is of
the form \eqref{eq:gw_main}.
The particular form of $F_1(\cdot)$ is not believed to be
sharp; in particular, the exponential dependence on $k$ in
\eqref{eq:gaussian in particular} follows from 
bounding the determinant of a matrix similar to $\Phi$ 
(see \eqref{eq:1st minor}) by a lower bound on the
minimum eigenvalue to the $k^{th}$ power.  The scaling with respect to
$\sigma^{-2}$ is a feature of $\delta'$ in Program \eqref{eq:mai} not
being normalized with respect to $\|y\|_2$ which, for $T$ and $S$ fixed,
decays with $\sigma$ due to the increased localisation of the
Gaussian.  Note that the $\epsilon$ dependence is a feature of the
proof and the $\epsilon$ which minimises the bound in
\eqref{eq:gw_main} is proportional to $\delta$ to some power as
determined by $C_1(\epsilon^{-1})$ from \eqref{eq:gaussian in particular}.
Theorem \ref{thm:main Gaussian} follows from the more general
result of Theorem \ref{thm:main noisy}, whose proof is given 
in Section \ref{sec:theory} and the appendices.

As an alternative to showing stability of Program \eqref{eq:mai} in
the Wasserstein distance, 
we also prove in Theorem \ref{thm:main Gaussian lambda} 
that any solution to Program \eqref{eq:mai} is
locally consistent with the discrete measure in terms of local
averages over intervals $T_{i,\epsilon}$ as given in \eqref{eq:T_epsilon}.
Moreover, for Theorem \ref{thm:main Gaussian lambda}, we 
make Property 2 of Conditions \ref{cond:conds thms 1 2} more
transparent by using the sample proximity $\lambda\Delta(T)$ from
Definition \ref{def:min_sep}; that is, $\eta$ defined in Conditions
\ref{cond:conds thms 1 2} is related to the sample proximity from
Definition \ref{def:min_sep} by $\lambda \Delta(T) \leq \eta/2$.

\begin{thm} \emph{\textbf{
    (Average stability of Program \eqref{eq:mai} for $\phi(t)$
    Gaussian: source proximity dependence)
  }} 
  \label{thm:main Gaussian lambda}
  Let $I = [0,1]$ and consider a k-sparse non-negative measure $x$
  supported on $T$ and sample locations $S$ as given in
  \eqref{eq:def_T_S} and 
  for positive $\sigma$, let $\{\phi_{j}(t)\}_{j=1}^{m}$ 
  as defined in \eqref{eq:gaussian window}.
  If the Conditions \ref{cond:conds thms 1 2} hold, 
  then, in the presence of additive noise, Program \eqref{eq:mai} 
  is stable in the sense that, for any solution $\hat{x}$ 
  of Program \eqref{eq:mai} with $\delta' = \delta$: 
  \begin{align}
    &\left|
      \int_{T_{i,\epsilon}} \hat{x}(\dif t) - a_i
    \right|
    \leq
    \left[
      (c_1 + F_2)
      \cdot \delta
      + c_2 \frac{\|\hat{x}\|_{TV} }{\sigma^2}
      \cdot \epsilon
    \right] F_3,
    \label{eq:main Gaussian lambda res}
    \\
    &
      \int_{T_{\epsilon}^C} \hat{x} (\dif t)
    \leq
    F_2 \cdot \delta,
    \label{eq:main Gaussian lambda res away}
  \end{align}
  where the exact expressions of 
  $F_2 = F_2(k,\Delta(T),\frac{1}{\sigma},\frac{1}{\epsilon})$
  and $F_3 = F_3(\Delta(T),\sigma,\lambda)$
  are given in the proof (see \eqref{eq:f2 final} in 
  Section \ref{sec:proof of main Gaussian lambda}), provided that 
  $\lambda$, $\Delta(T)$ and $\sigma$ satisfy \eqref{eq:condition_lambda}.
  In particular, for $\sigma < \frac{1}{\sqrt{3}}$,
  $ \Delta(T) > \sigma \sqrt{\log{\frac{5}{\sigma^2}}}$ and 
  $\lambda < 0.4$, we have 
  $F_3(\Delta(T),\sigma,\lambda) < c_5$ and:
  \begin{equation}
    F_2(k,\Delta(T),\frac{1}{\sigma}, \frac{1}{\epsilon}) <
    c_3 \frac{k C_2(\frac{1}{\epsilon})}{\sigma^2} 
    \left[
      \frac{c_4}{\sigma^6 (1-3\sigma^2)^2}
    \right]^k.
    \label{eq:in particular f2}
  \end{equation}
  Above, $c_1, c_2, c_3, c_4, c_5$ 
  are universal constants and $C_2(\frac{1}{\epsilon})$ 
  is given by \eqref{eq:c2 eps def} in Section \ref{sec:bounds coeffs}.
\end{thm}

Both Theorem \ref{thm:main Gaussian} and \ref{thm:main Gaussian lambda}
establish a notion of stability in the sense that they give similar bounds 
on the error (measured using a different metric in each theorem)
in terms of the magnitude $\delta$ of the noise and an additional 
term involving $\epsilon$.
However, despite the fact that their proofs make use of the
same bounds, they have some fundamental differences.  
While both \eqref{eq:gw_main} and \eqref{eq:main Gaussian lambda res}
have the same proportionality to $\delta$ and
$\epsilon$, the role of $\epsilon$ in particular differs substantially
in that Theorem \ref{thm:main Gaussian lambda} considers averages of
$\hat{x}$ over $T_{i,\epsilon}$.  Also different in their form is the
dependence on $\|x\|_{TV}$ and $\|\hat{x}\|_{TV}$ in Theorems
 \ref{thm:main Gaussian} and \ref{thm:main Gaussian lambda} 
respectively.  The presence of $\|\hat{x}\|_{TV}$ in Theorem
\ref{thm:main Gaussian lambda}  is a feature of the proof which we
expect can be removed and replaced with $\|x\|_{TV}$ by proving any
solution of Program \eqref{eq:mai} is necessarily bounded due to the
sampling proximity condition of Definition \ref{def:min_sep}.
It is also worth noting that \eqref{eq:main Gaussian lambda res}
avoids an unnatural $\eta^{-1}$ dependence present in \eqref{eq:gw_main}.
Theorem \ref{thm:main Gaussian lambda} follows from the more general
result of Theorem \ref{thm:noisy_expl}, whose proof is given 
in Section \ref{sec:proof of main Gaussian lambda}.


Lastly, we give a corollary of Theorems \ref{thm:main Gaussian}
and \ref{thm:main Gaussian lambda} where we show that, for $\delta>0$
but sufficiently small, one can equate the $\delta$ and $\epsilon$
dependent terms in Theorems \ref{thm:main Gaussian}
and \ref{thm:main Gaussian lambda} to show that their respective
errors approach zero as $\delta$ goes to zero.
\begin{corr}
  \label{corr:small eps}
  Under the conditions in Theorems \ref{thm:main Gaussian}
  and \ref{thm:main Gaussian lambda} and for
  $\sigma < \frac{1}{\sqrt{3}}$,
  $ \Delta(T) > \sigma \sqrt{\log{\frac{5}{\sigma^2}}}$ and
  $\lambda < 0.4$,
  there exists $\delta_0$ such that:
  \begin{align}
    &d_{GW}(x,\hat{x}) \le \bar{C}_1 \cdot \delta^{\frac17},
    \\
    &\left| \int_{T_{i,\epsilon}} \hat{x}(\dif t) - a_i \right|
    \le \bar{C}_2 \cdot \delta^{\frac16},
  \end{align}
  for all $\delta \in (0,\delta_0)$,
  where $\bar{C}_1$ and $\bar{C}_2$ are given in the proof in
  Section \ref{sec:proof of corr}.
\end{corr}

  \subsection{Organisation and summary of contributions}

\paragraph{Organisation:}

The majority of our contributions were presented in the context of
Gaussian windows in Section \ref{sec:example}.  These are particular
examples of a more general theory for windows that form a
\emph{Chebyshev system}, commonly abbreviated as \emph{T-system}, see
Definition \ref{def:(T-systems)-Real-valued-and}.  
A T-system  is a collection of continuous functions that loosely behave like algebraic monomials.
It is a general and widely-used concept in classical approximation theory \cite{karlin1966tchebycheff,karlin1968total,kreinMarkov} that has also found  applications in modern signal processing \cite{schiebinger2015superresolution,de2012exact}. 
The framework we
use for these more general results is presented in Section
\ref{sec:tsystems}, the results presented in Section \ref{sec:main_stability},
and their proof sketched in Section \ref{sec:theory}.  Proofs of the
lemmas used to develop the results are deferred to the appendices.

\paragraph{Summary of contributions:}

We begin discussing results for general window function $\phi$ with Proposition
\ref{prop:existence of dual}, which establishes that for exact samples,
namely $\delta=0$, $\{\phi_j\}_{j=1}^m$ a T-system, and from $m\ge 2k+1$
measurements, the unique solution 
to Program \eqref{eq:mai} with $\delta'=0$ 
is the $k$-sparse measure $x$ given in
\eqref{eq:def of x}. 
In other words, we show that the measurement operator $\Phi$ in \eqref{eq:Phi} is an injective map from $k$-sparse non-negative  measures on $I$ to $\mathbb{R}^m$ when $\{\phi_j\}_{j=1}^m$ form a T-system. 
No minimum separation between impulses is necessary here  and $\{\phi_j\}_{j=1}^m$  need only to be continuous. 
As detailed in Section \ref{sec:related}, Proposition
\ref{prop:existence of dual} is more general and its derivation is far
simpler and more intuitive than what the current literature offers.
Most importantly, no explicit regularisation is needed in Program
\eqref{eq:mai} to encourage sparsity: the solution is unique.

Our main contributions are given in Theorems \ref{thm:main noisy} and
\ref{thm:noisy_expl}, namely that solutions to Program \eqref{eq:mai} with 
$\delta'>0$ are proportionally close to the unique solution
\eqref{eq:def of x} with $\delta'=0$; these theorems consider nearness
in terms of the Wasserstein distance and local averages respectively. Furthermore,
Theorem \ref{thm:main noisy} allows $x$ to be a general non-negative
measure, and shows that solutions to Program \eqref{eq:mai} must be
proportional to both how well $x$ might be approximated by a
$k$-sparse measure, $\chi$, with minimum source separation $2\epsilon$,
and a $\delta$ proportional distance between $\chi$ and solutions to
Program \eqref{eq:mai}.  These theorems require $m\ge 2k+2$ and loosely-speaking
the measurement apparatus forms a T*-system, which is an extension of 
a T-system to allow the inclusion of an additional function which may
be discontinuous, and enforcing certain properties of minors of
$\Phi$. To derive the bounds in Theorems \ref{thm:main Gaussian} and
\ref{thm:main Gaussian lambda} we show that shifted Gaussians as given
in \eqref{eq:gaussian window} augmented with a particular piecewise
constant function form a T*-system.

Lastly, in Section \ref{sec:grouping}, we consider an extension of
Theorem \ref{thm:noisy_expl} where the minimum separation between
sources $\Delta(T)$ is smaller than $\epsilon$.  We extend the
intervals $T_{i,\epsilon}$ from \eqref{eq:T_epsilon} to 
$\tilde{T}_{i,\epsilon}$ in \eqref{grouping_partition}, where
intervals $T_{i,\epsilon}$ which overlap are combined.  The resulting 
Theorem \ref{thm:grouping_sources} establishes that, while sources
closer than $\epsilon$ may not be identifiable individually by Program
\eqref{eq:mai}, the local average over
$\tilde{T}_{i,\epsilon}$ of both $x$ in \eqref{eq:def of x} and 
any solution to Program \eqref{eq:mai} will be proportionally within
$\delta$ of each other. 

To summarise, the results and analysis in this work  simplify, generalise 
and extend the existing results for grid-free and non-negative super-resolution.  
These extensions follow by virtue of the non-negativity constraint in  
Program \eqref{eq:mai}, rather than the common approach based on the
TV norm as a sparsifying penalty.   We further put these results in
the context of existing literature in Section \ref{sec:related}.

  \subsection{Comparison with other techniques \label{sec:related}}

We show in Proposition~\ref{prop:existence of dual} that a \emph{non-negative} $k$-sparse discrete measure can be exactly reconstructed from $m\geq 2k+1$ samples (provided that the atoms form a $T$-system, a property satisfied by Gaussian windows for example) by solving a feasibility problem. This result is in contrast to earlier results in which a TV norm minimisation problem is solved. De Castro and Gamboa~\cite{de2012exact} proved exact reconstruction using TV norm minimisation, provided the atoms form a homogeneous T-system (one which includes the constant function)\footnote{
  Note that in \cite{de2012exact}, De Castro and Gamboa consider the more 
  general case of \textit{signed} measures, where the TV norm objective
  is required, and they focus on three extended examples:
  non-negative measures, generalised Chebyshev measures 
  and $\Delta$-interpolation. Therefore, our noise-free results in the
  present paper can only be compared with the results from Section 2 
  in \cite{de2012exact}, where the authors consider the non-negative measure
  reconstruction by solving a TV norm minimisation problem 
  over signed measures in the noise-free setting.
}. An analysis of TV norm minimisation based on T-systems was subsequently given by Schiebinger et al. in ~\cite{schiebinger2015superresolution}, where it was also shown that Gaussian windows satisfy the given conditions. We show in this paper that the TV norm can be entirely dispensed with in the case of non-negative super-resolution. Moreover, analysis of Program \eqref{eq:mai} is substantially simpler than its alternatives. In particular, Proposition \ref{prop:existence of dual} for noise-free super-resolution immediately follows from the standard results in the theory of T-systems. The fact that Gaussian windows form a T-system is immediately implied by well-known results in the T-system theory, in contrast to the heavy calculations involved in \cite{schiebinger2015superresolution}. 

While neither of the above works considers the noisy setting or model mismatch, Theorems~\ref{thm:main noisy} and~\ref{thm:noisy_expl} in our work show that solutions to the non-negative super-resolution problem which are both stable to measurement noise and model inaccuracy can also be obtained by solving a feasibility program. The most closely related prior work is by Doukhan and Gamboa~\cite{Doukhan1996}, in which the authors bound the maximum distance between a sparse measure and any other measure satisfying noise-corrupted versions of the same measurements. While~\cite{Doukhan1996} does not explicitly consider reconstruction using the TV norm, the problem is posed over probability measures, that is those with TV norm equal to one. Accuracy is captured according to the Prokhorov metric. It is shown that, for sufficiently small noise the Prokhorov distance between the measures is bounded by $\delta^c$, where $\delta$ is the noise level and $c$ depends upon properties of the window function. In contrast, we do not make any total variation restrictions on the underlying sparse measure, we extend to consider model inaccuracy and we consider different error metrics (the generalised Wasserstein distance and the local averaged error).

More recent results on noisy non-negative super-resolution all assume that an optimisation problem involving the TV norm is solved. Denoyelle et al.~\cite{Denoyelle2017} consider the non-negative super-resolution problem with a minimum separation $t$ between source locations. They analyse a TV norm-penalized least squares problem and show that a $k$-sparse discrete measure can be stably approximated provided the noise scales with $t^{2k-1}$, showing that the minimum separation condition exhibits a certain stability to noise.
In the gridded setting, stability results for noisy non-negative super-resolution were obtained in the case of Fourier convolution kernels in~\cite{morgenshtern2014stable} under the assumption that the spike locations satisfy a Rayleigh regularity property, and these results were extended to the case of more general convolution kernels in~\cite{bendory2017robust}. 

Super-resolution in the more general setting of \emph{signed} measures has been extensively studied. In this case, the story is rather different, and stable identification is only possible if sources satisfy some separation condition. The required minimum separation is dictated by the resolution of the sensing system, e.g., the Rayleigh limit of the optical system or the bandwidth of the radar receiver. Indeed, it is impossible to resolve extremely close sources with equal amplitudes of opposite signs; they nearly cancel out, contributing virtually nothing to the measurements. A non-exhaustive list of references is \cite{Donoho1992,Tang2013,Fyhn2013,
Demanet2013,Duval2015,Azais2015,bendory2016robust}. 

In Theorem \ref{thm:noisy_expl} we give an explicit dependence of the error 
on the sampling locations. This result relies on local windows, hence it requires
samples near each source, and we give a condition that this distance must satisfy.
The condition that there are samples near each source in order to guarantee reconstruction
also appears in a recent manuscript on sparse deconvolution \cite{Bernstein2017}. 
However, this work relies on the
minimum separation and differentiability of the convolution kernel, which we overcome
in Theorem \ref{thm:noisy_expl}.

\section{Stability of Program \eqref{eq:mai} to inexact samples for
  $\phi_j(t)$ T-systems}
  \label{sec:main results} 

The main results stated in the introduction, 
Theorems \ref{thm:main Gaussian} and \ref{thm:main Gaussian lambda}, 
are for Gaussian windows, which allows the results to omit technical
details of the more general results of Theorems \ref{thm:main
  noisy}-\ref{thm:grouping_sources}.   These more general results
apply to windows that form Chebyshev systems, see Definition
\ref{def:(T-systems)-Real-valued-and}, and an extension to
$T^*$-systems, see Definition \ref{def:(T-systems,-modified)-For}, 
which allows for explicit control of the stability of
solutions to Program \eqref{eq:mai}.  These Chebyshev systems and
other technical notions needed are introduced in
Section \ref{sec:tsystems} and our most general contributions are presented
using these properties in Section \ref{sec:main_stability}.

  \subsection{Chebyshev systems and sparse measures}\label{sec:tsystems}

Before establishing stability of Program \eqref{eq:mai} to inexact
samples, we show that solutions to
Program \eqref{eq:mai} with $\delta'=0$, that is with 
$y_i$ in \eqref{eq:def_of_y} having $\eta_i=0$, has $x$ from
\eqref{eq:def of x} as its unique solution once $m\ge 2k+1$.  This
result relies on $\phi_j(t)$ forming a Chebyshev system,
commonly abbreviated T-system \cite{karlin1966tchebycheff}.   
\begin{defn}
  \textbf{(Chebyshev, T-system \cite{karlin1966tchebycheff})}
  \label{def:(T-systems)-Real-valued-and}
  Real-valued and continuous functions $\{\phi_{j}\}_{j=1}^m$ form a T-system
  on the interval $I$ if the $m\times m$ matrix $[\phi_{j}(\tau_{l})]_{l,j=1}^m$
  is nonsingular for any increasing sequence $\{\tau_{l}\}_{l=1}^m\subset I$. 
\end{defn}
Example of T-systems include the monomials $\{1,t,\cdots,t^{m-1}\}$ on
any closed interval of the real line. In fact, T-systems generalise
monomials and in many ways preserve their  properties. For instance,
any ``polynomial'' $\sum_{j=1}^m b_j \phi_j$ of a T-system
$\{\phi_j\}_{j=1}^m$ has at most $m-1$ distinct zeros on $I$. Or,
given $m$ distinct points on $I$, there exists a unique polynomial in
$\{\phi_j\}_{j=1}^m$ that interpolates these points. Note also that
linear independence of $\{\phi_j\}$ is a necessary condition for
forming a T-system, but not sufficient. Let us emphasise that T-system
is a broad and general concept with a range of applications in
classical approximation theory and modern signal processing. In the
context of super-resolution for example, translated copies of the 
Gaussian window, as given in \eqref{eq:gaussian window}, and many
other measurement windows form a T-system on any interval.  
We refer the interested reader  
to \cite{karlin1966tchebycheff,kreinMarkov} for the role 
of T-systems in classical approximation theory and 
to \cite{pinkus1996spectral} for their relationship 
to \emph{totally positive kernels}.

  \subsubsection{Sparse non-negative measure uniqueness from exact samples}

Our analysis based on T-Systems has been inspired by the work by  
Schiebinger et al. \cite{schiebinger2015superresolution}, where the
authors use the property of 
T-Systems to construct the dual certificate 
for the spike deconvolution problem
and to show uniqueness of the solution to the TV norm minimisation 
problem without the need of a minimum separation.  
The theory of T-Systems has also been used in the same context  
by De Castro and Gamboa in \cite{de2012exact}. 
However, both \cite{schiebinger2015superresolution} 
and \cite{de2012exact} focus on the noise-free problem exclusively,
while we will extend 
the T-Systems approach to the noisy case as well, as we will see later.

Our work, in part, simplifies the prior analysis considerably by using readily available
results on T-Systems and we go one step further to show uniqueness of
the solution of the feasibility problem, which removes the need for TV
norm regularisation in the results of Schiebinger et al.
\cite{schiebinger2015superresolution}; this simplification in the
presence of exact samples is given in 
Proposition \ref{prop:existence of dual}.

\begin{prop}
\label{prop:existence of dual} \emph{\bf (Uniqueness of exactly
  sampled sparse non-negative measures)} Let $x$ be a non-negative
$k$-sparse discrete measure supported on $I$ as given in \eqref{eq:def of x}. 
Let $\{\phi_{j}\}_{j=1}^m$ form a T-system on $I$, and given $m\ge
2k+1$ measurements as in 
\eqref{eq:def_of_y}, then $x$ is the unique solution of Program
(\ref{eq:mai}) with $\delta' =0$.  
\end{prop} 
Proposition \ref{prop:existence of dual} states that Program
\eqref{eq:mai} successfully localises the $k$ 
impulses present in $x$ given only $2k+1$ measurements when
$\{\phi_j\}_{j=1}^m$ form a T-system on $I$. Note that
$\{\phi_j\}_{j=1}^m$ only need to be continuous and \emph{no} minimum
separation is required between the impulses.  Moreover, as discussed in
Section \ref{sec:related}, the noise-free analysis here is
substantially simpler as it avoids the introduction of the TV norm
minimisation and is more insightful in that it shows that it is
not the sparsifying property of TV minimisation which implies the
result, but rather it follows from the non-negativity constraint and 
the T-system property, see Section \ref{sec:In-the-Absence}.

\subsubsection{T*-systems in terms of source and sample
  configuration}

While Proposition \ref{prop:existence of dual} implies that T-systems
ensure unique non-negative solutions, more is needed to ensure
stability of these results to inexact samples; that is $\delta>0$.
This is to be expected as T-systems imply invertibility of the linear
system $\Phi$ in \eqref{eq:Phi} for any configuration of sources and samples
as given in \eqref{eq:def_T_S}, 
but do not limit the condition number of such a system.  We control
the condition number of $\Phi$ by imposing further 
conditions on the source and sample configuration, such as those stated in
Conditions \ref{cond:conds thms 1 2}, which is analogous to imposing conditions
that there exists a dual polynomial which is sufficiently bounded away
from zero in regions away from sources, see Section \ref{sec:main_stability}.  In
particular, we extend the notion of T-system in Definition
\ref{def:(T-systems)-Real-valued-and} to a T*-system which includes
conditions on samples at the boundary of the interval, 
additional conditions on the window function, and 
a condition ensuring that there exist
samples sufficiently near sources as given by the
notation \eqref{eq:T_epsilon} but stated in terms of a new 
variable $\rho$ so as to highlight its different role here.

\begin{defn}
  \textbf{(T*-system) }
  \label{def:(T-systems,-modified)-For}
  For an even integer $m$,
  real-valued functions $\{\phi_{j}\}_{j=0}^m$
  form a T*-system on $I=[0,1]$ if the following 
	holds for every $T = \{t_1,t_2,\ldots,t_k\} \subset I$
  when $\rho>0$ is sufficiently small.
  For any increasing sequence $\tau=\{\tau_{l}\}_{l=0}^m\subset I$
    such that  
    \begin{itemize}
      \item $\tau_0=0$, $\tau_m = 1$, 
      \item except exactly three points, namely $\tau_0$, $\tau_m$, 
        and say $\tau_{\underline{l}}\in\interior(I)$, the  other points 
        belong to $T_{\rho}$, 
      \item  every $T_{i,\rho}$ contains an even number of points,
    \end{itemize}
    we have that 
    \begin{enumerate}
      \item the determinant of the $(m+1)\times (m+1)$ 
        matrix  $M_{\rho}:=[\phi_{j}(\tau_{l})]_{l,j=0}^m$ is positive, and
      \item the magnitudes of all minors of $M_{\rho}$ along the row
        containing $\tau_{\underline{l}}$ approach zero at the same rate\footnote{A function $u: \mathbb{R} \rightarrow \mathbb{R}^{+}$ approaches zero at the rate
		$\rho^P$ when $u(\rho) = \Theta(\rho^P)$. See, for example \cite{Cormen2009}, page 44.}
        when $\rho\rightarrow 0$.
    \end{enumerate}
\end{defn}

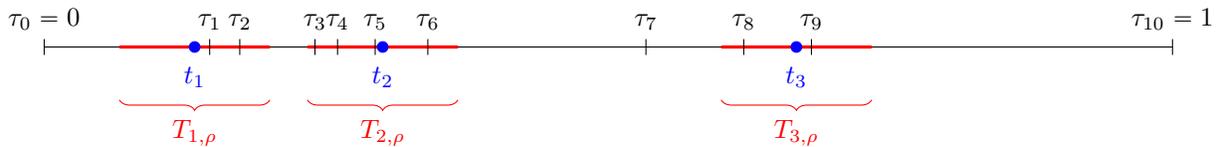
\begin{figure}[h]
  \centering
  \begin{tikzpicture}
    \draw
    (0,0) -- (15,0);

    \draw[very thick, red]
    (1,0) -- (3,0);
    \draw[red,decorate, decoration={brace,mirror,amplitude=4pt}]
    (1,-20pt) -- node[align=center,below=4pt] {$T_{1,\rho}$}
    (3,-20pt);

    \draw[very thick, red]
    (3.5,0) -- (5.5,0);
    \draw[red,decorate, decoration={brace,mirror,amplitude=4pt}]
    (3.5,-20pt) -- node[align=center,below=4pt] {$T_{2,\rho}$}
    (5.5,-20pt);

    \draw[very thick, red]
    (9,0) -- (11,0);
    \draw[red,decorate, decoration={brace,mirror,amplitude=4pt}]
    (9,-20pt) -- node[align=center,below=4pt] {$T_{3,\rho}$}
    (11,-20pt);

    \filldraw[blue]
    (2,0) circle (2pt) node[align=center, below=4pt] {$t_1$};

    \filldraw[blue]
    (4.5,0) circle (2pt) node[align=center, below=4pt] {$t_2$};

    \filldraw[blue]
    (10,0) circle (2pt) node[align=center, below=4pt] {$t_3$};

    \draw (0,3pt) node[align=center,above] {$\tau_0 = 0$}
    -- (0,-3pt);

    \draw (2.2,3pt) node[align=center,above] {$\tau_1$}
    -- (2.2,-3pt);

    \draw (2.6,3pt) node[align=center,above] {$\tau_2$}
    -- (2.6,-3pt);

    \draw (3.6,3pt) node[align=center,above] {$\tau_3$}
    -- (3.6,-3pt);

    \draw (3.9,3pt) node[align=center,above] {$\tau_4$}
    -- (3.9,-3pt);

    \draw (4.4,3pt) node[align=center,above] {$\tau_5$}
    -- (4.4,-3pt);

    \draw (5.1,3pt) node[align=center,above] {$\tau_6$}
    -- (5.1,-3pt);

    \draw (8,3pt) node[align=center,above] {$\tau_7$}
    -- (8,-3pt);

    \draw (9.3,3pt) node[align=center,above] {$\tau_8$}
    -- (9.3,-3pt);

    \draw (10.2,3pt) node[align=center,above] {$\tau_9$}
    -- (10.2,-3pt);

    \draw (15,3pt) node[align=center,above] {$\tau_{10} = 1$}
    -- (15,-3pt);

  \end{tikzpicture}
  \caption{Example of $\{\tau_l\}_{l=1}^m$ that satisfy the conditions
    in Definition \ref{def:(T-systems,-modified)-For} for $m=10$ and $k=3$.}
  \label{fig:taus_example}
\end{figure}

Let us briefly discuss T*-systems as an alternative to T-systems in
Definition \ref{def:(T-systems)-Real-valued-and}. The key property of
a T-system to our purpose is that an arbitrary polynomial
$\sum_{j=0}^m b_j \phi_j$ of a T-system $\{\phi_j\}_{j=0}^m$ on $I$
has at most $m$ zeros. Polynomials of a T*-system may not have such a 
property as T-systems allow arbitrary configurations of points $\tau$
while T*-systems only ensure the determinant in condition 1 of Definition
\ref{def:(T-systems,-modified)-For} be positive for configurations
where the majority of points in $\tau$ are paired in $T_{\rho}$.
However, as the analysis later shows, condition 1 in Definition
\ref{def:(T-systems,-modified)-For} is designed for constructing  dual
certificates for Program \eqref{eq:mai}.  We will also see later that
condition 2 in Definition \ref{def:(T-systems,-modified)-For} is meant
to exclude trivial polynomials that  do not qualify as dual
certificates. Lastly, rather than any increasing sequence
$\{\tau_l\}_{l\in[0:m]}\subset I$, Definition
\ref{def:(T-systems,-modified)-For} only considers subsets
$\tau$ that mainly cluster around the support $T$, whereas in our use
all but one entry in $\tau$ is taken from the set of samples $S$; this
is only intended to simplify the burden of verifying whether a family
of functions form a T*-system.  While the first and third bullet points in
Definition \ref{def:(T-systems,-modified)-For} require
that there need to be at least
two samples per interval $T_{i,\rho}$ as well as samples which define
the interval endpoints which gives a 
sampling complexity $m=2k+2$, we typically require $S$ to include
additional samples, $m>2k+2$, due to the location of $T$ being
unknown.   In fact, as $T$ is unknown, the third bullet point 
imposes a sampling density of $m$ being proportional to the inverse of
the minimum separation of the sources $\Delta(T)$.  The additional
point $\tau_{\underline{l}}$ is not taken from the set $S$, it instead
acts as a free parameter to be used in the dual certificate.
In Figure \ref{fig:taus_example}, we show an example of points
$\{\tau_l\}_{l=0}^{10}$ which  
satisfy the conditions in Definition
\ref{def:(T-systems,-modified)-For} for $k=3$ sources. 

%


We will state some of our more general stability results for solutions of
Program \eqref{eq:mai} in terms of the generalised Wasserstein 
distance \cite{piccoli2012generalized} between $x_1$ and $x_2$, 
both non-negative measures supported on $I$, defined as 
\begin{equation}
  d_{GW}\left(x_1,x_2\right)=\inf_{z_1,z_2} \l(\l\| x_1-z_1\r\|_{TV}+d_{W}\left(z_1,{z_2}\right) +\l\| {x_2}-{z_2}\r\|_{TV}\r),
\label{eq:def of gen EMD}
\end{equation}
where the infimum is over all non-negative Borel measures $z_1,{z}_2$
on $I$ such that $\|z_1\|_{TV}=\|{z}_2\|_{TV}$. Here, $\|z\|_{TV}=\int_I |z(\dif t)|$
is the \emph{total variation} of measure $z$, akin to the $\ell_1$-norm in finite dimensions,  
and $d_{W}$ is the standard
Wasserstein distance, namely 
\begin{equation}
  d_{W}\left(z_1,{z}_2\right)=\inf_{\gamma} \int_I\left|\tau_1-{\tau_2}\right|\cdot \gamma\left(\dif \tau_1,\dif{\tau_2}\right),
\label{eq:def of EMD}
\end{equation}
where the infimum is over all   measures $\gamma$
on $I\times I$ that produce $z_1$ and ${z}_2$ as marginals. In a sense, $d_{GW}$ extends $d_W$ to allow for calculating the distance between measures with different masses.
\footnote{
  In \cite{piccoli2012generalized}, the authors consider
  the p-Wasserstein distance, where popular choices of $p$
  are $1$ and $2$. In our work, we only use the 1-Wasserstein
  distance.
}

Moreover, in some of our most general results we consider the
extension to where $x$ need not be a discrete measure, see Theorem
\ref{thm:main noisy}.  In that setting, we introduce an intermediate
$k$-discrete measure which approximates $x$ in the $d_{GW}$
metric.  That is, given an integer $k$ and positive $\epsilon$, let
$x_{k,\epsilon}$ be a $k$-sparse  $2\epsilon$-separated measure 
supported on $T_{k,\epsilon} \subset \interior(I)$ of size $k$ and with
$\Delta(T_{k,\epsilon})\ge 2\epsilon$ such that, for $\beta \ge 1$,  
\begin{equation}
R(x,k,\epsilon):= d_{GW} (x,x_{k,\epsilon}) \le \beta \inf_{\chi} d_{GW}(x,\chi),
\label{eq:residual}
\end{equation}
where the infimum is over all $k$-sparse $2\epsilon$-separated
non-negative measures supported on $\interior(I)$ and the
parameter $\beta$ allows for near projections of $x$ onto the space of
$k$-sparse $2\epsilon$-separated measures. 

Lastly, we also assume that the measurement operator $\Phi$ in \eqref{eq:Phi}
 is Lipschitz continuous, namely there exists $L\ge 0$ such that 
\begin{equation}
\l\| \int_I \Phi(t) (x_1(\dif t)-x_2(\dif t))  \r\|_2
\le L \cdot d_{GW}(x_1,x_2),
\label{eq:Lipschitz assumption}
\end{equation}
for every pair of measures $x_1,x_2$ supported on $I$.

  \subsection{Stability of Program \eqref{eq:mai}}\label{sec:main_stability}

Equipped with the definitions of T and T*-systems, Definitions
\ref{def:(T-systems)-Real-valued-and} and
\ref{def:(T-systems,-modified)-For} respectively, we are able to characterise any
solution to Program \eqref{eq:mai} for $\phi_j(t)$ which form a
T-system and suitable source and sample configurations
\eqref{eq:def_T_S}.  We control the stability to 
inexact measurements by introducing two auxiliary functions in
Definition \ref{def:dual_separators}, which quantify the dual
polynomials $q(t)$ and $q^{\pi}(t)$ associated with 
Program \eqref{eq:mai} to be at least
$\bar{f}$ away from the necessary constraints for all values of $t$ at
least $\epsilon$ away from the sources. Specifically, for $F$ 
and $F^{\pi}$ defined below, we will require 
that $q(t) \geq F(t)$ and $q^{\pi}(t) \geq F^{\pi}(t)$ for all $t \in [0,1]$.

\begin{defn}
    \textbf{(Dual polynomial separators)}\label{def:dual_separators}
  Let $ f:\mathbb{R} \to \mathbb{R}_{+} $ be a bounded function with
  $f(0) = 0$, $\bar{f}, f_0, f_1$ be positive constants, 
  and $\{T_{i,\epsilon}\}_{i=1}^k$ the neighbourhoods as defined in
  (\ref{eq:T_epsilon}).  We then define 
  \begin{equation}
    F(t):=\begin{cases}
      f_0, & t=0, \\
      f_1, & t=1, \\
      f(t-t_{i}), & \text{when there exists }i\in[k]\text{ such that }t\in T_{i,{\epsilon}},\\
      \bar{f}, & \text{elsewhere on } int(I).
    \end{cases}
    \label{eq:function_F}
  \end{equation}
  Moreover, let $\pi\in \{\pm 1\}^k$ be an arbitrary sign pattern. We
  define $F^\pi$ as  
  \begin{equation}
    F^{\pi}(t):=
    \begin{cases}
      f_0, & t=0, \\
      f_1, & t=1, \\
      \pi_i - f\left(t-t_{i}\right), & \text{when there exists }
        i\in[k]\text{ such that }t\in T_{i,{\epsilon}},
        \\
      -\overline{f}, & \text{everywhere else on } int(I).
    \end{cases}
  \end{equation}
\end{defn}

We defer the introduction of dual polynomials $q$ and $q^{\pi}$ 
and the precise role of
the above dual polynomial separators to Section \ref{sec:theory}, but
state our most general results characterising the solutions to Program
\eqref{eq:mai} in terms of these separators.

\begin{thm}
  \label{thm:main noisy} \emph{\bf{(Wasserstein stability of Program
      \eqref{eq:mai} for $\phi_j(t)$ a T-system)}} Consider a
  non-negative measure $x$ supported on $\interior(I)=(0,1)$ and assume
  that the measurement operator $\Phi$   is $L$-Lipschitz, see
  (\ref{eq:Phi}) and  (\ref{eq:Lipschitz assumption}). Consider a
  $k$-sparse non-negative discrete measure $\chi$ supported on
  $T=\{t_i\}_{i=1}^k\subset \interior(I)$ and fix $\epsilon \le
  \Delta(T)/2$, see (\ref{eq:def of sep}), and
  consider functions $F(t)$ and $F^{\pi}(t)$ as defined in Definition
  \ref{def:dual_separators}.  For $m\ge2k+2$, suppose that
  \begin{itemize}
    \item  $\{\phi_{j}\}_{j=1}^m$ form a T-system on $I$,
    \item $\{F\}\cup\{\phi_{j}\}_{j=1}^m$ form a T*-system on $I$, and
    \item  $\{F^\pi\}\cup\{\phi_{j}\}_{j=1}^m$ form a T*-system on $I$ 
      for any sign pattern $\pi$.
  \end{itemize}
  Let $\widehat{x}$ be a solution of Program (\ref{eq:mai}) with 
  \begin{equation}
  \delta' =\delta+ L\cdot d_{GW}\left(x,\chi \right).
  \label{eq:deltap in thm}
  \end{equation}
    Then there exist vectors $b,\{b^\pi\}_\pi\subset \mathbb{R}^m$ such that 
  \begin{equation}
  d_{GW}\left( x,\widehat{x}\right)\le  \left(\left(6+\frac{2}{\bar{f}}\right)\|b\|_{2}+6\min_\pi \|b^{\pi}\|_{2}\right)\delta'+\epsilon\| \chi \|_{TV}+d_{GW}(x,\chi). 
  \label{eq:EMD thm}
  \end{equation}
  where the minimum is over all sign patterns $\pi$ and the vectors
  $b,\{b^\pi\}_\pi\subset\mathbb{R}^m$ above are the vectors of coefficients of the dual
  polynomials $q$ and $q^{\pi}$ associated with Program \eqref{eq:mai}, see Lemmas \ref{lem:dual}
  and \ref{lem:dual 3} in Section \ref{sec:theory} for their precise definitions. 
\end{thm}
%


Theorem \ref{thm:main Gaussian} follows from Theorem 
\ref{thm:main noisy} by considering $\phi_{j}(t)$ Gaussian as stated in
\eqref{eq:gaussian window} which is known to be a T-system
\cite{karlin1966tchebycheff}, and introducing Conditions
\ref{cond:conds thms 1 2} on the source and sample 
configuration \eqref{eq:def_T_S} such that the conditions of Theorem \ref{thm:main noisy} can be
proved and the dual coefficients $b$ and $b^{\pi}$ bounded; the
details of these proofs and bounds are deferred to Section
\ref{sec:theory} and the appendices.


The particular form of $F$ and $\{F^\pi\}_\pi$ in Theorem
\ref{thm:main noisy}, constant away from the support $T$ of
$x_{k,\epsilon}$, is purely to simplify the presentation and proofs.
Note also that the error $d_{GW}(x,\widehat{x})$ depends both
on the noise level $\delta$ and the residual $R(x,k,\epsilon)$, not
unlike the standard results in finite-dimensional sparse recovery and
compressed sensing
\cite{donoho2006compressed,candes2008introduction}. In particular,
when $\delta, \epsilon, R(x,k,\epsilon)\rightarrow 0$, we approach the
setting of Proposition \ref{prop:existence of dual}, where we have uniqueness
of $k$-sparse non-negative measures from exact samples.

Note that the noise level $\delta$ and the residual
$R(x,k,\epsilon)$ are not independent; that is, $\delta$ specifies
confidence in the samples and the model for how the samples are taken
while $R(x,k,\epsilon)$ reflects nearness to the model of $k$-discrete
measures.  Corollary \ref{corr:small eps} show that the parameter
$\epsilon$ can be removed, for $\phi_j(t)$ shifted Gaussians, in the
setting where $x$ is $k$-discrete, 
that is $R(x,k,\epsilon)=0$, in which case $d_{GW}(x,\hat{x})$ is
bounded by $\mathcal{O}(\delta^{1/7})$.

The more general variant of Theorem \ref{thm:main Gaussian lambda}
follows from Theorem \ref{thm:noisy_expl} by introducing alternative
conditions on the source and sample configuration and omitting the
need for the functions $F^{\pi}$, which is the cause of the unnatural
$\eta^{-1}$ dependence in Theorem \ref{thm:main Gaussian}. 

\begin{thm}
  \label{thm:noisy_expl} 
\emph{\bf{(Average stability for Program \eqref{eq:mai} for
    $\phi_j(t)$ a T-system)}} 
  Let $ \hat{x} $ be a solution of Program \eqref{eq:mai} and
  consider the function $F(t)$ as defined in Definition \ref{def:dual_separators}.
  Suppose that:
  \begin{itemize}
    \item  $\{\phi_{j}\}_{j=1}^m$ form a T-system on $I$,
    \item $\{F\}\cup\{\phi_{j}\}_{j=1}^m$ form a T*-system on $I$, and
    \item $\Delta=\Delta(T)$ and $\lambda = \lambda_0 \in (0,1/2)$ 
      from Definition \ref{def:min_sep} satisfy
    \begin{equation}
      \label{eq:condition_lambda}
      \phi(\lambda \Delta)
      =
      \phi(\Delta - \lambda \Delta) +
      \phi(\Delta + \lambda \Delta) +
      \frac{1}{\Delta}
        \int_{\Delta - \lambda\Delta}^{1/2 - \lambda\Delta} \phi(x) \dif x +
      \frac{1}{\Delta}
        \int_{\Delta + \lambda\Delta}^{1/2 + \lambda\Delta} \phi(x) \dif x.
    \end{equation} 
  \end{itemize}
  Then, for any $ \epsilon \in (0,\Delta/2) $ and for all $i \in [k]$,
  \begin{align}
    &\left|
      \int_{T_{i,\epsilon}} \hat{x}(\dif t) - a_i 
    \right|
    \leq
    \left(
      2 \left( 1 + \frac{\phi^{\infty}\|b\|_2}{\bar{f}} \right) 
        \cdot \delta
        + L \|\hat{x}\|_{TV} \cdot \epsilon
    \right)
    \sum_{j=1}^{k} (A^{-1})_{ij},
    \label{eq:bound_near_support}
    \\
    &\int_{T_\epsilon^C}  \hat{x}(\dif t) 
    \le \frac{2\|b\|_2 \delta}{\bar{f}},
    \label{eq:local avg away}
  \end{align}
  where:
  \begin{itemize}
    \item $b \in \mathbb{R}^m $ is the same vector of coefficients of the dual 
      certificate $q$ as in Theorem \ref{thm:main noisy}
      and $\bar{f}$ is given in Definition \ref{def:dual_separators}, 
      which is used to construct the dual certificate $q$, as described in 
      Lemma \ref{lem:dual} in Section \ref{sec:theory},

    \item $\phi^{\infty} = \max_{s,t \in I} |\phi(s-t)|$, 

    \item $L$ is the Lipschitz constant of $\phi$,

    \item $A \in \mathbb{R}^{k \times k}$ is the matrix
      \begin{equation}
        \label{eq:matrix_A}
		    A = 
		    \begin{bmatrix}
		      |\phi_1(t_1)|  & -|\phi_1(t_2)| & \ldots & -|\phi_1(t_k)| \\
		      -|\phi_2(t_1)| & |\phi_2(t_2)|  & \ldots & -|\phi_2(t_k)| \\
		      \vdots         & \vdots         & \ddots & \vdots         \\
		      -|\phi_k(t_1)| & -|\phi_k(t_2)| & \ldots & |\phi_k(t_k)|  \\
		    \end{bmatrix},
		  \end{equation}
      with $\phi_i(t_i)=\phi(t_i-s_{l(i)})$ evaluated at $s_{l(i)}$ 
      as defined in \eqref{eq:closest_sample_thm}.
  \end{itemize}
\end{thm}

Theorem \ref{thm:noisy_expl} bounds the difference between the average
over the interval $T_{i,\epsilon}$ of any solution to Program
\eqref{eq:mai} and the discrete measure whose average is simply $a_i$.
The condition on $\lambda$ to satisfy \eqref{eq:condition_lambda}
is used to ensure the matrix from \eqref{eq:matrix_A} is strictly diagonally dominant.
It relies on the windows $\phi_j(t)$ being sufficiently localised
about zero.  Though Theorem \ref{thm:noisy_expl}
explicitly states that the location of the closest samples to each source is
less than $\lambda_0\Delta(T)$, this can be achieved without knowing
the locations of the 
sources by placing the samples uniformly at interval
$2\lambda_0\Delta(T)$ which gives a sampling complexity of
$m=(2\lambda_0\Delta(T))^{-1}$. Lastly, a similar bound on the
integral of $\hat{x}$ over $T_{\epsilon}^C$ is given by Lemma
\ref{lem:dual} in Section \ref{sec:theory}.

\subsubsection{Clustering of indistinguishable sources}
\label{sec:grouping}

Theorems \ref{thm:main noisy} and \ref{thm:noisy_expl} give uniform
guarantees for all sources in terms of the minimum separation
condition $\Delta(T)$, which measures the worst proximity of sources.
One might imagine that, for example, if all but two sources are
sufficiently well separated, then Theorem \ref{thm:noisy_expl} might
hold for the sources that are well separated; moreover, assuming
$\epsilon$ is fixed, then if two sources 
$ t_i$ and $t_{i+1}$ with magnitudes $a_i$ and $a_{i+1}$ 
are closer than $2\epsilon$, 
namely $|t_i - t_{i+1}| < 2\epsilon$, we might imagine that a variant of
Theorem \ref{thm:noisy_expl} might hold but with sources $t_i$ and
$t_{i+1}$ approximated with source $ t_{\xi}$ near $t_i$ and
$t_{i+1}$ and with $ a_{\xi} = a_i + a_{i+1}$. 

In this section we extend Theorem \ref{thm:noisy_expl} to this setting
by considering $\epsilon$ fixed and alternative intervals 
$ \{\tilde{T}_{i,\epsilon} \}_{i=1}^{\tilde{k}} $ a partition of
${T_{\epsilon}}$ 
such that each $\tilde{T}_{i,\epsilon}$ contains a group of consecutive 
sources $ t_{i1},\ldots,t_{ir_i} $ 
(with weights $a_{i1},\ldots,a_{ir_i}$ respectively) which are within
at most $2\epsilon$ of each other. Define
\begin{equation}
  \label{grouping_partition}
  \tilde{T}_{i,\epsilon} = \bigcup_{l=1}^{k_i} T_{il,\epsilon}, 
  \quad
  \text{ where }
  t_{il} \in T_{il,\epsilon}
  \quad
  \text{ and }
  \quad
  | t_{i l+1} - t_{il} | < 2\epsilon,
  \quad
  \forall l \in [k_i-1],
\end{equation}
for $\sum_{i=1}^{\tilde{k}} k_i = k$,
so that we have
\begin{equation}
  T_{\epsilon} = \bigcup_{i=1}^{\tilde{k}} \tilde{T}_{i,\epsilon}
  \quad
  \text{ and }
  \quad
  \tilde{T}_{i,\epsilon} \bigcap \tilde{T}_{j,\epsilon} \ne \emptyset,
  \quad
  \forall i \ne j.
\end{equation}

\begin{thm}
  \label{thm:grouping_sources}
\emph{\bf{(Average stability for Program \eqref{eq:mai}: grouped
    sources)}}
Let $\hat{x}$ be a solution of Program \eqref{eq:mai} and $I=[0,1]$
be partitioned as described by \eqref{grouping_partition}.
  If the samples are placed uniformly at interval $2\lambda_0\epsilon$ 
  where $\lambda=\lambda_0$
  satisfies \eqref{eq:condition_lambda} with $\Delta=2\epsilon$, 
  then there exist $\{\xi_i\}_{i \in [\tilde{k}]}$ with $\xi_i \in \tilde{T}_{i,\epsilon}$
  such that
  \begin{equation}
    \left|
      \int_{\tilde{T}_{i,\epsilon}} \hat{x}(\dif t) - 
      \sum_{r=1}^{k_i} a_{ir} 
    \right|
    \leq
    \left(
    2 \left( 1 + \frac{\phi^{\infty}\|b\|_2}{\bar{f}} \right) 
      \cdot \delta
      + (2 k-1) L \|\hat{x}\|_{TV} \cdot \epsilon
    \right)
    \sum_{j=1}^{k} (\tilde{A}^{-1})_{ij},
  \end{equation}
  where the constants are the same as in \eqref{thm:noisy_expl} and the 
  matrix $ \tilde{A} \in \mathbb{R}^{\tilde{k} \times \tilde{k}} $ is
  \begin{equation*}
    \tilde{A} = 
    \begin{bmatrix}
      |\phi_1(\xi_1)|  & -|\phi_1(\xi_2)| & \ldots & -|\phi_1(\xi_{\tilde{k}})| \\
      -|\phi_2(\xi_1)| & |\phi_2(\xi_2)|  & \ldots & -|\phi_2(\xi_{\tilde{k}})| \\
      \vdots         & \vdots         & \ddots & \vdots         \\
      -|\phi_{\tilde{k}}(\xi_1)| & -|\phi_{\tilde{k}}(\xi_2)| & \ldots & |\phi_{\tilde{k}}(\xi_{\tilde{k}})|  \\
    \end{bmatrix}.
  \end{equation*}
\end{thm}

Note that Lemma \ref{lem:dual} still holds if we replace any group of sources from 
an interval $\tilde{T}_{i,\epsilon}$ with some $\xi_i \in \tilde{T}_{i,\epsilon}$,
so the bound from Lemma \ref{lem:dual} on $T_{\epsilon}^C$ remains
valid without modification.

As an exemplar source location where Theorem
\ref{thm:grouping_sources} might be applied, consider the situation
where the $k$ source locations comprising $T$ are drawn uniformly at
random in $(0,1)$, where  we have that
(from \cite{feller1968} page 42, Exercise 22)
\begin{equation*}
  P(\Delta(T) > \theta) = [1 - (k+1)\theta]^k,
  \quad
  \theta \in \left[0, \frac{1}{k+1}\right].
\end{equation*}
Then, the cumulative distribution function is
\begin{equation*}
  F(\theta) = P(\Delta(T) \leq \theta) = 1 - [1 - (k+1)\theta]^k,
\end{equation*}
and so the distribution of $\Delta(T)$ is
\begin{equation*}
  f(\theta) = F'(\theta) = (k+1)k[1 - (k+1)\theta]^{k-1},
\end{equation*}
with an expectation of
\begin{equation}\label{eq:expected_delta}
  E(\Delta(T)) = \int_0^{\frac{1}{k+1}} 
    P(\Delta(T) > \theta) \dif \theta 
  = \frac{1}{(k+1)^2}.
\end{equation}
That is, for $x$ from \eqref{eq:def of x} with sources $T$ drawn uniformly at
random in $(0,1)$, the expected value of $\Delta(T)$ is given by
\eqref{eq:expected_delta} and, in Theorems \ref{thm:main noisy} and
\ref{thm:noisy_expl}, the corresponding number of samples $m$ would scale
quadratically with the number of sources $k$ due to the scaling of
$m\sim \Delta(T)^{-1}$.  Alternatively, Theorem
\ref{thm:grouping_sources} allows meaningful results for $m$
proportional to $k$ by grouping the sources that are within $k^{-2}$
of one another.


\section{Dual polynomials for stability of non-negative measures}
  \label{sec:theory}

The results in Section \ref{sec:main results} 
 are developed by establishing dual polynomials of Program
 \eqref{eq:mai} which are non-negative except at the source locations,
 which implies that the solution to Program \eqref{eq:mai} is unique when $\delta'=0$, see
 Proposition \ref{prop:existence of dual}, and then showing that the
 dual polynomials are sufficiently non-negative away from the
 source locations and using this property to develop 
Theorems \ref{thm:main noisy}, \ref{thm:noisy_expl}, and \ref{thm:grouping_sources}.
In this section we state the key lemmas used to prove the
aforementioned results and then bound the quantities involving the
dual polynomials in order to establish 
Theorems \ref{thm:main Gaussian} and \ref{thm:main Gaussian lambda}
for the case of Gaussian windows.

\subsection{Uniqueness of non-negative sparse measures from exact
  samples: proof of Proposition \ref{prop:existence of dual} \label{sec:In-the-Absence}}

Proposition \ref{prop:existence of dual} states that if $x$ is a
non-negative $k$-sparse discrete measure supported on $I$, see
(\ref{eq:def of x}), provided $m\ge 2k+1$ and
$\{\phi_{j}\}_{j=1}^m$ are a T-system, then $x$ is the unique
non-negative solution to Program \eqref{eq:mai} with $\delta'=0$. This
follows from the existence of a dual polynomial as stated in Lemma
\ref{lem:dual noise free}, the proof of which is given in Appendix
\ref{sec:Proof-of-Lemma noise free dual cert}.
\begin{lem}{\bf (Dual polynomial and uniqueness of non-negative sparse
    measure equivalence)}
  \label{lem:dual noise free} Let $x$ be a non-negative $k$-sparse discrete 
  measure supported on $I$, see (\ref{eq:def of x}).  Then, $x$ is the unique solution of Program
  (\ref{eq:mai}) with $\delta'=0$ if 
  \begin{itemize}
    \item the  $k\times m$ matrix $[\phi_j(t_i)]_{i=1,j=1}^{i=k,j=m}$ is full rank, and  
    \item there exist  real coefficients $\{b_j\}_{j=1}^m$ 
      and $q(t)=\sum_{j=1}^m b_j\phi_j(t)$ such that   
      $q(t)$ is non-negative on $I$ and vanishes only on $T$. 
  \end{itemize}
\end{lem}

\begin{figure}[h]
  \centering
  \includegraphics[width=0.6\textwidth]{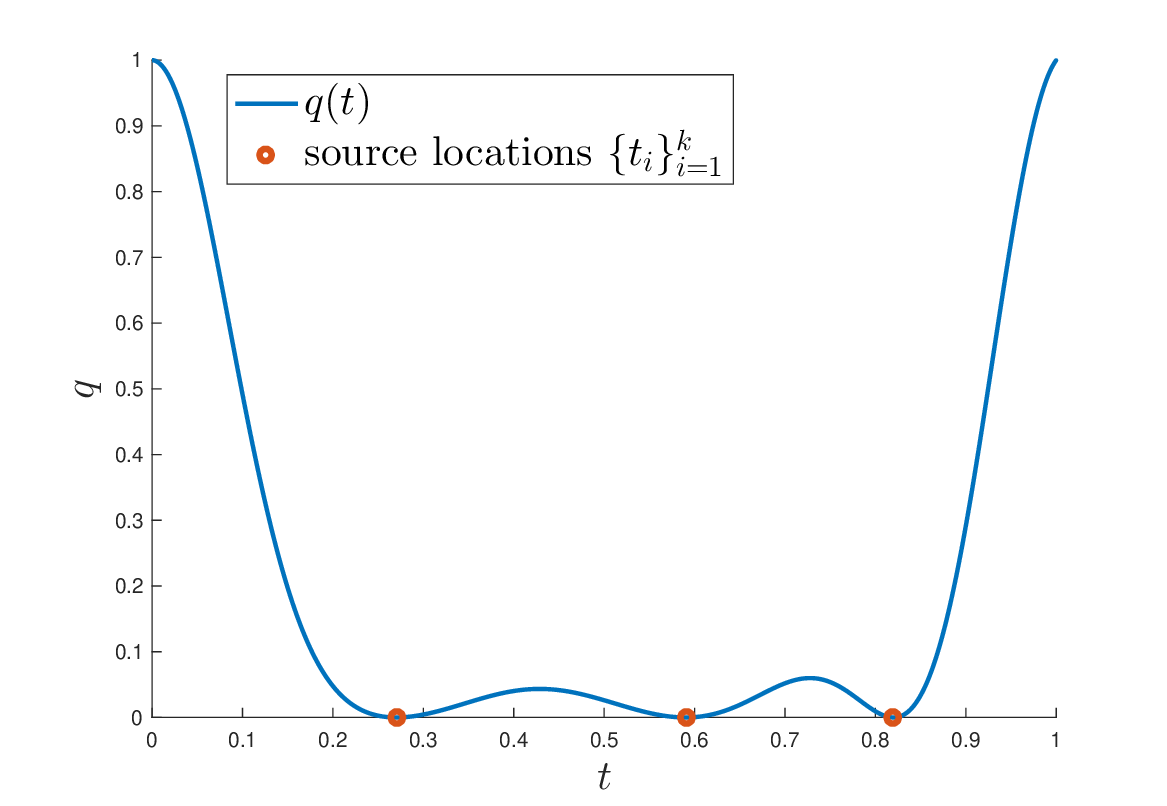}
  \caption{Example of dual certificate $q(t)$ required in Lemma \ref{lem:dual noise free}. Here, we 
  have $k=3$ and $t_1 = 0.27, t_2=0.59 $ and $ t_3 = 0.82$.}
  \label{fig:example_dual}
\end{figure}

Figure \ref{fig:example_dual} shows an example of such a dual
certificate using Gaussian $\phi_j$ as defined in \eqref{eq:gaussian
  window}.  It remains to show that such a dual polynomial exists.  To
do this, we employ the concept of T-system introduced in Definition
\ref{def:(T-systems)-Real-valued-and}. Of particular interest to us is
Theorem 5.1 in \cite{karlin1966tchebycheff}, slightly simplified
below, which immediately proves Proposition \ref{prop:existence of
  dual}.   
\begin{lem}
{\bf (Dual polynomial existence for T-systems)
  \cite[Theorem 5.1, pp. 28]{karlin1966tchebycheff}}
\label{lem:karlin-1}With $m\ge2k+1$, suppose that $\{\phi_{j}\}_{j=1}^m$
form a T-system on $I$. Then there exists a polynomial
$q(t)=\sum_{j=1}^m b_{j}\phi_{j}(t)$ that is non-negative on $I$ and
vanishes only on $T$.  
\end{lem}

\subsection{Stabilising dual polynomials for non-negative sparse measures: proof of Theorem \ref{thm:main noisy}} \label{sec:In-the-Presence}

We develop the proof of Theorem \ref{thm:main noisy} by using a dual
polynomial analogous to that in Lemma \ref{lem:dual noise free}, but
with further guarantees that away from the sources the dual polynomial
must be sufficiently bounded away from the constraint bounds.
However, first, let us bring the generality of Theorem \ref{thm:main
  noisy} to discrete measures by introducing an intermediate measure
$\chi$ which is $k$-discrete and whose support is $2\epsilon$
separated.  Noting that the measurement
operator $\Phi$ is $L$-Lipschitz, see \eqref{eq:Lipschitz assumption},
and using the triangle inequality, it follows that 
\begin{align}
\l\|y-\int_I \Phi(t) \chi (\dif t)\r\|_{2} & \le \l\|y-\int_I  \Phi(t) x(\dif t)\r\|_{2}+ \l\|\int_I  \Phi(t) (x(\dif t)-\chi (\dif t)) \r\|_{2} \nonumber \\
 & \le \delta+ L\cdot d_{GW}\left(x,\chi \right) =:\delta'. 
 \label{eq:def of deltap}
\end{align}
Therefore, a solution $\widehat{x}$ of Program \eqref{eq:mai} with
$\delta'$ specified above can be considered as an estimate of
$\chi$. In the rest of this section, we first bound the error
$d_{GW}(\chi,\widehat{x})$ and then use the triangle inequality to
control  $d_{GW}(x,\widehat{x})$.

\begin{figure}[h]
  \centering
  \begin{subfigure}{0.45\textwidth}
    \centering
    \includegraphics[width=\textwidth]{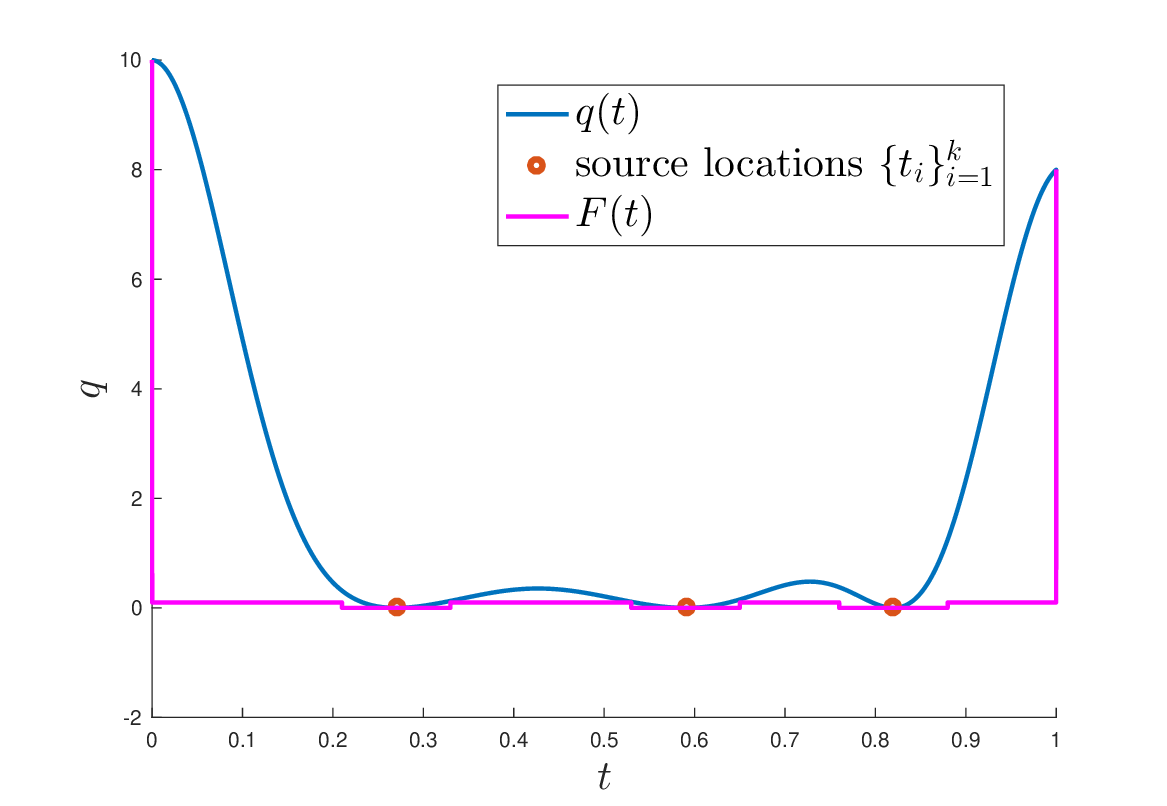}
    \caption{}
    \label{fig:example_dual_G_big}
  \end{subfigure}
  \begin{subfigure}{0.45\textwidth}
    \centering
    \includegraphics[width=\textwidth]{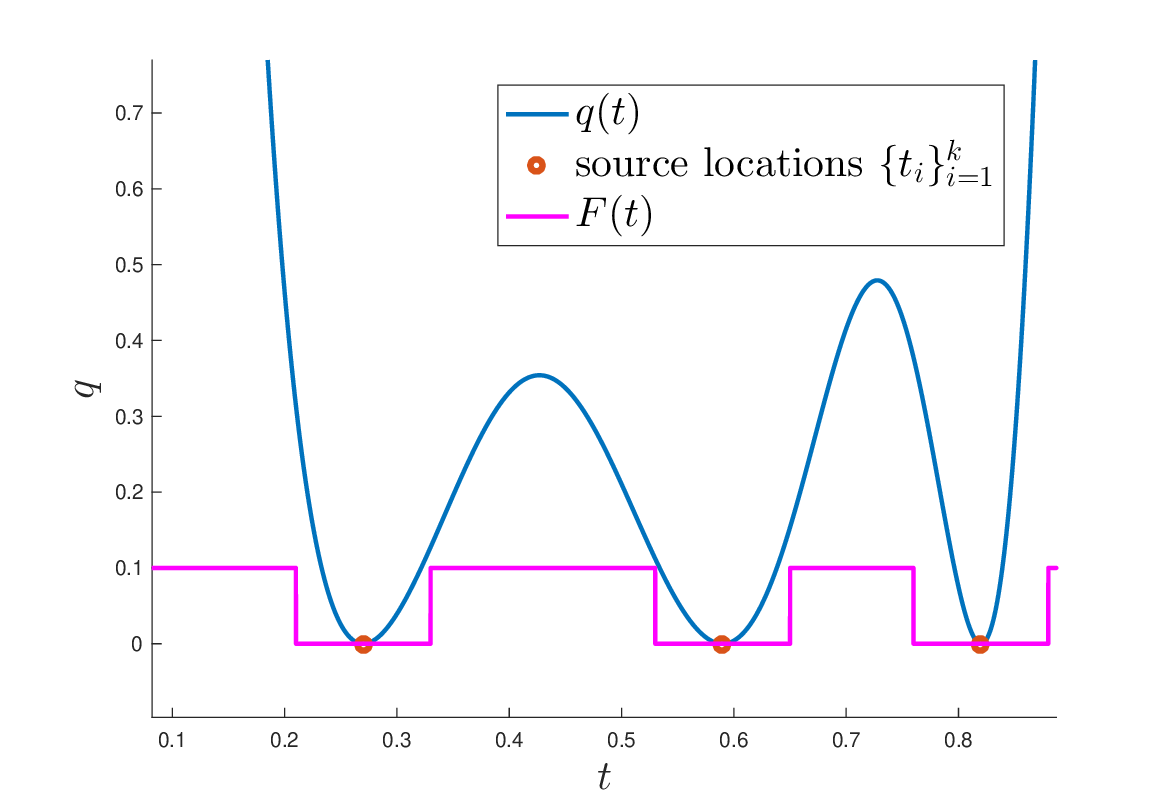}
    \caption{}
    \label{fig:example_dual_G_zoom}
  \end{subfigure}
  \caption{Example of dual certificate $q(t)$ that satisfies 
    the conditions
  in Lemma \ref{lem:dual} where the window function is the Gaussian
  kernel $\phi(t) = e^{-t^2 / \sigma^2}$. 
  We take $k=3$ and $ t_i \in \{0.27, 0.59,0.82\}$ and
  the function $F(t)$ such that $q(t) \geq F(t)$.}
  \label{fig:example_dual_G}
\end{figure}

To control $d_{GW}(\chi,\widehat{x})$ in turn, we will first show that
the existence of certain 
dual certificates leads to stability of Program (\ref{eq:mai}).
Then we see that these certificates exist under certain conditions on the measurement operator $\Phi$.  Turning now to the details, the following result is slightly more general than the one in \cite{candes2013super} and 
guarantees the stability of Program (\ref{eq:mai}) if a prescribed 
dual certificate $q$  exists. The proof is provided in Appendix \ref{sec:Proof-of-Lemma dual works}. 
\begin{lem}
  \label{lem:dual}\emph{\bf{(Error away from the support)}} Let $\widehat{x}$ be a solution of Program (\ref{eq:mai}) with $\delta'$ specified in (\ref{eq:def of deltap}) and set $h=\widehat{x}-\chi$ to be the error. 
  Consider $F(t)$ given in Definition \ref{def:dual_separators} and
  suppose that there exist a positive $\epsilon \le \Delta(T)/2$, 
  real coefficients $\{b_{j}\}_{j=1}^m$, and 
  a polynomial $q=\sum_{j=1}^m b_{j}\phi_{j}$ such that
  \[
    q(t)\ge F(t),
  \]
  where the equality holds on $T$. Then we have that 
  \begin{equation}
  \bar{f}\int_{T_{\epsilon}^{C}} h(\dif t)+\sum_{i=1}^k\int_{T_{i,\epsilon}}f\left(t-t_{i}\right) h(\dif t)\le2\|b\|_{2}\delta',\label{eq:o obar bound}
  \end{equation}
  where $b\in\mathbb{R}^{m}$ is the vector formed by the coefficients
  $\{b_{j}\}_{j=1}^m$. 
\end{lem}
There is a natural analogy here with the case of exact samples. In the
setting where $\eta_i=0$ in \eqref{eq:def_of_y}, the dual certificate
$q$ in Lemma \ref{lem:dual noise free} was required to be positive off
the support $T$. In the presence of inexact samples however, Lemma
\ref{lem:dual} loosely-speaking requires the dual certificate to be
bounded\footnote{Note the scale invariance of \eqref{eq:o obar bound} under scaling of $f$ and $\bar{f}$. Indeed, by changing $f,\bar{f}$ to $\alpha f,\alpha \bar{f}$ for positive $\alpha$, the proof dictates that $b$ changes to $\alpha b$ and consequently $\alpha$ cancels out from both sides of \eqref{eq:o obar bound}. Similarly, if we change $\Phi$ to $\alpha \Phi$ in \eqref{eq:Phi}, the proof dictates that  $b$ changes to $b/\alpha$ and $\alpha$ again cancels out, leaving \eqref{eq:o obar bound} unchanged.} \emph{away} from zero (see example in Figure
\ref{fig:example_dual_G}) for $t\in T_{\epsilon}^C$. 

Note also that Lemma \ref{lem:dual} controls the error $h$ away from the support $T$,
as it guarantees that 
\begin{equation}
  \int_{T_\epsilon^C}  h(\dif t) \le \frac{2\|b\|_2 \delta'}{\bar{f}},
  \label{eq:away spelled out}
\end{equation}
if the dual certificate $q$ exists. 
Indeed, \eqref{eq:away spelled out} follows directly from \eqref{eq:o obar bound} because the sum in \eqref{eq:o obar bound} is non-negative. This is in turn the case because $f(0)=0$ and the error $h$ is non-negative off the support $T$. 
Another key observation is that Lemma \ref{lem:dual} is  almost silent about the error near the impulses in $\chi$. Indeed,  because $f(0)=0$ by assumption, (\ref{eq:o obar bound}) completely fails to control the error on the support $T$. However, as the next result states, Lemma \ref{lem:dual} can be strengthened near the support provided that an additional dual certificate $q^0$ exists. The proof, given in Appendix \ref{sec:Proof-of-Lemma dual 3}, is not unsimilar to the analysis in \cite{fernandez2013support}.
\begin{lem}
  \label{lem:dual 3} \emph{\bf{(Error near the support)}} Suppose that the dual certificate in Lemma \ref{lem:dual}
  exists. Consider a function $F^0(t)$ defined 
  like $F^{\pi}(t)$ in Definition \ref{def:dual_separators} for 
  the sign pattern $\pi^0$ such that
  \[
    \pi^0_i = \begin{cases}
      1, & \text{when there exists }i\in[k]\text{ such that }
        t\in T_{i,{\epsilon}}\text{ and }
        \int_{T_{i,{\epsilon}}} h(\dif t)>0,\\
      -1, & \text{when there exists }i\in[k]\text{ such that }
        t\in T_{i,{\epsilon}}\text{ and }
        \int_{T_{i,{\epsilon}}} h(\dif t)\le0,
    \end{cases}
  \]
  and suppose also that there exist real coefficients 
  $\{b_{j}^{0}\}_{j\in [m]}$ and a polynomial 
  $q^{0}=\sum_{j=1}^m b_{j}^{0}\phi_{j}$
  such that 
  \[
    q^{0}(t)\ge F^{0}(t),
  \]
  where the equality holds on $T$. Then we have that 
  \begin{equation}
  \sum_{i=1}^k \left|\int_{T_{i,{\epsilon}}} h(\dif t)\right|\le2\left(\|b\|_{2}+\|b^0\|_{2}\right)\delta'.\label{eq:overall neighborhoods}
  \end{equation}
\end{lem}

\begin{figure}[h]
  \centering
  \begin{subfigure}{0.45\textwidth}
    \centering
    \includegraphics[width=\textwidth]{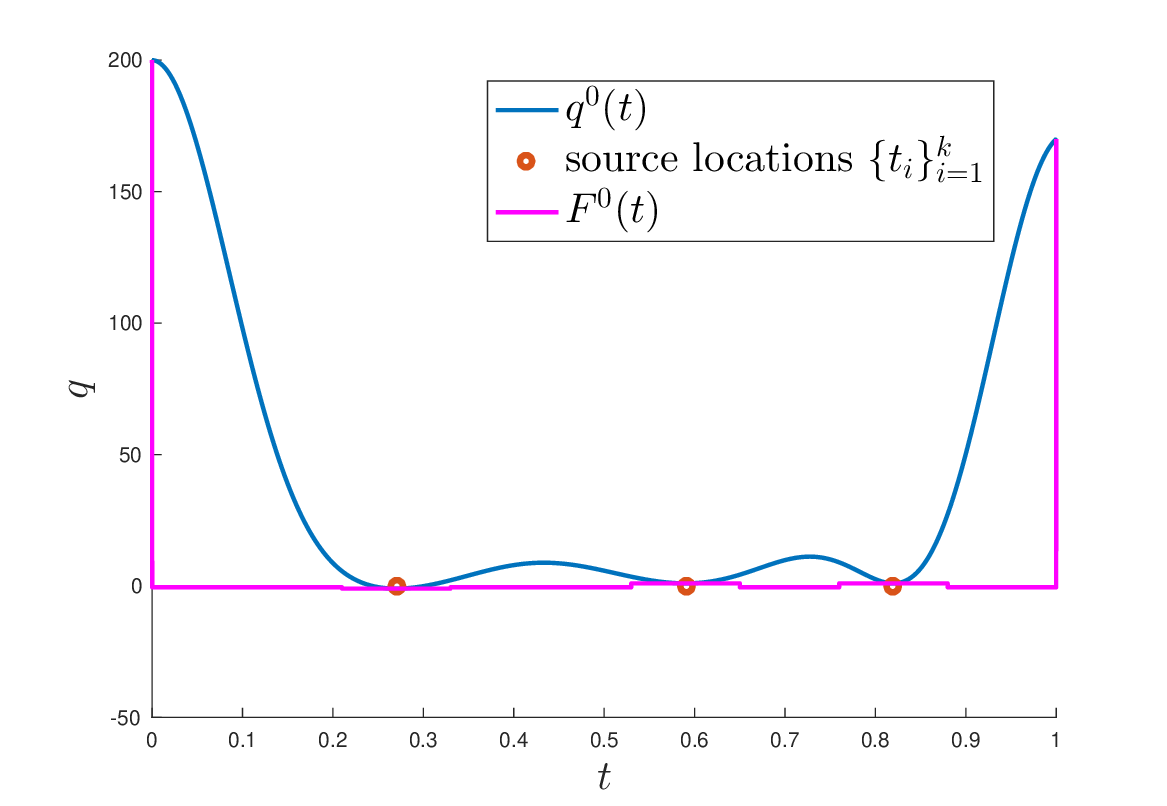}
    \caption{}
    \label{fig:example_dual_G_pi_big}
  \end{subfigure}
  \begin{subfigure}{0.45\textwidth}
    \centering
    \includegraphics[width=\textwidth]{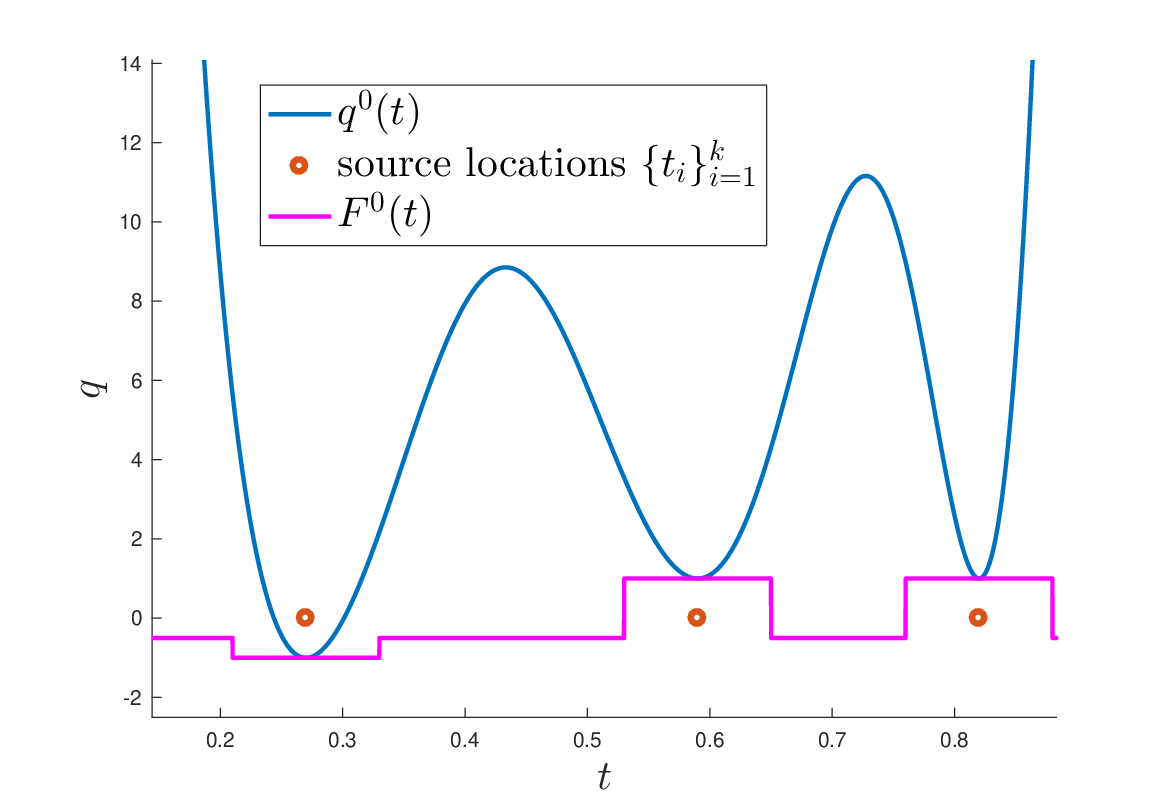}
    \caption{}
    \label{fig:example_dual_G_pi_zoom}
  \end{subfigure}
  \caption{Example of dual certificate $q^0(t)$ 
    that satisfies the conditions
  in Lemma \ref{lem:dual 3} where the window function is the Gaussian
  kernel $\phi(t) = e^{-t^2 / \sigma^2}$. 
  We take $k=3$ and $ t_i \in \{0.27, 0.59,0.82\}$ and
  the function $F^{0}(t)$ for the sign pattern $\pi^0 = \{-1, 1, 1\}$ 
  such that $q^0(t) \geq F^{0}(t)$.
  For the Gaussian kernel, the existence
  of $ q^{\pi}(t)$ for any sign pattern $\pi$ guarantees 
  the existence of $q^0(t)$ in Lemma \ref{lem:dual 3}.}
  \label{fig:example_dual_G_pi}
\end{figure}

In words, Lemma \ref{lem:dual 3} controls the error $ h$ near the support $T$, provided that a certain dual certificate $q^0$ exists (see example in Figure \ref{fig:example_dual_G_pi}). 
Note that \eqref{eq:overall neighborhoods} does not control the \emph{mass} of the error, namely $\sum_i \int_{T_{i,\epsilon}}   |h(\dif t)|= \int_{T_{\epsilon}}   |h(\dif t)|$, but rather it controls $\sum_i |\int_{T_{i,\epsilon}} h(\dif t)|$.  Of course, the latter is always bounded by the former, that is  
\begin{equation}
\label{eq:mass vs reality}
\sum_{i=1 }^k \l |\int_{T_{i,{\epsilon}}} h(\dif t)\r | \le  \int_{T_{{\epsilon}}}   |h(\dif t)|.
\end{equation}
However, the two sides of (\ref{eq:mass vs reality}) might differ significantly.  For instance, it might happen that the solution $\widehat{x}$ returns a slightly incorrect impulse at $t_i+\epsilon/4$ (rather than $t_i$) but with the correct amplitude of $a_i$. As  a result, the mass of the error is large in this case ($ \int_{T_{i,\epsilon}} |h(\dif t)| = 2a_i$) but the left-hand side of (\ref{eq:mass vs reality})  vanishes, namely $|\int_{T_{i,\epsilon}} h(\dif t)|=0$. 
Note that we cannot hope to strengthen \eqref{eq:overall neighborhoods} by replacing its left-hand side with the mass of the error, namely $\int_{T_\epsilon} |h(\dif t)|$.  This is the case mainly because the total variation is not the appropriate error metric for this context.

Indeed, while the mass of the error $\int_I |h(\dif t)|$ might not be small in general, we can instead control the generalised Wasserstein distance between the true and estimated measures, namely $x$ and $\widehat{x}$, see \eqref{eq:def of gen EMD} and \eqref{eq:def of EMD}.  The following result is proved by combining Lemmas \ref{lem:dual} and \ref{lem:dual 3}, see Appendix \ref{sec:Proof-of-Lemma EMD}.  
\begin{lem}
\label{lem:EMD} \emph{\bf{(Stability of Program \eqref{eq:mai} in the Generalised Wasserstein distance)}} Suppose that the dual certificates in Lemmas \ref{lem:dual}
and \ref{lem:dual 3} exist. Then it holds that 
\begin{equation}
d_{GW}\left( \chi,\widehat{x}\right)\le\left(\left(6+\frac{2}{\bar{f}}\right)\|b\|_{2}+6\|b^{0}\|_{2}\right)\delta'+\epsilon\| \chi \|_{TV}.
\label{eq:EMD}
\end{equation}
\end{lem}
An application of triangle inequality now yields that 
\begin{align}
d_{GW}(x,\widehat{x}) & \le d_{GW}(x,\chi) + d_{GW}(\chi,\widehat{x}) \nonumber\\
& \le d_{GW}(x,\chi) + \left(\left(6+\frac{2}{\bar{f}}\right)\|b\|_{2}+6\|b^{0}\|_{2}\right)\delta'+\epsilon\| \chi \|_{TV}.
\qquad \mbox{(see Lemma \ref{lem:EMD})}
\label{eq:final error before dual construction}
\end{align}
In words, Program \eqref{eq:mai} is stable $x$ if the certificates $q$ and $q^0$ exist. 
Let us now study the existence of these certificates. Proposition \ref{prop:dual construction}, proved in Appendix \ref{sec:Proof-of-Proposition dual construction},
guarantees the existence of the dual certificate $q$ required in Lemma
\ref{lem:dual} and heavily relies on the concept of T*-system in Definition \ref{def:(T-systems,-modified)-For}. We remark that the proof benefits from the ideas in \cite{karlin1966tchebycheff}. Similarly, Proposition
\ref{prop:dual construction 3},
stated without proof, ensures the existence of the certificate $q^0$ required in
Lemma \ref{lem:dual 3}.
\begin{prop}
  \label{prop:dual construction} \textbf{\emph{(Existence of $q$)}}  
  For $m\ge2k+2$, suppose that $\{\phi_{j}\}_{j=1}^m$
  form a T-system on $I$ and that $\{F\}\cup\{\phi_{j}\}_{j=1}^m$
  form a T*-system on $I$, where $F(t)$ is the function
  given in Definition \ref{def:dual_separators}.
  Then the dual certificate
  $q$ in Lemma \ref{lem:dual} exists and consequently 
  Program (\ref{eq:mai}) is stable in the sense 
  that (\ref{eq:o obar bound}) holds.
\end{prop}
Note that to ensure the success of Program \eqref{eq:mai}, 
it suffices  that there exists a polynomial $q=\sum_{j=1}^m b_j\phi_j$   
such that $q(t) \ge F(t)$ with the equality met on the support $T$, 
see Lemma \ref{lem:dual}. 
Equivalently, it suffices that there exists a non-negative 
polynomial $\dot{q}=-b_0 F + \sum_{j=1}^{m} b_j \phi_j $ 
that vanishes on $T$ such that $b_0>0$ and at least one 
other coefficient, say $b_{j_0}$, is nonzero. 
This situation is  reminiscent of Lemma \ref{lem:karlin-1}. 
In contrast to Lemma \ref{lem:karlin-1}, however, such $\dot{q}$ exists 
when $\{F\}\cup\{\phi_j\}_{j=1}^m$ is a T*-system rather than a T-system. The more subtle T*-system requirement is to avoid trivial or 
unbounded polynomials.  

\begin{prop}
  \label{prop:dual construction 3} \textbf{\emph{(Existence of $q^0$)}} 
  For $m\ge2k+2$, suppose that $\{\phi_{j}\}_{j=1}^{m}$
  form a T-system on $I$ and that $\{F^{0}\}\cup\{\phi_{j}\}_{j=1}^{m}$
  form a T*-system on $I$, where $F^0(t)$ is the function 
  defined in Lemma \ref{lem:dual 3}.
  Then the dual certificate
  $q^0$ in Lemma \ref{lem:dual 3} exists and consequently 
  Program (\ref{eq:mai}) is stable in the sense 
  that (\ref{eq:overall neighborhoods}) holds.
\end{prop}
Having constructed the necessary dual certificates in Propositions \ref{prop:dual construction} and \ref{prop:dual construction 3}, the proof of Theorem \ref{thm:main noisy} is now complete in light of  \eqref{eq:final error before dual construction}.

\subsection{Proof of Theorem \ref{thm:noisy_expl} (Average stability for Program \eqref{eq:mai})}

In this section we give an overview of the main ideas involved in proving Theorem \ref{thm:noisy_expl}.
To start with, let $A \in \mathbb{R}^{k \times k} $ be defined 
as in \eqref{eq:matrix_A}:
      \begin{equation*}
		    A = 
		    \begin{bmatrix}
		      |\phi_1(t_1)|  & -|\phi_1(t_2)| & \ldots & -|\phi_1(t_k)| \\
		      -|\phi_2(t_1)| & |\phi_2(t_2)|  & \ldots & -|\phi_2(t_k)| \\
		      \vdots         & \vdots         & \ddots & \vdots         \\
		      -|\phi_k(t_1)| & -|\phi_k(t_2)| & \ldots & |\phi_k(t_k)|  \\
		    \end{bmatrix},
		  \end{equation*}
where $\phi_i(t_i) = \phi(t_i-s_{l(i)})$ is evaluated at the source $t_i$ and 
the closest sample to it, as defined in \eqref{eq:closest_sample_thm}.

The proof of Theorem \ref{thm:noisy_expl} consists of two steps. We first show that we 
can bound the error if the matrix $A$ is strictly diagonally dominant.
It is easy to see that, if the window function $\phi$ is localised, then the entries on the
main diagonal are larger in absolute value than the off-diagonal entries. If, 
moreover, we choose the sampling locations $\{s_i\}_{i \in [m]}$ such that $A$ 
is strictly diagonally dominant (which means that for each source, there is a 
sampling location that is "close enough" to it), then the bound \eqref{eq:bound_near_support}
is guaranteed.

\begin{prop}
  \label{prop:A_diagonally_dominant}
For each source $t_i$, select $s_{l(i)}$ to be the closest sample as defined 
in \eqref{eq:closest_sample_thm}, 
and define the matrix A in \eqref{eq:matrix_A} using
the sequences $\{t_i\}_{i=1}^{k}$, $\{s_{l(i)}\}_{i=1}^{k}$. If $A$ is strictly diagonally
dominant, then the error around the support is bounded according to \eqref{eq:bound_near_support}.
\end{prop}

Then, we want to go further and see what it means exactly for $A$ to be strictly
diagonally dominant, so the second step in the proof of Theorem \ref{thm:noisy_expl} 
is to give an upper bound for the distance between the sources $\{t_i\}_{i \in [k]}$ 
and the closest sampling locations $\{s_{l(i)}\}_{i \in [k]}$ 
such that $A$ is strictly diagonally dominant.

Given an even positive function $\phi$ that is localised at $0$ and with fast decay, 
let $\Delta$ and $\lambda$ as given in Definition \ref{def:min_sep}, so 
\begin{equation}
  |t_i - s_{l(i)}| \le \lambda \Delta 
  \label{sitibound}
\end{equation}
We want to find $\lambda_0$ such that 
\begin{equation}
  \label{eq:cond_A_dd}
  \phi(s_{l(i)} - t_i) \geq 
  \sum_{j \ne i} \phi(s_{l(i)} - t_j),
  \quad \forall \lambda \in (0, \lambda_0),
  \quad \forall i \in [k],
\end{equation}
namely, we want the matrix $A$ to be strictly diagonally dominant. 
From the conditions \eqref{eq:cond_A_dd}, we can obtain a more
general equality, depending on $\phi$ and $\Delta$, that $\lambda_0$ 
must satisfy such that, for $\lambda < \lambda_0$, $A$ is 
strictly diagonally dominant. The equality is given by \eqref{eq:condition_lambda}:
  \begin{equation*}
    \phi(\lambda_0 \Delta)
    =
    \phi(\Delta - \lambda_0 \Delta) +
    \phi(\Delta + \lambda_0 \Delta) +
    \frac{1}{\Delta}
      \int_{\Delta - \lambda_0\Delta}^{1/2 - \lambda_0\Delta} \phi(x) \dif x +
    \frac{1}{\Delta}
      \int_{\Delta + \lambda_0\Delta}^{1/2 + \lambda_0\Delta} \phi(x) \dif x.
  \end{equation*}

\begin{prop}
  \label{prop:sampling_locations}
  Let $\lambda_0 \in (0,\frac12) $ such that $\left|t_i - s_{l(i)}\right| \leq \lambda_0 \Delta$ 
  for all $i \in [k]$. If $\lambda_0$ satisfies \eqref{eq:condition_lambda},
  then the matrix A defined in \eqref{eq:matrix_A} is strictly diagonally dominant.
\end{prop}

Finally, we note that the proof of Theorem \ref{thm:grouping_sources} involves the same ideas 
as the ones discussed in this section, with a few modifications.
The detailed proofs of Proposition \ref{prop:A_diagonally_dominant} and
Proposition \ref{prop:sampling_locations} are given in 
Appendices \ref{sec:proof_A_diagonally_dominant} and \ref{sec:sampling_locations} respectively. 
The proof of Theorem \ref{thm:grouping_sources} is similar to the proof presented in the current section, so
we only show the differences in Appendix \ref{sec:proof_grouping_sources}.

\subsection{Proofs of Theorems \ref{thm:main Gaussian} and \ref{thm:main Gaussian lambda}
 (Gaussian with sparse measure)}

In this section we give the main steps taken to obtain the explicit 
bounds in Theorems \ref{thm:main Gaussian} and \ref{thm:main Gaussian lambda}
for the Gaussian window function.
These are particular cases of the more general Theorems \ref{thm:main noisy} 
and \ref{thm:noisy_expl} respectively, where the window function is 
taken to be the Gaussian 
$\phi_j(t) = e^{-(t - s_j)^2 / \sigma^2}$,
given in \eqref{eq:gaussian window}, and the true measure $x$
is a k-discrete non-negative measure as in \eqref{eq:def of x}.

\subsubsection{Bounds on the coefficients of the dual certificates for Gaussian window}
\label{sec:bounds coeffs}

We will now give explicit bounds on the vectors of 
coefficients $\|b\|_2$ and $\|b^{\pi}\|_2$ of the
dual certificates $q$ and $q^{\pi}$ 
from Lemmas \ref{lem:dual} and \ref{lem:dual 3}
in terms of the parameters of the problem $k, T, S$ 
and $\sigma$ (the width of the Gaussian window).

Firstly, we introduce a more specific form of the 
dual polynomial separators $F(t)$, $F^{\pi}(t)$ 
from Definition \ref{def:dual_separators}. 
Here, we take $f(t)=0$ for $t \in (-\epsilon, \epsilon)$
and $\bar{f}, f_1$ positive constants
with $\bar{f} < 1$. Then, $f_0$ is defined to be greater that 
both $\bar{f}$ and $f_1$, with the exact relationship between $f_0$ 
and $\bar{f}$ given in the proof of 
Lemma \ref{lem:Gaussian is Tstar sys}.
Therefore, for $\epsilon > 0$ and a sign 
pattern $\pi\in\{\pm 1\}^k$, $F(t)$ and
$F^{\pi}(t)$ are:
\begin{align}
  F(t)&=\begin{cases}
    {f}_{0}, &  t=0,\\
    {f}_{1}, & t = 1,\\
    0, & \text{when there exists }i\in[k]\text{ such that }t\in {T}_{i,{\epsilon}},\\
    \bar{f}, & \text{elsewhere on } I.
  \end{cases}
  \label{eq:f example}
  \\
  F^{\pi}(t)&=\begin{cases}
    f_{0}, & t=0,\\
    f_{1}, & t=1,\\
    \pi_i, & \text{when there exists }i\in[k]
      \text{ such that }t\in {T}_{i,\epsilon},
      \\
    -\bar{f}, & \mbox{elsewhere on } I.
  \end{cases}
  \label{eq:g pi}
\end{align}
With the above definitions, both Theorem \ref{thm:main noisy}
and Theorem \ref{thm:noisy_expl} require that 
$\{F\}\cup\{\phi_j\}_{j=1}^m$ form a T*-system on $I$. 
Likewise, Theorem \ref{thm:main noisy} requires that
$\{F^\pi\}\cup\{\phi_j\}_{j=1}^m$ form a T*-system for 
any sign pattern $\pi$. We show 
in Lemma \ref{lem:Gaussian is Tstar sys}
that both these requirements
are satisfied for the choice in \eqref{eq:gaussian window}
of $\phi_j(t) = e^{-(t - s_j)^2 / \sigma^2}$. 
The proof is given in 
Appendix \ref{sec:proof of lemma Gaussian is Tstar sys}. 
\begin{lem}\label{lem:Gaussian is Tstar sys}
  Consider the function $F(t)$ defined in \eqref{eq:f example} 
  and suppose that  $m \geq 2k+2$. 
  Then $\{F\}\cup\{\phi_j\}_{j=1}^m$ form a T*-system on $I$, 
  with $\phi$ extended totally positive, 
  even Gaussian and $\phi_j$ defined as in \eqref{eq:gaussian window}, 
  provided that $f_{0} \gg \bar{f}$, $f_{0} \gg f_{1}$ 
  and $\bar{f}, f_{0},f_{1} \gg 0$.
  These requirements are made precise in the proof 
  and are dependent on $\epsilon$.  Moreover, for an 
  arbitrary sign pattern $\pi$ and $F^\pi$ as defined 
  in \eqref{eq:g pi}, $\{F^\pi\}\cup \{\phi_j\}_{j=1}^m$ form  
  a T*-system on $I$ when, in addition, $f_0 \gg 1$.
\end{lem}

In this setting, consider a subset of $m=2k+2$ 
samples $\{s_j\}_{j=1}^{m} \subset S$
(since in the proof of Lemma \ref{lem:Gaussian is Tstar sys} we select
the $2k+2$ samples that are the closest to the sources) such that
they satisfy Conditions \ref{cond:conds thms 1 2}.
Therefore, we have 
that $s_1 = 0$, $s_m = s_{2k+2}=1$, and
\begin{equation}
  | s_{2i} -t_i | \le \eta, \qquad s_{2i+1}-s_{2i} =\eta,
  \qquad \forall i\in [k],
  \label{eq:assump display}
\end{equation} 
for a small $\eta \leq \sigma^2$, see \eqref{eq:gaussian window}.
That is, we collect two samples close to each impulse $t_i$ in $x$, 
one on each side of $t_i$. 

Suppose also that 
\begin{equation}
  \sigma \leq \sqrt{2}, \quad 
  \Delta > \sigma \sqrt{\log\left(3 + \frac{4}{\sigma^2} \right)}, 
  \quad \eta \leq \sigma^2,
  \label{eq:assump sigma sep}
\end{equation}
namely the width of the Gaussian is much smaller than the separation of 
the impulses in $x$. Lastly, assume that the impulse 
locations $T=\{t_i\}_{i=1}^k$ and  sampling 
points $\{s_j\}_{j=2}^{m-1}$ are away from the boundary of $I$, namely 
\begin{equation*}
  \sigma\sqrt{ \log(1/\eta^3)} \le t_i \le 1-\sigma\sqrt{\log(1/\eta^3)},
  \qquad \forall i\in [k],
\end{equation*}
\begin{equation}
  \sigma\sqrt{ \log(1/\eta^3)} \le s_j \le 1-\sigma\sqrt{ \log(1/\eta^3)},
  \qquad \forall j\in [2:m-1].
  \label{eq:away frm bndr display}
\end{equation}

\textit{Remark 1.}
The Property
\textit{2. Samples near sources} in Conditions \ref{cond:conds thms 1 2}
states that for each two samples $s'$ and $s$
near each source, we have that:
\begin{equation}
  C_1 \eta \leq s'-s \leq C_2 \eta,
  \label{eq:C1etaC2eta}
\end{equation}
but in \eqref{eq:assump display} above we simplified the condition
to $s_{2i+1}-s_{2i}=\eta$. Throughout the proofs in this paper we will
use the simplified condition in \eqref{eq:assump display} instead 
of the more general \eqref{eq:C1etaC2eta} so that we do not obscure
the central issues of the proof with extra indices and separate treatment
of the upper and lower bounds for $s'-s$. This has implications for 
Lemma \ref{lem:bounds on b for Gaussian} below, and we will
point out the places in its proof where using the more general
condition \eqref{eq:C1etaC2eta} would require separate treatment
(see footnotes \ref{fn:C1etaC2_1} and \ref{fn:C1etaC2_2}).

\textit{Remark 2.}
The Property \textit{3. Sources away from the boundary}
in Conditions \ref{cond:conds thms 1 2} is necessary due to our method
of proof, which imposes that the sources and samples are in 
the interval $I=[0,1]$ independently from the width $\sigma$ of the 
convolution kernel and the minimum separation of sources $\Delta$. 
While the exact form of the boundary conditions depend on the specific
approach that we took in the proof, 
there is an interplay between $\sigma$ and $\Delta$ in the noisy setting
and due to the fact that the interval $I=[0,1]$ is fixed, some form of scaling
is required, which is what this condition achieves. For example, 
for $\sigma=1$, the Gaussian kernel defined in \eqref{eq:conv} 
varies by at most $\frac{1}{e}$ (between two sources located
at $t_1=0$ and $t_2=1$) and therefore our results for the noisy 
setting are not meaningful for large $\sigma$.

We can now give explicit bounds on $\|b\|_2$
and $\|b^{\pi}\|_2$ for the Gaussian window function,
as required by Theorems \ref{thm:main noisy} 
and \ref{thm:noisy_expl}. The following result is proved 
in Appendix \ref{sec:proof of lemma bounds on b for Gaussian}.

\begin{lem}\label{lem:bounds on b for Gaussian}
  Suppose that the window function $\phi$ is Gaussian, as 
  defined in \eqref{eq:gaussian window}, 
  the assumptions \eqref{eq:assump display}, 
  \eqref{eq:assump sigma sep} and \eqref{eq:away frm bndr display} 
  (namely Conditions \ref{cond:conds thms 1 2}) hold 
  and $\eta$ satisfies:
  \begin{align}
      \eta \leq \min \left\{
      \frac{
        8 F_{\min}(\Delta,\frac{1}{\sigma})
      }{
        34(2k+2)
        \left(
          80k + 8 + k P\left(\frac{1}{\sigma}\right)
          \frac{2}{
            1 - e^{-\frac{\Delta^2}{\sigma^2}}
          }
        \right)^{\frac12}
      }
      ,
      \frac{
        \bar{C}(f_0,f_1)^{\frac16}
      }{
        \left(4k+4 + \frac{4k}{\sigma^2}\right)^{\frac13}
      }
      \right\}.
      \label{eq:main cond eta} 
  \end{align}
  Then we have the following bounds:
  \begin{align}
    \|b\|_2  &\le 
    \frac{
    \sqrt{
      (2k+2)
      \left(
        4k + 5 + \frac{4k}{\sigma^4}
      \right)
    }
    }{
      1 - \frac{\sqrt{e}}{2}
    } 
    \bar{C}(f_0,f_1)^{\frac54}
    \left[
      \frac{
        F_{\max}\left(\Delta,\frac{1}{\sigma}\right)
      }{
        F_{\min}\left(\Delta,\frac{1}{\sigma}\right)^2
      }
    \right]^k,
    \label{eq:bound b}
    \\
    \|b^{\pi}\|_2 & \leq
    \frac{
      \sqrt{
      2k+2
      }
    }{
      \eta \left(1 - \frac{\sqrt{e}}{2}\right)
    }
    \left(\bar{C}(f_0,f_1) + 2k\right)^{\frac32}
    \left[
      \frac{
        F_{\max}\left(\Delta,\frac{1}{\sigma}\right)
      }{
        F_{\min}\left(\Delta,\frac{1}{\sigma}\right)^2
      }
    \right]^k,
    \label{eq:bound b pi}
  \end{align}
  where 
  \begin{align}
    &\bar{C}(f_0,f_1) = f_0^2 + f_1^2 + 2f_0 + 2f_1 + 2,
    \label{eq:c bar def}
    \\
    &P\left(\frac{1}{\sigma}\right) = 
    \frac{4}{\sigma^4} 
    + \frac{13}{4} \left(
        \frac{2}{\sigma^2} + \frac{4}{\sigma^4}
      \right)^2
    + \frac{9}{4} \left(
        \frac{12}{\sigma^4} + \frac{8}{\sigma^6}
      \right)^2.
    \label{eq:P def}
    \\
    &F_{\max}\left(\Delta,\frac{1}{\sigma}\right) =
      \left(
        8 + \left(1 + \frac{4}{\sigma^4} \right)
          \frac{2}{1-e^{-\frac{\Delta^2}{\sigma^2}}}
      \right)^{\frac12}
      \left(
        32 + \left( 
          \frac{1}{\sigma^4} + \frac{2}{\sigma^6} + \frac{2}{\sigma^8}  
        \right)
        \frac{16}{1-e^{-\frac{\Delta^2}{\sigma^2}}}
      \right)^{\frac12},
    \label{eq:Fmax def}\\
    &F_{\min}\left(\Delta,\frac{1}{\sigma}\right) =  
        1 - \left(1 + \frac{2}{\sigma^2} \right)
        \frac{
          2 e^{-\frac{\Delta^2}{\sigma^2}}
        }{
          1 - e^{-\frac{\Delta^2}{\sigma^2}}
        }.
    \label{eq:Fmin def}
  \end{align}
\end{lem}
To obtain the final bounds for the Gaussian window function,
we will substitute the above bounds in the right hand side
of \eqref{eq:EMD thm} in Theorem \ref{thm:main noisy} 
and in \eqref{eq:bound_near_support} and \eqref{eq:local avg away} 
in Theorem \ref{thm:noisy_expl}. We will then 
obtain $F_1$ in Theorem \ref{thm:main Gaussian} 
(see \eqref{eq:final F1}) and $F_2$
in Theorem \ref{thm:main Gaussian lambda} 
(see \eqref{eq:final F2}).

For more clarity, in the following lemma
we simplify $F_1$ and $F_2$ further, in the case 
when stronger conditions apply to $\sigma$, $\Delta(T)$
and $\lambda$.
\begin{lem}
  \label{lem:in particular}
  If the conditions in Lemma \ref{lem:bounds on b for Gaussian}
  hold and, in addition, 
  $\sigma < \frac{1}{\sqrt{3}}$,
  $\Delta > \sigma \sqrt{\log{\frac{5}{\sigma^2}}}$
  and $\bar{f}<1$, then 
  \begin{align}
    &\frac{
      F_{\max}\left(\Delta,\frac{1}{\sigma}\right)
    }{
      F_{\min}\left(\Delta,\frac{1}{\sigma}\right)^2   
    }
    < 
    \frac{c_1}{\sigma^6 (1-3\sigma^2)^2},
    \label{eq:in particular 0} 
    \\
    &\left(
      (6 + \frac{2}{\bar{f}})
      \sqrt{4k + 5 + \frac{4k}{\sigma^4}}
      \bar{C}(f_0,f_1)^\frac54
      +
      \frac{6}{\eta}(\bar{C}(f_0,f_1)+2k)^\frac32
    \right)
    \frac{\sqrt{2k+2}}{1-\frac{\sqrt{e}}{2}}
    <
    c_2 \cdot 
    \frac{k C_1(\frac{1}{\epsilon})}{\eta\sigma^2},
    \label{eq:in particular 1}
    \\
    &\frac{
      \sqrt{
        (2k+2)
        \left(4k + 5 + \frac{4k}{\sigma^4}\right)
      }
    }{
      1 - \frac{\sqrt{e}}{2}
    }
    \frac{
      \bar{C}(f_0,f_1)^\frac54}{
      \bar{f}
    }
    <
    c_3 \cdot \frac{
      k
      C_2(\frac{1}{\epsilon})
    }{
      \sigma^2
    },
    \label{eq:in particular 12} 
  \end{align}
  where 
  \begin{align}
    C_1\left(\frac{1}{\epsilon}\right) &= 
    \frac{(\bar{C}(f_0,f_1)+2k)^\frac32}{\bar{f}},
    \label{eq:c1 eps def}
    \\
    C_2\left(\frac{1}{\epsilon}\right) &= 
    \frac{\bar{C}(f_0,f_1)^\frac54}{\bar{f}},
    \label{eq:c2 eps def}
  \end{align}
  and, similarly, the condition \eqref{eq:main cond eta} is 
  simplified to condition \eqref{eq:cond eta partic}.
  Moreover, if $\lambda$ in Theorem \ref{thm:noisy_expl} 
  satisfies $\lambda < \frac24$, then
  \begin{equation}
    \frac{1}{
      e^{-\frac{\Delta^2 \lambda^2}{\sigma^2}}
      -e^{-\frac{\Delta^2 \lambda^2}{\sigma^2}}
      \cdot \frac{
        e^{-\frac{\Delta^2}{\sigma^2}}
        + e^{-\frac{2 \Delta^2}{\sigma^2}}
      }{
        1 - e^{-\frac{\Delta^2}{\sigma^2}}
      }
      -e^{-\frac{\Delta^2 (1 - \lambda)^2}{\sigma^2}}
    }
    <
    c_4.
    \label{eq:in particular 2}
  \end{equation}
  Above, $c_1, c_2, c_3, c_4$ are universal constants.
\end{lem}
More specifically, \eqref{eq:in particular 0}, 
\eqref{eq:in particular 1} 
will be used to bound $F_1$ to 
obtain \eqref{eq:gaussian in particular}
and \eqref{eq:in particular 0}, \eqref{eq:in particular 12}, 
will be used to bound $F_2$ to 
obtain \eqref{eq:in particular f2}.
Lastly, \eqref{eq:in particular 2} will be used to 
bound $F_3$, which appears in the bound 
given by Theorem \ref{thm:main Gaussian lambda}.

The proof of Lemma \ref{lem:in particular} is given 
in Appendix \ref{sec:in particular}. Note that we 
give $C_1$ and $C_2$ as functions of $\frac{1}{\epsilon}$ because,
as $\epsilon \to 0$, $C_1$ and $C_2$ grow at a rate dependent
on $\epsilon$, as indicated in the Lemma \ref{lem:small epsilon}
in Section \ref{sec:proof of corr}.

\subsubsection{Proof of Theorem \ref{thm:main Gaussian} 
  (Wasserstein stability of Program
      \eqref{eq:mai} for $\phi(t)$ Gaussian)
  \label{sec:proof of main Gaussian}}

We consider Theorem \ref{thm:main noisy} and restrict ourselves to the 
case $x = \chi$, a k-discrete non-negative measure as 
in \eqref{eq:def of x} with support $T$.
Let $\Delta = \Delta(T) \ge 2\epsilon$, with $T$ as the support of $x$. 
We begin with estimating the Lipschitz constant of the measurement operator $\Phi$ according to  \eqref{eq:Lipschitz assumption}, see Appendix \ref{sec:proof of L cte} for the proof of the next result.
\begin{lem}\label{lem:L cte}
Consider $S=\{s_j\}_{j=1}^m\subset \mathbb{R}$ and $\{\phi_j\}_{j=1}^m$ specified in \eqref{eq:gaussian window}. Then the operator $\Phi:I \rightarrow \mathbb{R}^m$ defined in (\ref{eq:Phi}) 
is $\frac{2\sqrt{m}}{\sigma \sqrt{2e}}$-Lipschitz with respect to the generalised 
Wasserstein distance, namely (\ref{eq:Lipschitz assumption}) holds 
with $L=\frac{2\sqrt{m}}{\sigma \sqrt{2e}}$. 
\end{lem} 

We may now invoke Theorem \ref{thm:main noisy} to conclude that, for an
arbitrary $k$-sparse $2\epsilon$-separated 
non-negative measure $x$ and arbitrary sampling points 
$\{s_j\}_{j=1}^m \subset\mathbb{R}$, Program (\ref{eq:mai}) 
with $\delta'=\delta$ is stable  in the sense that there exist 
vectors $b,\{b^\pi\}_\pi \subset\mathbb{R}^m$ such that 
\begin{equation}
  d_{GW}\left( x,\widehat{x}\right)
  \le  
  \left(  
    \left(
      6+\frac{2}{\bar{f}}
    \right)\|b\|_{2}
    +6\min_\pi \|b^{\pi}\|_{2}
  \right) \cdot \delta
  + \|x\|_{TV} \cdot \epsilon
\label{eq:bnd on err Gaussians 1}
\end{equation}
provided that $m\ge 2k+2$ and $f_0\gg f_1 \gg \bar{f} \gg 0$. 
The exact relationships between $\bar{f}, f_0, f_1$ are given
in the proof of Lemma \ref{lem:Gaussian is Tstar sys}
in Appendix \ref{sec:proof of lemma Gaussian is Tstar sys}.

Combining Lemma \ref{lem:bounds on b for Gaussian} with \eqref{eq:bnd on err Gaussians 1} yields a final bound on the stability of Program \eqref{eq:mai} with Gaussian window and completes the proof of the first part 
of Theorem \ref{thm:main Gaussian}:
  \begin{equation}
    d_{GW}(x,\widehat{x}) 
    \le
    F_1(k,\Delta(T),\frac{1}{\sigma},\frac{1}{\epsilon},\eta)
    \cdot \delta
    + \|x\|_{TV} \cdot \epsilon,
  \end{equation}
  where
  \begin{align}
    &F_1(k,\Delta(T),\frac{1}{\sigma},\frac{1}{\epsilon},\eta)=
      \left(
        (6 + \frac{2}{\bar{f}})
        \sqrt{4k + 5 + \frac{4k}{\sigma^4}}
        \bar{C}^\frac54
        +
        \frac{6}{\eta}(\bar{C}+2k)^\frac32
      \right)
      \frac{\sqrt{2k+2}}{1-\frac{\sqrt{e}}{2}}
      \left(
      \frac{
        F_{\max}\left(\Delta,\frac{1}{\sigma}\right)
      }{
        F_{\min}\left(\Delta,\frac{1}{\sigma}\right)^2
      }
      \right)^k, 
      \label{eq:final F1}
  \end{align}
with $\bar{C},F_{\max},F_{\min}$ given 
in \eqref{eq:c bar def}, \eqref{eq:Fmax def},
\eqref{eq:Fmin def} respectively
and $f_0=f_0(\frac{1}{\epsilon})$ depends on $\epsilon$ 
(see the proof of Lemma \ref{lem:small epsilon}).

Finally, to show that \eqref{eq:gaussian in particular} holds 
when $\sigma < \frac{1}{\sqrt{3}}$,
and $\Delta > \sigma\sqrt{\log{\frac{5}{\sigma^2}}}$, 
we apply the first part of Lemma \ref{lem:in particular} 
with $\bar{f} < 1$
and the proof of Theorem \ref{thm:main Gaussian} is complete.

\textit{Remark.}
 In particular, $f_0$ increases as $\epsilon \to 0$ and $f_0$ also 
depends on the other parameters of the problem, 
namely $\eta,\sigma,\Delta(T),k$.
See Section \ref{sec:discussion} for a detailed discussion.
Furthermore, $f_1$ and $\bar{f}$ are considered fixed positive 
constants with $f_1 < f_0$.

\subsubsection{Proof of Theorem \ref{thm:main Gaussian lambda} 
  (Average stability of Program \eqref{eq:mai} for $\phi(t)$
    Gaussian: source proximity dependence)
  \label{sec:proof of main Gaussian lambda}
}

We now apply Theorem \ref{thm:noisy_expl} 
with $\phi(t) = g(t) = e^{-t^2 / \sigma^2}$.
We have that $\phi^{\infty} = 1$, the Lipschitz constant $L$ of $g$ 
on $[-1,1]$ is $L=\frac{2}{\sigma^2}$, and
\begin{equation}
  \sum_{j=1}^{k}(A^{-1})_{ij} \leq \| A^{-1} \|_{\infty}
  <
  \frac{1}{\min_{j}
    \left(
      g(s_j - t_j) - \sum_{i \ne j} g(s_j - t_i)
    \right)
  }.
  \label{eq:infty norm bound}
\end{equation}
The last inequality comes from the definition of $A$
in \eqref{eq:matrix_A} with $\phi(t) = g(t)$ and
$s_j := s_{l(j)}$ as given in Definition \ref{def:min_sep},
and \cite{varah1975}.
Then, by assumption, $|s_j - t_j| \leq \lambda\Delta$ for an
arbitrary $j \in [k]$ and $g$ is decreasing, so
\begin{equation*}
  g(s_j - t_j) \geq g(\lambda\Delta) 
  = e^{-\frac{\lambda^2\Delta^2}{\sigma^2}}.
\end{equation*}
We now assume without loss of generality that $s_j < t_j$. Then, it follows that
\begin{equation*}
  |s_j - t_i| \geq |j-i|\Delta - \lambda\Delta, 
    \text{ if } i < j
  \quad\quad \text{ and } \quad \quad
  |s_j - t_i| \geq |j-i|\Delta + \lambda\Delta,
    \text{ if } i > j.
\end{equation*}
This, in turn, leads to
  \begin{align}
    \sum_{i \ne j} g(s_j - t_i) &= 
      \sum_{i=1}^{j-1} g(s_j - t_i) +
      \sum_{i=j+1}^{k} g(s_j - t_i) \nonumber \\
    &\leq
      \sum_{i=1}^{j-1} g((j-i)\Delta - \lambda\Delta) +
      \sum_{i=j+1}^{k} g((i-j)\Delta + \lambda\Delta) \nonumber \\
    &\leq
      \sum_{i=1}^{\infty} g((i-\lambda)\Delta) +
      \sum_{i=1}^{\infty} g((i+\lambda)\Delta).
      \label{eq:sum two terms}
  \end{align}
We now bound each sum in \eqref{eq:sum two terms} as follows
  \begin{align}
    \sum_{i=1}^{\infty} g((i-\lambda)\Delta) &=
    \sum_{i=1}^{\infty} e^{-\frac{(i-\lambda)^2\Delta^2}{\sigma^2}} 
    \nonumber \\
    &= 
      e^{-\frac{\Delta^2(1-\lambda)^2}{\sigma^2}} + 
      e^{-\frac{\Delta^2\lambda^2}{\sigma^2}}
      \sum_{i=2}^{\infty} 
        \left(e^{-\frac{\Delta^2}{\sigma^2}}\right)^{i^2 - 2i\lambda} 
    \nonumber \\
    &\leq
      e^{-\frac{\Delta^2(1-\lambda)^2}{\sigma^2}} + 
      e^{-\frac{\Delta^2\lambda^2}{\sigma^2}}
      \sum_{i=2}^{\infty}
        \left(e^{-\frac{\Delta^2}{\sigma^2}}\right)^{i} 
        \quad\quad\quad (i^2 - 2i\lambda > i \text{ for } i \geq 2)
    \nonumber \\
    &=
      e^{-\frac{\Delta^2(1-\lambda)^2}{\sigma^2}} + 
      e^{-\frac{\Delta^2\lambda^2}{\sigma^2}} \cdot
      \frac{e^{-\frac{2\Delta^2}{\sigma^2}}}{
        1 - e^{-\frac{\Delta^2}{\sigma^2}} 
      },
    \label{eq:first sum}
  \end{align}
and similarly, we have that
\begin{equation}
  \sum_{i=1}^{\infty} g((i+\lambda)\Delta) =
  \sum_{i=1}^{\infty} e^{-\frac{(i+\lambda)^2\Delta^2}{\sigma^2}} 
  \leq
  e^{-\frac{\Delta^2\lambda^2}{\sigma^2}} \cdot
  \frac{e^{-\frac{\Delta^2}{\sigma^2}}}{
    1 - e^{-\frac{\Delta^2}{\sigma^2}} 
  }.
  \label{eq:second sum}
\end{equation}
By combining \eqref{eq:sum two terms}, \eqref{eq:first sum} and \eqref{eq:second sum},
we obtain:
\begin{equation}
  g(s_j - t_j) - \sum_{i \neq j} g(s_j - t_i) 
  \geq
  e^{-\frac{\Delta^2\lambda^2}{\sigma^2}} -
  e^{-\frac{\Delta^2\lambda^2}{\sigma^2}} \cdot
  \frac{
     e^{-\frac{\Delta^2}{\sigma^2}} 
     + e^{-\frac{2\Delta^2}{\sigma^2}} 
  }{
    1 - e^{-\frac{\Delta^2}{\sigma^2}} 
  }
  -e^{-\frac{\Delta^2(1-\lambda)^2}{\sigma^2}}.
\end{equation}
The above inequality also holds when we take the minimum 
over $j \in [k]$ and, inserting it in \eqref{eq:infty norm bound}
and using this result and the bound on $\|b\|_2$ from 
Lemma \ref{lem:bounds on b for Gaussian} 
in \eqref{eq:bound_near_support},
we obtain \eqref{eq:main Gaussian lambda res}:
\begin{equation}
    \left|
      \int_{T_{i,\epsilon}} \hat{x}(\dif t) - a_i
    \right|
    \leq
    \left[
      (c_1 + F_2)
      \cdot \delta
      + c_2 \frac{\|\hat{x}\|_{TV} }{\sigma^2}
      \cdot \epsilon
    \right] F_3,
    \label{eq:f2 final}
\end{equation}
where 
\begin{align}
    &F_2(k,\Delta(T),\frac{1}{\sigma},\frac{1}{\epsilon})=
      \frac{
        \sqrt{
          (2k+2)
          \left(4k + 5 + \frac{4k}{\sigma^4}\right)
        }
      }{
        1 - \frac{\sqrt{e}}{2}
      }
      \frac{
        \bar{C}(f_0,f_1)^\frac54}{
        \bar{f}
      }
      \left[\frac{
        F_{\max}\left(\Delta,\frac{1}{\sigma}\right)
      }{
        F_{\min}\left(\Delta,\frac{1}{\sigma}\right)^2
      }
      \right]^k,
    \label{eq:final F2}
    \\
    &F_3(\Delta(T),\sigma,\lambda) =
    \frac{1}{
      e^{-\frac{\Delta^2 \lambda^2}{\sigma^2}}
      -e^{-\frac{\Delta^2 \lambda^2}{\sigma^2}}
      \cdot \frac{
        e^{-\frac{\Delta^2}{\sigma^2}}
        + e^{-\frac{2 \Delta^2}{\sigma^2}}
      }{
        1 - e^{-\frac{\Delta^2}{\sigma^2}}
      }
      -e^{-\frac{\Delta^2 (1 - \lambda)^2}{\sigma^2}}
    },
    \label{eq:final F3}
\end{align}
and $\bar{C},F_{\max},F_{\min}$ are given 
in \eqref{eq:c bar def}, \eqref{eq:Fmax def},
\eqref{eq:Fmin def} respectively. The error
bound away from the sources \eqref{eq:main Gaussian lambda res away}
is obtained by applying Lemma \ref{lem:dual}
with the same bounds on $\|b\|_2$.

Then, by using Lemma \ref{lem:in particular} with $\bar{f} < 1$, 
we obtain \eqref{eq:in particular f2}.
Note that we can apply Lemma \ref{lem:in particular} because,
for $\sigma < \frac{1}{\sqrt{3}}$, we have that $\frac{5}{\sigma^2} > 3 + \frac{4}{\sigma^2}$
and, therefore, 
$\Delta > \sigma\sqrt{\log{\frac{5}{\sigma^2}}} > \sigma\sqrt{\log{(3+\frac{4}{\sigma^2})}}$.

\subsubsection{Proof of Corollary \ref{corr:small eps}}
\label{sec:proof of corr}

First we give an explicit dependence of $C_1(\frac{1}{\epsilon})$
and $C_2(\frac{1}{\epsilon})$ on $\epsilon$ for small $\epsilon>0$ 
in the following lemma, proved in Appendix \ref{sec:proof of small epsilon}.
\begin{lem}
  \label{lem:small epsilon}
  If $f_1 < f_0$, $1 < f_0$ and $\bar{f} < 1$, 
  then there exists $\epsilon_0 > 0$ such that:
  \begin{align}
    C_1\left(\frac{1}{\epsilon}\right) 
    &< \frac{(\bar{c}_1 C_{\epsilon}^2 + 2k \epsilon^4)^\frac32}{\bar{f}} \cdot \frac{1}{\epsilon^6}
    \\
    C_2\left(\frac{1}{\epsilon}\right) 
    &<
    \bar{c}_2 C_{\epsilon}^\frac52 
    \cdot \frac{1}{\epsilon^5},
  \end{align}
  for all $\epsilon \in (0,\epsilon_0)$,
  where $\bar{c}_1$ and $\bar{c}_2$ are universal constants
  and $C_{\epsilon}$ is defined in the proof, 
  see \eqref{eq:f0 fbar eps}.
\end{lem}

To prove the first part of the corollary, we first 
let $\epsilon = \delta^{\frac17}$ in the bound 
on $C_1(\frac{1}{\epsilon})$ in Lemma \ref{lem:small epsilon}:
\begin{equation}
  C_1\left(\frac{1}{\epsilon}\right) 
  < \frac{(\bar{c}_1 C_{\epsilon}^2 + 2k \epsilon^4)^\frac32}{\bar{f}} \cdot \frac{1}{\delta^{\frac67}},
  \quad \forall \delta < \epsilon_0^7,
\end{equation}
and we substitute the above inequality in the 
bound \eqref{eq:gaussian in particular} in Theorem \ref{thm:main Gaussian} to obtain:
\begin{equation*}
  d_{GW}(x,\hat{x}) < 
  \bar{C}_1 \cdot \delta^{\frac17},
  \quad \forall \delta < \epsilon_0^7,
\end{equation*}
where
\begin{equation}
  \bar{C}_1 = 
  \frac{
    c_1 k 
    (\bar{c}_1 C_{\epsilon}^2 + 2k \epsilon^4)^\frac32 
  }{
    \eta \sigma^2 \bar{f}
  }
  \left[
    \frac{c_2}{\sigma^6 (1-3\sigma^2)^2}
  \right]^k
  +
  \|x\|_{TV}.
\end{equation}
Similarly, let $\epsilon = \delta^{\frac16}$ in the 
bound on $C_2(\frac{1}{\epsilon})$ in Lemma \ref{lem:small epsilon}:
\begin{equation}
  C_2\left(\frac{1}{\epsilon}\right) 
  <
  \bar{c}_2 C_{\epsilon}^\frac52 
  \cdot \frac{1}{\delta^{\frac56}},
  \quad \forall \delta < \epsilon_0^6,
\end{equation}
which we substitute in the bound \eqref{eq:in particular f2}
in Theorem \ref{thm:main Gaussian lambda} to obtain:
\begin{equation*}
  \left| \int_{T_{i,\epsilon}} \hat{x}(\dif t) - a_i \right|
    \le \bar{C}_2 \cdot \delta^{\frac16},
  \quad \forall \delta < \epsilon_0^6,
\end{equation*}
where 
\begin{equation}
  \bar{C}_2 =
  c_1
  + c_3 \frac{
      k \bar{c}_2 C_{\epsilon}^\frac52
    }{
      \sigma^2
    } 
    \left[
      \frac{c_4}{\sigma^6 (1-3\sigma^2)^2}
    \right]^k
  + \frac{c_2 \|\hat{x}\|_{TV}}{\sigma^2}.
\end{equation}
Note that we apply Lemma \ref{lem:small epsilon} 
with $\bar{f} < 1$ and that both inequalities in the 
corollary hold for $\delta < \delta_0 = \epsilon_0^7$,
where $\epsilon_0$ is given by Lemma \ref{lem:small epsilon}.

\subsection{Discussion \label{sec:discussion}}

In this section, we discuss a few issues regarding the robustness of our
construction of the dual certificate from Appendix \ref{sec:Proof-of-Proposition dual construction}.
There are two points that
need to be raised: the construction itself and the proof that we indeed have
a T*-System.

At the moment, we do not use any samples that are away from sources in the 
the construction of the dual certificate. 
If the sources are close enough compared to $\sigma$, then this is not an issue. 
However, for $\sigma$ small relative to the distance between samples,
in light of the proof of Lemma \ref{lem:Gaussian is Tstar sys} 
(see Appendix \ref{sec:proof of lemma Gaussian is Tstar sys}),
if we consider the dual certificate as 
the expansion of the determinant $N$ of $M^{\rho}$ 
in \eqref{eq:m rho}
along the $\tau_{\underline{l}}$ row:
\begin{equation}
  N = -F(\tau_{l}) \beta_0 + 
    \sum_{j=1}^{m} (-1)^{j+1} \beta_j g(\tau_l - s_j),
\end{equation}
then the terms $g(\tau_{\underline{l}})$ become exponentially small 
(as $\tau_{\underline{l}}$ is far from all samples $s_j$)
and, therefore, the value of $N$ is close to $-F(\tau_l)$ (which is $-\bar{f}$ 
if $\tau_l \in T_{\epsilon}^C$). This is problematic, as we require 
that $N > 0$.
We can overcome this by adding ``fake'' sources at intervals $\eta^{-1}$ 
so that they cover the regions where we have no true sources, together with
two close samples for each extra source. The determinant $N$ becomes:
\begin{equation}
  N = \begin{vmatrix}
    f_0 & f(s_1) & \cdots & g(s_m) \\
    0 & g(t_1 - s_1) & \cdots & g(t_1 - s_m) \\
    0 & g'(t_1 - s_1) & \cdots & g'(t_1 - s_m) \\
    \vdots & \vdots & & \vdots \\
    \bar{f} & g(\tau_j - s_1) & \cdots & g(\tau_j - s_m) \\
    \bar{f} & g'(\tau_j - s_1) & \cdots & g'(\tau_j - s_m) \\
    \vdots & \vdots & & \vdots \\
    F(\tau_{\underline{l}}) & g(\tau_{\underline{l}} - s_1) & \cdots & g(\tau_{\underline{l}} - s_m) \\
    \vdots & \vdots & \cdots & \vdots \\
    0 & g(t_k - s_1) & \cdots & g(t_k - s_m) \\
    0 & g'(t_k - s_1) & \cdots & g'(t_k - s_m) \\
    f_1 & g(1-s_1) & \cdots & g(1-s_m)   
  \end{vmatrix}.
  \label{eq:new construction}
\end{equation}
Here, the rows are ordered according to the ordering of the set 
containing $t_i, \tau_j, \tau_{\underline{l}}$. The terms in the expansion 
of \eqref{eq:new construction} along the row
with $\tau_{\underline{l}}$ do not approach $0$ exponentially with this construction,
since for any $\tau_{\underline{l}}$ there exists $s_i$ close enough so 
that $g(\tau_{\underline{l}} - s_i) > f^*$ for some $f^* > 0$.

More specifically, consider also the expansion of $N$ along the first column:
\begin{equation}
  N = f_0 N_{1,1} + f_1 N_{m+1,1} 
    - F(\tau_{\underline{l}}) N_{\tau_{\underline{l}},1}
    - \bar{f} \sum_{j < \tau_{\underline{l}}} (N_{j,1} - N_{j+1,1})
    + \bar{f} \sum_{j > \tau_{\underline{l}}} (N_{j,1} - N_{j+1,1}).
  \label{eq:new constr expansion}
\end{equation}
We use this expansion in the proof of Lemma \ref{lem:Gaussian is Tstar sys} in Appendix \ref{sec:proof of lemma Gaussian is Tstar sys} to show that the 
functions $F \cup \{g_j\}_{j=1}^{m}$ form a T*-System.
For $\tau_{\underline{l}} \in T_{\epsilon}^C$, $F(\tau_{\underline{l}})=\bar{f}$
and the setup in Lemma \ref{lem:Gaussian is Tstar sys}, 
we require that (see \eqref{eq:dominated}):
\begin{equation}
  \frac{f_0}{\bar{f}} \geq 
    \frac{N_{\tau_{\underline{l},1}}}{
      \min_{\tau_{\underline{l}} \in T_{\epsilon}^C} N_{1,1}
    }. 
  \label{eq:fs condition}
\end{equation} 
In the construction \eqref{eq:new construction}, if we upper bound the 
pairs $N_{j,1} - N_{j+1,1}$ in the two sums in \eqref{eq:new constr expansion}
(a separate problem by itself),
then we can impose a similar condition to \eqref{eq:fs condition} 
for $f_0$ and $\bar{f}$. From here, we obtain that
  $f_0 = C \bar{f}$
where finding
$
  C \geq 
    \frac{N_{\tau_{\underline{l},1}}}{
      \min_{\tau_{\underline{l}} \in T_{\epsilon}^C} N_{1,1}
    }
$
involves finding a lower bound on $N_{1,1}$:
\begin{equation}
  N = \begin{vmatrix}
    g(t_1 - s_1) & \cdots & g(t_1 - s_m) \\
    g'(t_1 - s_1) & \cdots & g'(t_1 - s_m) \\
    \vdots & & \vdots \\
    g(\tau_j - s_1) & \cdots & g(\tau_j - s_m) \\
    g'(\tau_j - s_1) & \cdots & g'(\tau_j - s_m) \\
    \vdots & & \vdots \\
    g(\tau_{\underline{l}} - s_1) & \cdots & g(\tau_{\underline{l}} - s_m) \\
    \vdots & & \vdots \\
    g(t_k - s_1) & \cdots & g(t_k - s_m) \\
    g'(t_k - s_1) & \cdots & g'(t_k - s_m) \\
    g(1-s_1) & \cdots & g(1-s_m)   
  \end{vmatrix}.
  \label{eq:new construction N11}
\end{equation}
The structure of the above determinant is similar to the denominator in 
Appendix \ref{sec:proof of lemma bounds on b for Gaussian}
but only up to the row with $\tau_{\underline{l}}$. The rows after it do not
preserve the diagonally dominant structure of the matrix, as each source becomes
associated with one close sample to it and the first sample corresponding 
to the next source. This is an issue in both the construction
described in the proof of Proposition \ref{prop:dual construction} 
in Appendix \ref{sec:Proof-of-Proposition dual construction} 
(and detailed in the proof of Lemma \ref{lem:Gaussian is Tstar sys})
and the construction described in the current section 
(which would result from considering a determinant 
with ``fake'' sources like \eqref{eq:new construction}). However, 
by adding extra ``fake'' sources, one could argue that 
the determinant \eqref{eq:new construction} is better behaved,
as the distance between a source and the first sample corresponding to the next 
source is smaller, which we leave for further work.

\section*{Acknowledgements}

AE is grateful to Gongguo Tang and Tamir Bendory for enlightening discussions 
and helpful comments. 
This work was done while HT was affiliated to the Alan Turing Institute, London, 
and the School of Mathematics, University of Edinburgh, UK.
This publication is based on work supported by 
the EPSRC Centre For Doctoral Training in Industrially Focused Mathematical 
Modelling (EP/L015803/1) in collaboration with the National Physical Laboratory 
and by the Alan Turing Institute under the EPSRC grant EP/N510129/1 
and the Turing Seed Funding grant SF019. The authors thank Stephane Chretien, 
Alistair Forbes and Peter Harris from the National Physical Laboratory
for the support and insight regarding the associated applications of super-resolution.
Lastly, the authors thank the reviewers for their suggestions to improve the
generality of the results.

\bibliographystyle{unsrt}
\bibliographystyle{plain}
\bibliography{References}

\appendix

\section{Proof of Lemma \ref{lem:dual noise free} \label{sec:Proof-of-Lemma noise free dual cert}}

Let $\hat{x}$ be a solution of Program \ref{eq:mai} with $\delta=0$
and let $h=\hat{x}-x$ be the error. Then, by feasibility of both
$x$ and $\widehat{x}$ in Program (\ref{eq:mai}), we have that 
\begin{equation}
\int_I \phi_{j}(t)h(\dif t)=0,\qquad j\in[m].\label{eq:feasibility-1}
\end{equation}
Let $T^{C}$ be the complement of $T=\{t_{i}\}_{i=1}^k$ with respect
to $I$. By assumption, the existence of a dual certificate allows
us to write that 
\begin{align*}
\int_{T^{C}}q(t) h(\dif t) & =\int_{I} q(t)h(\dif t)-\int_{T}q(t)h(\dif t)\\
& = \int_{I} q(t)h(\dif t) \qquad \l( q\l(t_i\r) = 0,\,\, i\in[k]\r) \\ 
 & =\sum_{j=1 }^m b_{j}\int_{I} \phi_{j}(t)h(\dif t)\\
 & =0. \qquad \mbox{(see \eqref{eq:feasibility-1})}
\end{align*}
Since $x=0$ on $T^C$, then $h=\hat{x}$ on $T^C$, so the last equality is
equivalent to
\begin{equation}
  \int_{T^C} q(t) \hat{x}(\dif t) = 0.
\end{equation}
But $q$ is strictly positive on $T^C$, so it must be 
that $h = \hat{x} = 0$ on $T^C$ and, therefore, 
$h=\sum_{i=1}^k c_{i}\delta_{t_{i}}$
for some coefficients $\{c_{i}\}$. Now (\ref{eq:feasibility-1})
reads $\sum_{i=1}^k c_{i}\phi_{j}(t_{i})=0$ for every $j\in[m]$.
This gives $c_{i}=0$ for all $i\in[k]$ because $[\phi_{j}(t_{i})]_{i,j}$
is, by assumption, full rank. Therefore $h=0$  and $\hat{x}=x$ on $I$,
which completes the proof of Lemma \ref{lem:dual noise free}.

\section{Proof of Lemma \ref{lem:dual}\label{sec:Proof-of-Lemma dual works}}

Let $\hat{x}$ be a solution of Program \eqref{eq:mai} and set $h=\hat{x}-\chi$ to  be the error. Then,
by feasibility of both $\chi$ and $\widehat{x}$ in Program (\ref{eq:mai})
and using the triangle inequality, we have that 
\begin{equation}
\left\Vert \int_I  \Phi(t) h(\dif t)\right\Vert _{2}\le2\delta'.\label{eq:feasibility}
\end{equation}
Next, the existence of the dual certificate $q$ allows us to write
that 
\begin{align*}
& \bar{f}\int_{T_{\epsilon}^{C}} h(\dif t)+\sum_{i=1}^k \int_{T_{i,\epsilon}}f\left(t-t_{i}\right) h(\dif t) \\
& \le\int_{T_{\epsilon}^{C}}q(t) h(\dif t)+\sum_{i=1}^k \int_{T_{i,\epsilon}}q(t) h(\dif t)\\
 & =\int_{T_{\epsilon}^{C}}q(t) h(\dif t)+\int_{T_{\epsilon}}q(t) h(\dif t)
\qquad \l( T_\epsilon = \cup_{i=1}^k T_{i,\epsilon} \r) 
 \\
 & =\int_{I}q(t) h(\dif t)\\
 & =\sum_{j=1}^m b_{j}\int_I \phi_{j}(t) h(\dif t)\\
 & \le\|b\|_{2}\cdot\left\Vert \int_{I} \Phi(t) h(\dif t)\right\Vert _{2}
\qquad \mbox{(Cauchy-Schwarz inequality)}
 \\
 & \le\|b\|_{2}\cdot2\delta',
 \qquad \mbox{(see \eqref{eq:feasibility})}
\end{align*}
which completes the proof of Lemma \ref{lem:dual}. 

\section{Proof of Lemma \ref{lem:dual 3}\label{sec:Proof-of-Lemma dual 3}}

The existence of the dual certificate $q^0$ allows us to write
that 
\begin{align*}
 & \sum_{i=1}^k \left|\int_{T_{i,\epsilon}} h(\dif t)\right|\\
 & =\sum_{i=1}^k \int_{T_{i,\epsilon}}s_{i} h(\dif t)\qquad s_{i}=\text{sign}\left(\int_{T_{i,\epsilon}} h(\dif t)\right)\\
 & =\sum_{i=1}^k \int_{T_{i,\epsilon}}\left(s_{i}-q^{0}(t)\right) h(\dif t)+\sum_{i =1}^k \int_{T_{i,\epsilon}}q^{0}(t) h(\dif t)\\
 & =\sum_{i=1}^k \int_{T_{i,\epsilon}}\left(s_{i}-q^{0}(t)\right) h(\dif t)+\int_I q^{0}(t) h(\dif t)-\int_{T_{\epsilon}^{C}}q^{0}(t) h(\dif t)\\
 & =\sum_{s_{i}=1}\int_{T_{i,\epsilon}}\left(1-q^{0}(t)\right) h(\dif t)+\sum_{s_{i}=-1}\int_{T_{i,\epsilon}}\left(-1-q^{0}(t)\right) h(\dif t)\\
 & \qquad+\int_I q^{0}(t) h(\dif t)-\int_{T_{\epsilon}^{C}}q^{0}(t) h(\dif t)\\
 & \le\sum_{s_{i}=1}\int_{T_{i,\epsilon}}f\left(t-t_{i}\right) h(\dif t)+\sum_{s_{i}=-1}\int_{T_{i,\epsilon}}f\left(t-t_{i}\right) h(\dif t)+\int_I q^{0}(t) h(\dif t)+\bar{f}\int_{T_{\epsilon}^{C}} h(\dif t)\\
 & =\sum_{i=1}^k \int_{T_{i,\epsilon}}f\left(t-t_{i}\right) h(\dif t)+\bar{f}\int_{T_{\epsilon}^{C}} h(\dif t)+\int_I q^{0}(t) h(\dif t)\\
 & \le2\|b\|_{2}\delta'+\int_I q^{0}(t) h(\dif t)\qquad\text{(see Lemma \ref{lem:dual})}\\
 & =2\|b\|_{2}\delta'+\sum_{j=1}^m b_{i}^{0}\int_I \phi_{j}(t) h(\dif t)\\
 & \le2\|b\|_{2}\delta'+\|b^{0}\|_{2}\cdot2\delta',
 \qquad \mbox{(see \eqref{eq:feasibility})}
\end{align*}
which completes the proof of Lemma \ref{lem:dual 3}.

\section{Proof of Lemma \ref{lem:EMD} \label{sec:Proof-of-Lemma EMD}}

Our strategy is as follows. We first argue that 
\begin{equation*}
d_{GW} \l( \chi, \widehat{x} \r) \approx d_{GW} \l(  \chi,  \widetilde{x}\r),
\end{equation*}
where $\widetilde{x}$ is the restriction of $\widehat{x}$ to the $\epsilon$-neighbourhood of the support of $\chi$, namely  $T_{\epsilon}$ defined in \eqref{eq:T_epsilon}.  This,  loosely speaking, reduces the problem to that of  computing the distance between two discrete measures supported on $T$. We control the latter distance using a particular suboptimal choice of measure $\gamma$ in \eqref{eq:def of EMD}. Let us turn to the details.  

Let $\widetilde{x}$ be the restriction
of $\widehat{x}$ to $T_{\epsilon}$, namely $\widetilde{x}=\widehat{x}|_{T_{\epsilon}}$, or more specifically
\begin{equation*}
\widetilde{x}(\dif t) = 
\begin{cases}
\widehat{x}(\dif t) & t\in T_\epsilon,\\
0 & t\in T_\epsilon^C.
\end{cases} 
\end{equation*}
 Then,
using the triangle inequality, we observe that 
\begin{equation}
d_{GW}\left( \chi,\widehat{x}\right)\le d_{GW}\left( \chi,\widetilde{x}\right)+d_{GW}\left(\widetilde{x},\widehat{x}\right).\label{eq:brk}
\end{equation}
The last distance above is easy to control: We write that 
\begin{align}
  d_{GW}\left(\widetilde{x},\widehat{x}\right) 
  & =\inf_{\tilde{z},\hat{z}}
    \left( 
      \l\|\widetilde{x}-\widetilde{z}\r\|_{TV}
      +\l\|\widehat{x}-\widehat{z}\r\|_{TV}
      +d_{W}\left(\widetilde{z},\widehat{z}\right)
    \right),
\qquad \mbox{(see \eqref{eq:def of gen EMD})}
\nonumber \\
 & \le \l\|\widetilde{x}-\widehat{x}\r\|_{TV}+\l\|\widehat{x}-\widehat{x}\r\|_{TV}+d_{W}\left(\widehat{x},\widehat{x}\right)\qquad\left(\widetilde{z}=\widehat{z}=\widehat{x}\right)\nonumber \\
 & =\l\|\widetilde{x}-\widehat{x}\r\|_{TV}\nonumber \\
 & =\l\|\widehat{x}|_{T_{\epsilon}^{C}}\r\|_{TV}\nonumber \\
 & =\int_{T_{\epsilon}^{C}} \widehat{x}(\dif t)
\qquad \l( \widehat{x} \mbox{ is non-negative} \r) 
 \nonumber \\
 & =\int_{T_{\epsilon}^{C}} h(\dif t)
\qquad \l( h(\dif t) = \widehat{x}(\dif t)-\chi(\dif t) = \widehat{x}(\dif t) \mbox{ when } t\in T^C  \r) 
 \nonumber \\
 & \le\frac{2\|b\|_{2}\delta'}{\bar{f}}.
\qquad \mbox{(see Lemma \ref{lem:dual})} 
 \label{eq:leg 1}
\end{align}
We next control the
term $d_{GW}\left( \chi,\widetilde{x}\right)$ in (\ref{eq:brk})
by writing that 
\begin{align}
  d_{GW}\left( \chi,\widetilde{x}\right) 
  & =\inf_{\tilde{z}, \hat{z}}
    \left(
      \l\| \chi-z\r\|_{TV}
      +\l\|\widetilde{x}-\widetilde{z}\r\|_{TV}
      +d_{W}\left(z,\widetilde{z}\right)
    \right)
\qquad \mbox{(see \eqref{eq:def of gen EMD})}
\nonumber \\
 & \le\l\| \chi- \chi\r\|_{TV}+\l\|\widetilde{x}-\widetilde{x}\cdot \frac{\l\| \chi\r\|_{TV}}{\l\|\widetilde{x}\r\|_{TV}}\r\|_{TV}+d_{W}\left(\chi,\widetilde{z}\right)\qquad\left(z= \chi,\,\widetilde{z}=\widetilde{x}\cdot \frac{\l\| \chi\r\|_{TV}}{\l\|\widetilde{x}\r\|_{TV}}\right)\nonumber \\
 & =\l| \l\|\widetilde{x}\r\|_{TV}-\l\| \chi\r\|_{TV} \r|+d_{W}\left(\chi,\widetilde{z}\right)\nonumber \\
 & =\l| \l\|\widetilde{x}|_{T_{\epsilon}^{C}}\r\|_{TV}+\l\|\widetilde{x}|_{T_{\epsilon}}\r\|_{TV}-\l\| \chi|_{T_{\epsilon}}\r\|_{TV} \r| +d_{W}\left(\chi,\widetilde{z}\right)\nonumber \\
 & =\l| \l\|\widehat{x}|_{T_{\epsilon}}\r\|_{TV}-\l\| \chi |_{T_{\epsilon}}\r\|_{TV} \r| +d_{W}\left(\chi,\widetilde{z}\right)
\qquad
 \left(\widetilde{x}=\widehat{x}|_{T_{\epsilon}}\right)
 \nonumber \\
 & =\l| \int_{T_{\epsilon}} \widehat{x}(\dif t)-\int_{T_{\epsilon}} \chi(\dif t) \r| +d_{W}\left(\chi,\widetilde{z}\right)
 \qquad
\left( \chi \mbox{ and } \widehat{x} \mbox{ are non-negative} \right)  
\nonumber \\
 & =\l| \int_{T_{\epsilon}} h(\dif t) \r| +d_{W}\left(\chi,\widetilde{z}\right)\nonumber \\
 & \le\sum_{i=1}^k \left|\int_{T_{i,\epsilon}} h(\dif t)\right|+d_{W}\left(\chi,\widetilde{z}\right)
\qquad \mbox{(triangle inequality)} 
 \nonumber \\
 & \le2\left(\|b\|_{2}+\|b^0\|_{2}\right)\delta'+d_{W}\left(\chi,\widetilde{z}\right).
\qquad \mbox{(see Lemma \ref{lem:dual 3})} 
 \label{eq:leg 2 pre}
\end{align}
For future use, we record an intermediate result that is obvious from studying \eqref{eq:leg 2 pre}, namely  
\begin{equation}
\l | \l\| \wt{x} \r\|_{TV}-\|\chi\|_{TV} \r| \le 2\l( \|b\|_2+\l\| b^0\r\|_2 \r) \delta'. 
\label{eq:intermediate result}
\end{equation}
It remains to control
$d_{W}(\chi,\widetilde{z})$ above where
\begin{equation}
d_W\l(\chi,\widetilde{z}\r)  = \inf \int_{I^2} \l |\tau-\widetilde{\tau} \r| \gamma\l(\dif \tau,\dif \widetilde{\tau} \r) ,
\label{eq:EMD repeated}
\end{equation}
is the Wasserstein distance between the measures $\chi$ and $\widetilde{z}$. 
The infimum above is over all  measures $\gamma$ on $I^2 = I\times I$ that produce $\chi$ and $\widetilde{z}$ as marginals, namely we have: 
\begin{equation}
  \int_{A \times I}\gamma(\dif \tau, \dif \widetilde{\tau}) = \chi(A)
  \quad \text{ and } \quad
  \int_{I \times B}\gamma(\dif \tau, \dif \widetilde{\tau}) = \widetilde{z}(B),
  \quad \forall A,B \subset I
  \label{eq:marginals_def}
\end{equation}
For every $i\in [k]$, let also $\widetilde{x}_{i}=\widetilde{x}|_{T_{i,\epsilon}}$ and $\widetilde{z}_{i}=\widetilde{z}|_{T_{i,\epsilon}}$ 
be the restrictions of $\widetilde{x}$ and  $\widetilde{z}$ to $T_{i,\epsilon}$, respectively. Because of our choice of $z,\widetilde{z}$ in the second line of \eqref{eq:leg 2 pre}, note that $\widetilde{z}_i =  \widetilde{x}_i \cdot \|z\|_{TV}/\| \widetilde{x}\|_{TV}$.  
Recalling that $\chi$ is supported on $T=\{t_i\}_{i=1}^k$, we write that $\chi = \sum_{i=1}^k a_i \delta_{t_i}$ for non-negative amplitudes $\{a_i\}_{i=1}^k$.  
Then, noting that $\chi$ is supported on $T$ and $\tilde{z}$ is supported on $T_{\epsilon}$, 
then any feasible $\gamma$ in \eqref{eq:EMD repeated} is supported on $T \times T_{\epsilon}$ and we can construct a feasible but 
suboptimal $\gamma(\dif \tau, \dif \widetilde{\tau})$ in a two-step
approach, where we first extract up to $a_i$ weight of $\widetilde{z}_i$
on each $\delta_{t_i}$, for example let
\begin{equation}
  \gamma_1 = \sum_{i=1}^{k} \widetilde{z}_i(\dif \widetilde{\tau})
  1_{T_i}
  \quad \text{where} \quad
  1_{T_i} = 
  \begin{cases}
    T_{i,\epsilon} \quad &\text{if} \quad
      \int_{T_{i,\epsilon}} z_i (\dif \tau) \le a_i, \\
    [t_i - \xi_i, t_i + \xi_i] \quad &\text{if otherwise},
  \end{cases}
\end{equation}
where $\xi_i$ is defined such 
that $\int_{t_i-\xi_i}^{t_i+\xi_i} z_i(\dif \tau) = a_i$.
As a result, $\gamma_1$ has $\dif\tau$ marginal equal 
to $\widetilde{z}_i$ on the support of $\gamma_1$ and 
the $\dif \widetilde{\tau}$ marginal no more than the desired $a_i$.
We then construct $\gamma_2$ by partitioning the 
remaining $\widetilde{z}_i$ into the $\dif \widetilde{\tau}$
subsets in order to make up the $\tilde{z}$ marginal, which is 
exactly achievable using all of $\tilde{z}$ 
due to $\int \tilde{z} (\dif \tau) = \sum_{i=1}^{k} a_i$. Then,
we take $\gamma = \gamma_1 + \gamma_2$.

Intuitively, this is a transport plan according to which we move as much mass as possible 
inside each $T_{i,\epsilon}$ (from $\delta_{a_i}$ to $\wt{z_i}$, the minimum of the masses
the two) and the remaining mass
is moved outside $T_{i,\epsilon}$.
Therefore, for this choice of $\gamma$, we have that
\begin{align}
  d_W(\chi,\widetilde{z}) 
  &\leq
  \int_{I^2} 
    |\tau - \widetilde{\tau}| \gamma (\dif \tau, \dif \widetilde{\tau}) 
  \nonumber \\
  &=
  \sum_{i=1}^{k} \int_{\{t_i\} \times T_{i, \epsilon}} 
    |\tau - \widetilde{\tau}| \gamma (\dif \tau, \dif \widetilde{\tau}) 
  + \sum_{i=1}^{k} \int_{\{t_i\} \times T_{i, \epsilon}^C} 
    |\tau - \widetilde{\tau}| \gamma (\dif \tau, \dif \widetilde{\tau}) 
    \nonumber \\
  &\leq
  \epsilon \sum_{i=1}^{k} \int_{\{t_i\} \times T_{i, \epsilon}} 
    \gamma_1 (\dif \tau, \dif \widetilde{\tau})
  + \sum_{i=1}^{k} \int_{\{t_i\} \times T_{i, \epsilon}^C} 
    \gamma_2 (\dif \tau, \dif \widetilde{\tau}) 
\end{align}
The third line above uses the fact that $ \tau = t_i $ and,
if $\widetilde{\tau}\in T_{i,\epsilon}$, then $|\tau - \widetilde{\tau}|\le \epsilon$
and $|\tau - \widetilde{\tau}|\le 1$ otherwise. 
By evaluating the integrals in the last line above, we find that 
\begin{align}
  d_W\l(\chi, \wt{z} \r)
  &\le \epsilon \sum_{i=1}^{k} \min\{a_i, \|\wt{z_i}\|_{TV}\}
    + \sum_{i=1}^{k} \left( a_i - \min\{a_i, \|\wt{z_i}\|_{TV}\} \right)
    \nonumber \\
  &\le \epsilon \sum_{i=1}^{k} a_i
    + \sum_{i=1}^{k} \left| a_i - \|\wt{z_i}\|_{TV} \right|
    \nonumber \\
  &= \epsilon \| \chi\|_{TV}  
    + \sum_{i=1}^{k} \left| a_i - \|\wt{z_i}\|_{TV} \right|
    \qquad \l( \sum_{i=1}^k a_i = \|\chi\|_{TV} \r).
\end{align}
Then, we have that
\begin{align}
   d_W\l(\chi, \wt{z} \r) 
& \le \epsilon \l\|  \chi \r\|_{TV}+\sum_{i=1}^k \l| \int_{T_{i,\epsilon}}  \chi(\dif t)  - \frac{\l\|  \chi\r\|_{TV}}{\l\|  \widetilde{x}\r\|_{TV}} \int_{T_{i,\epsilon}}  \wt{x}(\dif t) \r| 
\qquad \l(  \mbox{see the second line of \eqref{eq:leg 2 pre}} \r)
\nonumber\\
& \le \epsilon \l\|  \chi\r\|_{TV} + \sum_{i=1}^k \l| \int_{T_{i,\epsilon}}  \chi(\dif t) -  \wt{x}(\dif t)\r| 
+ \sum_{i=1}^{k} \l\|   \wt{x}_{i}\r\|_{TV}\l| 1-\frac{\l\|  \chi\r\|_{TV}}{\l\|  \wt{x}\r\|_{TV}} \r|
\qquad \mbox{(triangle inequality)}\nonumber\\
& = \epsilon \l\|  \chi\r\|_{TV} +\sum_{i=1}^k \l| \int_{T_{i,\epsilon}}  h(\dif t) \r| 
+ \l| \l\|  \wt{x}\r\|_{TV} - \l\|  \chi\r\|_{TV}\r|
\qquad \l(\wt{x} = \widehat{x}|_{T_\epsilon} \r) 
\nonumber\\
& \le \epsilon \l\|  \chi\r\|_{TV} +\sum_{i=1}^k \l| \int_{T_{i,\epsilon}}  h(\dif t) \r| 
+ 2\l(\l\|b\r\|_2+\|b^0\|_2\r) \delta' 
\qquad \mbox{(see \eqref{eq:intermediate result})}
\nonumber\\
& \le \epsilon \l\|  \chi\r\|_{TV} + 2\l(\l\|b\r\|_2+\|b^0\|_2\r) \delta' 
+ 2\l(\l\|b\r\|_2+\|b^0\|_2\r) \delta'. 
\qquad \mbox{(see Lemma \ref{lem:dual 3})}
\end{align}

Substituting the above bound back into (\ref{eq:leg 2 pre}), we find
that 
\begin{equation}
d_{GW}\left( \chi,\widetilde{x}\right)\le6\left(\|b\|_{2}+\|b^0\|_{2}\right)\delta'+\epsilon\| \chi\|_{TV}.\label{eq:leg 2}
\end{equation}
Then combining (\ref{eq:leg 1}) and (\ref{eq:leg 2}) yields 
\begin{align*}
d_{GW}\left( \chi,\widehat{x}\right) & \le d_{GW}\left( x,\widetilde{x}\right)+d_{GW}\left(\widetilde{x},\widehat{x}\right)\\
 & \le\frac{2\|b\|_{2}\delta'}{\bar{f}}+6\left(\|b\|_{2}+\|b^0\|_{2}\right)\delta'+\epsilon\| \chi\|_{TV},
\end{align*}
which completes the proof of Lemma \ref{lem:EMD}.

\section{Proof of Proposition \ref{prop:dual construction}\label{sec:Proof-of-Proposition dual construction}}

Without loss of generality and for better clarity, suppose that $T=\{t_{i}\}_{i=1}^k$
 is an increasing sequence.
  Consider a positive scalar $\rho$ such that $\rho \le \epsilon \le \Delta/2$.
Consider also an increasing sequence $\{\tau_{l}\}_{l=1}^m\subset {I}=[0,1]$ such that $\tau_1 = 0$, $\tau_m = 1$, and every $T_{i,\rho}$ contains an even and nonzero number of
the remaining points. Let us define the polynomial
\begin{equation}
q^{\rho}(t)=\left|\begin{array}{cccc}
-F(t) & \phi_{1}(t) & \cdots & \phi_{m}(t)\\
-F(\tau_{1}) & \phi_{1}(\tau_{1}) & \cdots & \phi_{m}(\tau_{1})\\
-F(\tau_{2}) & \phi_{1}(\tau_{2}) & \cdots & \phi_{m}(\tau_{2})\\
\vdots & \vdots & \vdots & \vdots\\
-F(\tau_{m}) & \phi_{1}(\tau_{m}) & \cdots & \phi_{m}(\tau_{m})
\end{array}\right|,\qquad t\in I.
\label{eq:constructed poly}
\end{equation}
Note that $q^{\rho}(t)=0$ when $t\in \{\tau_{l}\}_{l=1}^m$. By assumption, 
$\{F\}\cup\{\phi_{j}\}_{j=1}^m$ form a T*-system on $I$. Therefore, invoking the first part of Definition \ref{def:(T-systems,-modified)-For}, we find that  $q^{\rho}$ is non-negative on ${T}_{\rho}^{C}$. We represent this
polynomial with $q^{\rho}=-\beta_{0}^{\rho}F+\sum_{j=1}^m (-1)^{j} \beta_{j}^{\rho}\phi_{j}$
and note that $\beta^\rho_0 = |\phi_j(\tau_i)|_{i,j=1}^m$.
By assumption also $\{\phi_j\}_{j=1}^m$ form a T-system on $I$ and therefore $\beta_{0}^{\rho}> 0$. This observation allows us to form the normalized polynomial 
\[
\dot{q}^{\rho}:=\frac{q^{\rho}}{\beta_0^{\rho}}=-F+\sum_{j=1}^m(-1)^j \frac{\beta_{j}^{\rho}}{\beta_0^{\rho}}\phi_{j}=:-F+\sum_{j=1}^m (-1)^j b_{j}^{\rho}\phi_{j}.
\]
Note also that the coefficients  $\{\beta^\rho_j\}_{j\in[0:m]}$ correspond to the minors in the second part of Definition \ref{def:(T-systems,-modified)-For}. Therefore,  
for each $j\in [0:m]$, we have that $|\beta^\rho_j|$ approaches zero at the same rate, 
as $\rho\rightarrow 0$. 
So for sufficiently small $\rho_0$, every $b_j^\rho$ is bounded in 
magnitude when $\rho \le \rho_0$; in particular, $|b_j^\rho| = \Theta(1)$.
This means that for sufficiently small
$\rho_0$, $\{\dot{q}^\rho: \rho\le \rho_0\}$ is bounded. 
Therefore, we can 
find a subsequence $\{\rho_l\}_l\subset [0,\rho_0]$ such that $\rho_l\rightarrow 0$ and
the subsequence $\{\dot{q}^{\rho_l}\}_l$ converges to the  polynomial 
\[
\dot{q}:=-F+\sum_{j=1}^m b_{j}\phi_{j}.
\]
Note that $b_j \ne 0$ for every $j\in[m]$; in particular, $|b_j| = \Theta(1)$. Hence 
the polynomial $\sum_{j=1}^m {b}_j\phi_j$ is nontrivial, namely does not uniformly vanish on $I$. (It would have sufficed to have some nonzero coefficient, say $b_{j_0}$, rather than  requiring all $\{b_j\}_j$ to be nonzero. However that would have made the statement of Definition \ref{def:(T-systems,-modified)-For} more cumbersome.) Lastly observe that $\dot{q}$ is non-negative on $I$ and vanishes on $T$ (as well as on the boundary of $I$). 
This completes the proof of Proposition \ref{prop:dual construction}.

\section{Proof of Proposition \ref{prop:A_diagonally_dominant}}
\label{sec:proof_A_diagonally_dominant}

In this proof, we will use the following result for strictly diagonally dominant matrices 
from \cite{plemmons1977}:
\begin{lem}
  If $ A $ is a strictly diagonally dominant matrix with positive entries on the main
  diagonal and negative entries otherwise, then $A$ is invertible and $ A^{-1} $ has
  non-negative entries.
  \label{lem:m-matrix}
\end{lem}

\subsubsection*{Proof of Proposition \ref{prop:A_diagonally_dominant}}

Let $ \hat{x} $ be a solution of \eqref{eq:mai} and $ h = x - \hat{x}$. 
Then, with $\phi_j(t) = \phi(t-s_j)$ for some $j$, by reverse triangle 
inequality we have
\begin{subequations}
  \begin{align*}
    \delta \geq 
    \left( 
      \sum_{j=1}^{m} \left(y(s_j) - \int_0^1 \phi_j(t) \hat{x}(\dif t)\right)^2 
    \right)^{1/2}
    &=
    \left(
      \sum_{j=1}^{m} \left( \phi_j(t)h(\dif t) + \delta_j   \right)^2
    \right)^{1/2} \\
    &\geq
    \left(
      \sum_{j=1}^{m} \left( \phi_j(t)h(\dif t)  \right)^2
    \right)^{1/2}
    - \left(
      \sum_{j=1}^m \delta_j^2
    \right)^{1/2}
    \\
    &\geq
    \left(
      \sum_{j=1}^{m} \left( \phi_j(t)h(\dif t)  \right)^2
    \right)^{1/2}
    - \delta,
  \end{align*}
\end{subequations}
and so
\begin{equation}
  \sum_{j=1}^{m} \left( \int_0^1 \phi_j(t) h(\dif t) \right)^2 \leq 4\delta^2 
  \implies 
  \left| \int_0^1 \phi_j(t) h(\dif t) \right| \leq 2\delta, \quad 
    \forall j \in [m].
  \label{eq:tolowerbound}
\end{equation}
We apply the reverse triangle inequality again to find a lower
bound of the left-hand side term in \eqref{eq:tolowerbound}:
\begin{subequations}
  \begin{align}
    \left| \int_0^1 \phi_j(t) h(\dif t) \right| 
    &\geq
      \left| \int_{T_{i,\epsilon}} \phi_j(t) \hat{x}(\dif t) 
        - a_i \phi_j(t_i)\right| 
      \label{term1} \\
    &- \sum_{l \neq i} 
      \left| \int_{T_{l,\epsilon}} \phi_j(t) \hat{x}(\dif t) 
        - a_l \phi_j(t_l) \right| 
      \label{term2} \\
    &- \left| \int_{T_{\epsilon}^C} \phi_j(t) \hat{x}(\dif t) \right|.
      \label{term3}
  \end{align}
  \label{eq:allterms}
\end{subequations}
We now need to lower bound the term in \eqref{term1} and upper bound
the terms in \eqref{term2}, \eqref{term3}. For the first one, we obtain:
\begin{align}
    \left| 
      \int_{T_{i,\epsilon}} \phi_j(t) \hat{x}(\dif t) 
      - a_i \phi_j(t_i)\right| 
    &= 
    \left| 
      \int_{T_{i,\epsilon}} \phi_j(t) \hat{x}(\dif t) 
      - a_i \phi_j(t_i)
      + \phi_j(t_i) \int_{T_{i,\epsilon}} \hat{x}(\dif t)  
      - \phi_j(t_i) \int_{T_{i,\epsilon}} \hat{x}(\dif t)  
    \right| 
  \nonumber \\
    &\geq
    \left| \phi_j(t_i) \right| 
    \left| \int_{T_{i,\epsilon}} \hat{x}(\dif t) - a_i \right|
    - \int_{T_{i,\epsilon}} \left| \phi_j(t) - \phi_j(t_i)\right| \hat{x}(\dif t)
  \nonumber \\
    &\geq
    \left| \phi_j(t_i) \right| 
    \left| \int_{T_{i,\epsilon}} \hat{x}(\dif t) - a_i \right|
    - L \int_{T_{i,\epsilon}} \left| t - t_i \right| \hat{x}(\dif t) 
  \nonumber \\
    &\geq
    \left| \phi_j(t_i) \right| 
    \left| \int_{T_{i,\epsilon}} \hat{x}(\dif t) - a_i \right|
    - L \epsilon \int_{T_{i,\epsilon}} \hat{x}(\dif t).
  \label{boundterm1}
\end{align}
Therefore, from \eqref{boundterm1}, we obtain:
\begin{equation}
    \left| 
      \int_{T_{i,\epsilon}} \phi_j(t) \hat{x}(\dif t) 
      - a_i \phi_j(t_i)\right|
    \geq
    \left|\phi_j(t_i) \right| 
    \left| \int_{T_{i,\epsilon}} \hat{x}(\dif t) - a_i \right|
    - L  \epsilon \|\hat{x}\|_{TV}.
    \label{bound1}
\end{equation}
For the term \eqref{term2}, we have:
\begin{align*}
    \left| \int_{T_{l,\epsilon}} \phi_j(t) \hat{x}(\dif t) - a_l \phi_j(t_l) \right|
    &= 
    \left| 
      \int_{T_{l,\epsilon}} \phi_j(t) \hat{x}(\dif t) - a_l \phi_j(t_l) 
      + \phi_j(t_l) \int_{T_{l,\epsilon}} \hat{x}(\dif t)
      - \phi_j(t_l) \int_{T_{l,\epsilon}} \hat{x}(\dif t)
    \right| \\
    &\leq
    \left| \phi_j(t_l)\right| 
    \left| \int_{T_{l,\epsilon}} \hat{x}(\dif t) - a_l \right|
    + \int_{T_{l,\epsilon}} \left| \phi_j(t) - \phi_j(t_l) \right| \hat{x}(\dif t) \\
    &\leq
    \left| \phi_j(t_l)\right| 
    \left| \int_{T_{l,\epsilon}} \hat{x}(\dif t) - a_l \right|
    + L \int_{T_{l,\epsilon}} \left| t - t_l \right| \hat{x}(\dif t) \\
    &\leq
    \left| \phi_j(t_l)\right| 
    \left| \int_{T_{l,\epsilon}} \hat{x}(\dif t) - a_l \right|
    + L \epsilon \int_{T_{l,\epsilon}} \hat{x}(\dif t),
\end{align*}
so
\begin{equation}
  \sum_{l \neq i} 
    \left| \int_{T_{l,\epsilon}} \phi_j(t) \hat{x}(\dif t) 
        - a_l \phi_j(t_l) \right| 
  \leq
  \sum_{l \ne i}
    \left(
      \left| \phi_j(t_l) \right| 
      \left| \int_{T_{l,\epsilon}} \hat{x}(\dif t) - a_l \right|
    \right)
  + L \epsilon
  \sum_{l \ne i} 
    \int_{T_{l,\epsilon}} 
       \hat{x}(\dif t).
  \label{bound2}
\end{equation}
Finally, for the term \eqref{term3}, we have:
\begin{equation}
  \left| \int_{T_{\epsilon}^C} \phi_j(t) \hat{x}(\dif t) \right|
  \leq 
  \max_{t \in T_{\epsilon}^C} \left| \phi_j(t) \right|
  \int_{T_{\epsilon}^C} \hat{x}(\dif t)
  \leq
  \phi^{\infty} 
  \left(
    \frac{2 \|b \|_2 \delta}{\bar{f}}
  \right).
  \label{bound3}
\end{equation}
Let us denote
\begin{equation*}
  z_i = \left| \int_{T_{i,\epsilon}} \hat{x}(\dif t) - a_i \right|,
\end{equation*}
for all $i \in [k]$. Then, by combining \eqref{eq:allterms} 
with the bounds \eqref{boundterm1},\eqref{bound2}
and \eqref{bound3}, we obtain the $j$-th row of a linear system:
\begin{equation}
  2\delta \geq
  \left| \phi_j(t_i) \right| z_i - L \epsilon \int_{T_{i,\epsilon}} \hat{x}(\dif t)
  - \sum_{l \ne i} \left| \phi_j(t_l) \right| z_l - L\epsilon \sum_{l \ne i} \int_{T_{l,\epsilon}} \hat{x}(\dif t)
  - \phi^{\infty} \frac{2 \|b\|_2 }{\bar{f}} \delta.
  \label{row_unfinished}
\end{equation}
By using \eqref{row_unfinished} along with
\begin{equation*}
  \int_{T_{i,\epsilon}} \hat{x}(\dif t)
  + \sum_{l \ne i} \int_{T_{l,\epsilon}} \hat{x}(\dif t)
  = \int_{T_{\epsilon}} \hat{x}(\dif t) 
  \leq \int_{I} \hat{x}(\dif t)
  = \|\hat{x}\|_{TV},
\end{equation*}
we obtain: 
\begin{equation}
  2 \left( 1 + \frac{\phi^{\infty}\|b\|_2}{\bar{f}} \right) \delta
  + \epsilon L \|\hat{x}\|_{TV} 
  \geq
  \left| \phi_j(t_i) \right| z_i
  - \sum_{l \ne i} \left| \phi_j (t_l) \right| z_l.
  \label{ineqline}
\end{equation}
Now, for all $ i $, we select the $ j=l(i) $,
the index corresponding to the closest sample 
as defined in Definition \ref{def:min_sep}.
The inequalities in \eqref{ineqline} can be written as
\begin{equation}
  Az \leq v,
  \label{linsys}
\end{equation}
where $ A, z $ and $ v $ are defined as
\begin{equation*}
  A = 
  \begin{bmatrix}
    |\phi_1(t_1)|  & -|\phi_1(t_2)| & \ldots & -|\phi_1(t_k)| \\
    -|\phi_2(t_1)| & |\phi_2(t_2)|  & \ldots & -|\phi_2(t_k)| \\
    \vdots         & \vdots         & \ddots & \vdots         \\
    -|\phi_k(t_1)| & -|\phi_k(t_2)| & \ldots & |\phi_k(t_k)|  \\
  \end{bmatrix},
  \quad
  z = 
  \begin{bmatrix}
    z_1 \\
    z_2 \\
    \vdots \\
    z_k 
  \end{bmatrix},
  \quad
  v = 
  \left(
    2 \left( 1 + \frac{\phi^{\infty}\|b\|_2}{\bar{f}} \right) \delta
      + \epsilon L \|\hat{x}\|_{TV} 
  \right)
  \begin{bmatrix}
    1 \\ 1 \\ \vdots \\ 1
  \end{bmatrix}.
\end{equation*}

Because $A$ is strictly diagonally dominant, Lemma \ref{lem:m-matrix} holds and 
therefore $A^{-1}$ exists and has non-negative entries, so when we
multiply \eqref{linsys} by $A^{-1}$, the sign does not change:
\begin{equation}
  z \leq A^{-1} v,
  \label{eq:final_ineq}
\end{equation}
where we can bound the entries of $A^{-1}$ \cite{varah1975}:
\begin{equation*}
  \| A^{-1} \|_{\infty} < \frac{1}{\min_j (\phi_j(t_j) - \sum_{i \ne j} \phi_j(t_i))}.
\end{equation*}
The proof of Proposition \ref{prop:A_diagonally_dominant} is now complete, 
since \eqref{eq:final_ineq} is equivalent to our error 
bound \eqref{eq:bound_near_support}.

\section{Proof of Proposition \ref{prop:sampling_locations}}
\label{sec:sampling_locations}
For the sake of simplicity, let $s_i$ be the closest sample $s_{l(i)}$
to the source $t_i$, as defined in Definition \ref{def:min_sep}.
For a fixed $i$, assume (without loss of generality, as we will see later) that $ s_i < t_i $. 
It follows that
\begin{align*}
  |s_i - t_l| &= s_i - t_l \geq (i-l)\Delta - \lambda\Delta, 
    \quad
    \forall l < i, \\
  |s_i - t_l| &= t_l - s_i \geq (l-i)\Delta + \lambda\Delta,
    \quad
    \forall l > i,
\end{align*}
and so 
\begin{align*}
  \phi(|s_i - t_l|) \leq \phi((i-l)\Delta - \lambda\Delta), 
    \quad \forall l < i, \\
  \phi(|s_i - t_l|) \leq \phi((l-i)\Delta + \lambda\Delta),
    \quad \forall l > i.
\end{align*}
Then we have
\begin{align}
  \sum_{l \ne i} \phi(|s_i - t_l|) 
  &= 
  \sum_{l=1}^{i-1} \phi(|s_i - t_l|) + \sum_{l=i+1}^{k} \phi(|s_i - t_l|) 
  \nonumber \\
  &\leq
  \sum_{l=1}^{i-1} \phi((i-l)\Delta - \lambda \Delta) 
    + \sum_{l=i+1}^{k} \phi( (l-i)\Delta + \lambda \Delta) 
  \nonumber \\
  &=
  \sum_{l=1}^{i-1} \phi(l\Delta - \lambda \Delta)
    + \sum_{l=1}^{k-i} \phi(l\Delta + \lambda \Delta).
    \label{sumsldelta}
\end{align}
We now want to find upper bounds for each of the two sums 
in \eqref{sumsldelta}. We will derive the bound for the first term,
as the second one is similar. We have that
\begin{equation}
    \sum_{l=1}^{i-1} \phi(l\Delta - \lambda \Delta)
    = \phi(\Delta - \lambda \Delta)
      + \frac{1}{\Delta}
        \sum_{l=2}^{i-1} \phi(l\Delta - \lambda \Delta) \Delta,
    \label{sumeqr}
\end{equation}
and the sum in the previous equation is a lower Riemann sum
(note that $ \phi $ is decreasing in $[0,1]$)
\begin{equation*}
  S = \sum_{l=2}^{i-1} \phi(x_l^*) (x_l - x_{l-1})
\end{equation*}
of $ \phi(x) $ over $ [\Delta - \lambda\Delta, (i-1)\Delta - \lambda\Delta ] $, 
with partition and $ x_l^*$ chosen as follows:
\begin{align*}
  [x_l, x_{l-1}] &= [l\Delta - \lambda\Delta, (l-1)\Delta - \lambda\Delta],
    \quad &l = 2\ldots,i-1,
  \\
  x_l^* &= l\Delta - \lambda\Delta,
    \quad &l = 2,\ldots,i-1.
\end{align*}
Therefore, the sum $ S $ is less than or equal to the integral:
\begin{equation}
  \sum_{l=2}^{i-1} \phi(l\Delta - \lambda \Delta) \Delta
  \leq
  \int_{\Delta - \lambda\Delta}^{(i-1)\Delta - \lambda\Delta} \phi(x) \dif x.
  \label{intbound}
\end{equation}
By substituting \eqref{intbound} into \eqref{sumeqr}, we obtain
\begin{equation*}
    \sum_{l=1}^{i-1} \phi(l\Delta - \lambda \Delta)
    \leq \phi(\Delta - \lambda \Delta)
      + \frac{1}{\Delta}
        \int_{\Delta - \lambda\Delta}^{(i-1)\Delta - \lambda\Delta} \phi(x) \dif x.
\end{equation*}
We can obtain a similar upper bound for the second sum in \eqref{sumsldelta} and then
\begin{equation}
  \sum_{l \ne i} \phi(|s_i - t_l|) 
  \leq
  \phi(\Delta - \lambda \Delta) +
  \phi(\Delta + \lambda \Delta) +
  \frac{1}{\Delta}
    \int_{\Delta - \lambda\Delta}^{(i-1)\Delta - \lambda\Delta} \phi(x) \dif x +
  \frac{1}{\Delta}
    \int_{\Delta + \lambda\Delta}^{(k-i)\Delta + \lambda\Delta} \phi(x) \dif x,
  \quad
  \forall i \in [k].
  \label{libound2}
\end{equation}
We can further upper bound the right hand side over all $ i \in [k]$
and this bound corresponds to the case when the source $ t_i$ is
in the middle of the unit interval (at $\frac12$) and the 
sources $t_1$ and $t_k$ are at $0$ and $1$ respectively:
\begin{equation*}
  \int_{\Delta - \lambda\Delta}^{(i-1)\Delta - \lambda\Delta} \phi(x) \dif x +
  \int_{\Delta + \lambda\Delta}^{(k-i)\Delta + \lambda\Delta} \phi(x) \dif x
  \leq
  \int_{\Delta - \lambda\Delta}^{1/2 - \lambda\Delta} \phi(x) \dif x +
  \int_{\Delta + \lambda\Delta}^{1/2 + \lambda\Delta} \phi(x) \dif x,
  \quad
  \forall i \in [k],
\end{equation*}
and therefore we have
\begin{equation*}
  \sum_{l \ne i} \phi(|s_i - t_l|) 
  \leq
  \phi(\Delta - \lambda \Delta) +
  \phi(\Delta + \lambda \Delta) +
  \frac{1}{\Delta}
    \int_{\Delta - \lambda\Delta}^{1/2 - \lambda\Delta} \phi(x) \dif x +
  \frac{1}{\Delta}
    \int_{\Delta + \lambda\Delta}^{1/2 + \lambda\Delta} \phi(x) \dif x,
  \quad
  \forall i \in [k].
\end{equation*}
In order to find $\lambda_0$, we solve \eqref{eq:condition_lambda}
since $ |s_i - t_i| \leq \lambda\Delta $ implies $\phi(|s_i - t_i|) \geq \phi(\lambda\Delta)$.

We note that if we only have three sources, then the integral terms should not be included, 
and if we have four sources, then the last integral term should not be included.

\section{Proof of Theorem \ref{thm:grouping_sources}}
\label{sec:proof_grouping_sources}

The proof of Theorem \ref{thm:grouping_sources} involves the same ideas as Theorem \ref{thm:noisy_expl}. 
The differences are in the analysis of Proposition \ref{prop:A_diagonally_dominant}. 
We continue this analysis from \eqref{eq:allterms}, where we lower bound the left-hand side term 
of \eqref{eq:tolowerbound}:
\begin{subequations}
  \label{eq:allterms_2}
  \begin{align}
    2\delta \geq \left| \int_0^1 \phi_j(t) h(\dif t) \right| 
    &\geq
      \left| \int_{\tilde{T}_{i,\epsilon}} \phi_j(t) \hat{x}(\dif t) 
        - \sum_{r \in [k_i]} a_{ir} \phi_j(t_{ir})\right| 
      \label{term1_2} \\
    &- \sum_{l \neq i} 
      \left| \int_{\tilde{T}_{l,\epsilon}} \phi_j(t) \hat{x}(\dif t) 
        - \sum_{r \in [k_l]} a_{lr} \phi_j(t_{lr}) \right| 
      \label{term2_2} \\
    &- \left| \int_{T_{\epsilon}^C} \phi_j(t) \hat{x}(\dif t) \right|,
      \label{term3_2}
  \end{align}
\end{subequations}
where, in each term, the sum from $r=1$ to $k_i$ is over all the 
true sources in $\tilde{T}_{i,\epsilon}$ (for all $ i \in [\tilde{k}]$). In order to obtain bounds for 
the terms \eqref{term1_2} and \eqref{term2_2}, we need the following fact:
\begin{equation}
  \exists \xi_i \in [\argmin_{r\in[k_i]} \phi_j(t_{ir}), \argmax_{r\in [k_i]} \phi_j(t_{ir})]
  \quad
  \text{such that} 
  \quad
  \phi_j(\xi_i) \sum_{r=1}^{k_i} a_{ir} 
  = \sum_{r=1}^{k_r} a_{rk} \phi_j(t_{ir}), 
  \quad
  \forall r \in [\tilde{k}].
  \label{ivt}
\end{equation}
This comes from the continuity of $\phi_j$ and intermediate value theorem, since:
\begin{equation*}
  \min_{k} \phi_j(t_{ir})
  \leq
  \frac{\sum_{r=1}^{k_i} a_{ir} \phi_j(t_{ir})}
  {\sum_{r=1}^{k_i} a_{ir}}
  \leq
  \max_{k} \phi_j(t_{ir}).
\end{equation*}
We proceed as before to find a lower bound for \eqref{term1_2} and an upper bound
for \eqref{term2_2}, while the upper bound for \eqref{term3_2} is the same. For \eqref{term1_2}:
\begin{align*}
  \left| 
    \int_{\tilde{T}_{i,\epsilon}} \phi_j(t) \hat{x}(\dif t) 
    - \sum_{r=1}^{k_i} a_{ir} \phi_j(t_{ir})
  \right| 
  &= 
  \left| 
    \int_{\tilde{T}_{i,\epsilon}} \phi_j(t) \hat{x}(\dif t) 
    - \phi_j(\xi_i) \sum_{r=1}^{k_i} a_{ir} 
  \right| \\
  &=\left| 
    \int_{\tilde{T}_{i,\epsilon}} \phi_j(t) \hat{x}(\dif t) 
    - \phi_j(\xi_i) \sum_{r=1}^{k_i} a_{ir} 
    + \phi_j(\xi_i) \int_{\tilde{T}_{i,\epsilon}} \hat{x}(\dif t)  
    - \phi_j(\xi_i) \int_{\tilde{T}_{i,\epsilon}} \hat{x}(\dif t)  
  \right| \\
  &\geq
  \left| \phi_j(\xi_i) \right| 
  \left| \int_{\tilde{T}_{i,\epsilon}} \hat{x}(\dif t) - \sum_{r=1}^{k_i} a_{ir} \right|
  - \int_{\tilde{T}_{i,\epsilon}} \left| \phi_j(t) - \phi_j(\xi_i)\right| \hat{x}(\dif t)\\
  &\geq
  \left| \phi_j(\xi_i) \right| 
  \left| \int_{\tilde{T}_{i,\epsilon}} \hat{x}(\dif t) - \sum_{r=1}^{k_i} a_{ir} \right|
  - L \int_{\tilde{T}_{i,\epsilon}} \left| t - \xi_i \right| \hat{x}(\dif t) \\
  &\geq
  \left|\phi_j(\xi_i) \right| 
  \left| \int_{\tilde{T}_{i,\epsilon}} \hat{x}(\dif t) - \sum_{r=1}^{k_i} a_{ir} \right|
  - L  (2 k_i-1) \epsilon \int_{\tilde{T}_{i,\epsilon}} \hat{x}(\dif t),
\end{align*}
where the width of $\tilde{T}_{i,\epsilon}$ 
is at most $2 k_i \epsilon$ 
and $ \xi_i \in \tilde{T}_{i,\epsilon}$ 
is chosen according to \eqref{ivt}, 
so the distance $|t - \xi_i|$ for $ t \in \tilde{T}_{i,\epsilon}$ is 
at most $(2 k_i-1) \epsilon$.
For the second term \eqref{term2_2}:
\begin{align*}
  \left| \int_{\tilde{T}_{l,\epsilon}} \phi_j(t) \hat{x}(\dif t) 
    - \sum_{r=1}^{k_l} a_{lr} \phi_j(t_{lr}) 
  \right|
  &=
  \left| \int_{\tilde{T}_{l,\epsilon}} \phi_j(t) \hat{x}(\dif t) 
    - \phi_j(\xi_{l}) \sum_{r=1}^{k_l} a_{lr} 
  \right| \\
  &= 
  \left| 
    \int_{\tilde{T}_{l,\epsilon}} \phi_j(t) \hat{x}(\dif t) 
    - \phi_j(\xi_{l}) \sum_{r=1}^{k_l} a_{lr} 
    + \phi_j(\xi_l) \int_{\tilde{T}_{l,\epsilon}} \hat{x}(\dif t)
    - \phi_j(\xi_l) \int_{\tilde{T}_{l,\epsilon}} \hat{x}(\dif t)
  \right| \\
  &\leq
  \left| \phi_j(\xi_l)\right| 
  \left| \int_{\tilde{T}_{l,\epsilon}} \hat{x}(\dif t) - \sum_{r=1}^{k_l} a_{lr} \right|
  + \int_{\tilde{T}_{l,\epsilon}} \left| \phi_j(t) - \phi_j(\xi_l) \right| \hat{x}(\dif t) \\
  &\leq
  \left| \phi_j(t_l)\right| 
  \left| \int_{\tilde{T}_{l,\epsilon}} \hat{x}(\dif t) - \sum_{r=1}^{k_l} a_{lr} \right|
  + L \int_{\tilde{T}_{l,\epsilon}} \left| t - \xi_l \right| \hat{x}(\dif t) \\
  &\leq
  \left| \phi_j(t_l)\right| 
  \left| \int_{\tilde{T}_{l,\epsilon}} \hat{x}(\dif t) - \sum_{r=1}^{k_l} a_{lr} \right|
  + L (2k_l-1) \epsilon \int_{\tilde{T}_{l,\epsilon}} \hat{x}(\dif t) 
\end{align*}
so
\begin{equation*}
  \sum_{l \neq i} 
    \left| \int_{\tilde{T}_{l,\epsilon}} \phi_j(t) \hat{x}(\dif t) 
      - \sum_{r=1}^{k_l} a_{lr} \phi_j(t_{lr}) 
    \right| 
  \leq
  \sum_{l \ne i}
    \left(
      \left| \phi_j(\xi_l) \right| 
      \left| \int_{\tilde{T}_{l,\epsilon}} \hat{x}(\dif t) - \sum_{r=1}^{k_l} a_{lr} \right|
    \right)
  + L \epsilon 
  \sum_{l \ne i} 
    (2 k_l-1) 
    \int_{\tilde{T}_{l,\epsilon}} \hat{x}(\dif t).
\end{equation*}
Let
\begin{equation*}
  \tilde{z}_i = 
  \left| 
    \int_{\tilde{T}_{i,\epsilon}} \hat{x}(\dif t) - \sum_{r=1}^{k_i} a_{ir} 
  \right|
\end{equation*}
and we obtain an inequality as before:
\begin{equation*}
  2 \left( 1 + \frac{\phi^{\infty}\|b\|_2}{\bar{f}} \right) \delta
  + (2k-1) \epsilon L \|\hat{x}\|_{TV}
  \geq
  \left| \phi_j(\xi_i) \right| \tilde{z}_i
  - \sum_{l \ne i} \left| \phi_j (\xi_l) \right| \tilde{z}_l,
\end{equation*}
where we obtained the second constant as follows:
\begin{equation*}
  L\epsilon (2 k_i-1) \int_{\tilde{T}_{i,\epsilon}} \hat{x}(\dif t) 
  + L\epsilon  \sum_{l \ne i} (2 k_l-1) \int_{\tilde{T}_{l,\epsilon}} \hat{x}(\dif t)
  = L\epsilon \sum_{l=1}^{\tilde{k}} (2 k_l-1) \int_{\tilde{T}_{l,\epsilon}} \hat{x}(\dif t)
  \leq L \epsilon \|\hat{x}\|_{TV} \sum_{l=1}^{\tilde{k}} (2 k_l-1)
  \leq (2k-1)  L \epsilon \|\hat{x}\|_{TV} 
\end{equation*}
and, for all $i \in [\tilde{k}]$, we select $j(i) = \argmin_j |s_j - \xi_i|$.
The linear system is 
\begin{equation*}
  \tilde{A} \tilde{z} \leq \tilde{v},
\end{equation*}
with:
\begin{equation*}
  \tilde{A} = 
  \begin{bmatrix}
    |\phi_1(\xi_1)|  & -|\phi_1(\xi_2)| & \ldots & -|\phi_1(\xi_{\tilde{k}})| \\
    -|\phi_2(\xi_1)| & |\phi_2(\xi_2)|  & \ldots & -|\phi_2(\xi_{\tilde{k}})| \\
    \vdots         & \vdots         & \ddots & \vdots         \\
    -|\phi_{\tilde{k}}(\xi_1)| & -|\phi_{\tilde{k}}(\xi_2)| & \ldots & |\phi_{\tilde{k}}(\xi_{\tilde{k}})|  \\
  \end{bmatrix},
  \quad
  \tilde{z} = 
  \begin{bmatrix}
    z_1 \\
    z_2 \\
    \vdots \\
    z_{\tilde{k}}
  \end{bmatrix},
  \quad
  \tilde{v} = 
  \left(
    2 \left( 1 + \frac{\phi^{\infty}\|b\|_2}{\bar{f}} \right) \delta
      + (2k-1) L\epsilon \|\hat{x}\|_{TV}
  \right)
  \begin{bmatrix}
    1 \\ 1 \\ \vdots \\ 1
  \end{bmatrix}.
\end{equation*}

In Appendix \ref{sec:sampling_locations} we discuss what the choice of $\lambda$ should be so that,
if $ |t_i - s_i| \leq \lambda \Delta $, the matrix $ A $ is strictly diagonally dominant.
Here, the matrix $\tilde{A} $ is similar to $A$ except that we evaluate $ \phi $ 
at $ |\xi_i - s_i| $, where $\xi_i$ corresponds to a group of sources in $ \tilde{T}_{i,\epsilon} $
that are located within distances smaller than $2\epsilon$ 
and $s_i$ is the closest sample
to $\xi_i$. Given that the minimum separation between $\xi_i$ 
sources is $2\epsilon$, the 
analysis in Appendix \ref{sec:sampling_locations} is the same, so $\tilde{A}$ is strictly 
diagonally dominant if
\begin{equation*}
  |\xi_i - s_i| \leq 2 \lambda \epsilon, 
  \quad
  \forall i \in [\tilde{k}]
\end{equation*}
and $\lambda$ is chosen to satisfy \eqref{eq:condition_lambda} where 
we take $\Delta = 2 \epsilon$. For the value of $\lambda$ found
this way, we select the sampling locations uniformly at 
intervals of $2\lambda \epsilon$.

\section{Proof of Lemma \ref{lem:L cte} \label{sec:proof of L cte}}

For real $s$, let us first study the Lipschitz constant of the operator that takes a non-negative measure $x$ supported on $I$ to $\int_I g(s-t) x(\dif t)$, where $g(t) = e^{-\frac{t^2}{\sigma^2}}$.
To that end,  consider a pair of non-negative measures $z_1,z_2$ supported on $I$ such that $\|z_1\|_{TV}=\|z_2\|_{TV}$. Their Wasserstein distance $d_W(z_1,z_2)$ was defined in \eqref{eq:def of EMD}. The dual of Program \eqref{eq:def of EMD} is in fact (see, for example \cite{santambrogio2015optimal}, Chapter 5)
\begin{equation}
  \max_{\omega} \{ \int_I \omega(t) (z_1-z_2)(\dif t) :
  \omega \text{ is } 1\text{-Lipschitz} \}.
  \label{eq:dual of EMD}
\end{equation}
where the supremum is over all (measurable) functions $\omega:I\rightarrow \mathbb{R}$. 
Consider the particular choice of $\omega(t)= \alpha\cdot  g (s-t)$ for positive $\alpha$ to be set shortly. 
Let us see if this choice of $\omega$ is feasible for Program \eqref{eq:dual of EMD}. For $t_1,t_2\in I$, we write that 
\begin{align}
  \omega(t_1)-\omega(t_2) & = \alpha \l( g(s-t_1) - g(s-t_2) \r) \nonumber\\
  & = \alpha \int_{t_1}^{t_2} g'(s-t) \,\,dt \nonumber\\
  & \le \alpha \max_{t \in [-1,1]} |g'(t)| \cdot |t_2-t_1| \nonumber\\
  & = \frac{\alpha}{\sigma}\sqrt{\frac{2}{e}} \cdot |t_2-t_1| \nonumber\\
  & = |t_2-t_1|, 
\end{align}
where we set $\alpha=\frac{\sigma \sqrt{2e}}{2}$ in the last line above, for small enough $\sigma$ 
(specifically, for $\sigma < \sqrt{2}$).
Therefore, with this choice of $\alpha$, the function $\omega$ specified above is feasible 
for Program \eqref{eq:dual of EMD}. From this observation and the strong duality 
between Programs \eqref{eq:def of EMD} and \eqref{eq:dual of EMD}, it follows that 
\begin{align}
  \frac{\sigma\sqrt{2e}}{2} \int_I g(s-t) (z_1(\dif t)-z_2(\dif t)) & \le d_W(z_1,z_2).
\end{align}
Combined with the other direction, we find that 
\begin{align}
  \frac{\sigma \sqrt{2e}}{2} \l| \int_I g(s-t) (z_1(\dif t)-z_2(\dif t)) \r| & \le d_W(z_1,z_2),
\label{eq:dW is L}
\end{align}
namely that the map $z\rightarrow \int_I g(s-t)z(\dif t)$ is $(2/\sigma\sqrt{2e})$-Lipschitz with respect to the Wasserstein distance. 

It is easy to extend the above conclusion to the generalised Wasserstein distance. Consider a pair of non-negative measures $x_1,x_2$ supported on $I$, and a pair of non-negative measures $z_1,z_2$ on $I$ such that $\|z_1\|_{TV}=\|z_2\|_{TV}$.  Then, using the triangle inequality, we write that 
\begin{align}
  & \l| \int_I g(s-t) (x_1(\dif t)-x_2(\dif t)) \r| \nonumber\\
  & \le \l| \int_I g(s-t) (x_1(\dif t)-z_1(\dif t)) \r|
  + \l| \int_I g(s-t) (z_1(\dif t)-z_2(\dif t)) \r| 
  + \l| \int_I g(s-t) (z_2(\dif t)-x_2(\dif t)) \r| \nonumber\\
  & \le \int_I \l| x_1(\dif t)-z_1(\dif t)\r| + \l| \int_I g(s-t) (z_1(\dif t)-z_2(\dif t)) \r| + \int_I \l| z_2(\dif t)-x_2(\dif t)\r| \qquad \l(g(t)\le 1 \r) \nonumber\\
  & = \| x_1-z_1\|_{TV}+ \l| \int_I g(s-t) (z_1(\dif t)-z_2(\dif t)) \r| + \|z_2-x_2\|_{TV} \nonumber\\
  & \le \| x_1-z_1\|_{TV}+ \frac{2}{\sigma\sqrt{2e}} d_W(z_1,z_2) + \|z_2-x_2\|_{TV}
  \qquad \mbox{(see \eqref{eq:dW is L})} \nonumber\\
  & \le \frac{2}{\sigma \sqrt{2e}}
    \left( \| x_1-z_1\|_{TV}+ d_W(z_1,z_2)+ \|z_2-x_2\|_{TV} \right)
    \qquad \left(\text{for } \sigma \le \frac{2}{\sqrt{2e}}\right).
\end{align}
The choice of $z_1,z_2$ above was arbitrary and therefore, recalling the definition of $d_{GW}$ in \eqref{eq:def of gen EMD}, we find that
\begin{equation}
  \l| \int_I g(s-t) (x_1(\dif t)-x_2(\dif t)) \r| \le \frac{2}{\sigma \sqrt{2e}} d_{GW}(x_1,x_2),
\end{equation}
namely the map $x\rightarrow\int_I g(s-t) x(\dif t)$ is a 
  $\frac{2}{\sigma \sqrt{2e}}$-Lipschitz operator.
It immediately follows that the operator $\Phi$ that was formed using the sampling points $S=\{s_j\}_{j=1}^m$ in \eqref{eq:Phi} is 
  $\frac{2\sqrt{m}}{\sigma \sqrt{2e}}$-Lipschitz.
This completes the proof of Lemma \ref{lem:L cte}. 

\section{Proof of Lemma \ref{lem:Gaussian is Tstar sys} \label{sec:proof of lemma Gaussian is Tstar sys}}

Following the definition of T*-systems in  Definition \ref{def:(T-systems,-modified)-For},  consider an increasing
sequence $\{\tau_{l}\}_{l=0}^m\subset I$ 
such that $\tau_0 = 0$, $\tau_m = 1$, and  except one more point 
(say $\tau_{\underline{l}}$),  the rest of points belong 
to $T_{\rho}$, the $\rho$-neighbourhood of the support $T\subset \interior(I)$. 

We also select the subset of samples $S$ of size $2k+2$ that is closest to the 
support $T$ (in Hausdorff distance), so without loss of generality, we set $m=2k+2$.
With this assumption, the setup in Definition \ref{def:(T-systems,-modified)-For} 
forces that  every neighbourhood ${T}_{i,\rho}$ contain exactly two points, 
say $t_{i}$ and $t_{i,\rho}:=t_i +\rho$ to simplify the presentation.  
Then the determinant in part 1 of Definition \ref{def:(T-systems,-modified)-For} can be written as  
\begin{align}
  M^{\rho}  & =\left|\begin{array}{cccc}
  F\left(0\right) & g\left(s_{1}\right) & \cdots & g\left(s_{m}\right)\\
  F\left(t_{1}\right) & g\left(t_1-s_{1}\right) & \cdots & g\left(t_1-s_{m}\right)\\ 
  F\left(t_{1,\rho}\right) & g\left(t_{1,\rho}-s_{1}\right) & \cdots & g\left(t_{1,\rho}-s_{m}\right)\\
  \vdots  \\
  F\left(\tau_{\underline{l}}\right) & g\left(\tau_{\underline{l}}-s_{1}\right) & \cdots & g\left(\tau_{\underline{l}}-s_{m}\right)\\
  \vdots \\
  F\left(t_{k}\right) & g\left(t_k-s_{1}\right) & \cdots & g\left(\tau-s_{m}\right)\\
  F\left(t_{k,\rho}\right) & g\left(t_{k,\rho}-s_{1}\right) & \cdots & g\left(t_{k,\rho}-s_{m}\right)\\
  F\left(1\right) & g\left(1-s_{1}\right) & \cdots & g\left(1-s_{m}\right)
  \end{array}\right|\nonumber\\
  & = \left|\begin{array}{cccc}
  f_{0} & g\left(s_{1}\right) & \cdots & g\left(s_{m}\right)\\
  0 & g\left(t_1-s_{1}\right) & \cdots & g\left(t_1-s_{m}\right)\\ 
  0 & g\left(t_{1,\rho}-s_{1}\right) & \cdots & g\left(t_{1,\rho}-s_{m}\right)\\
  \vdots  \\
  F\l(\tau_{\underline{l}} \r) & g\left(\tau_{\underline{l}}-s_{1}\right) & \cdots & g\left(\tau_{\underline{l}}-s_{m}\right)\\
  \vdots \\
  0 & g\left(t_k-s_{1}\right) & \cdots & g\left(\tau-s_{m}\right)\\
  0 & g\left(t_{k,\rho}-s_{1}\right) & \cdots & g\left(t_{k,\rho}-s_{m}\right)\\
  f_{1} & g\left(1-s_{1}\right) & \cdots & g\left(1-s_{m}\right)
  \end{array}\right|. 
  \qquad (\text{evaluating } F(t) \text{ as in } \eqref{eq:f example})
  \label{eq:m rho}
\end{align}
We will now need the following lemma in order to simplify the determinant above. 
The result is proved in Appendix \ref{sec:proof of determinant bounds lem}.

\begin{lem}
  \label{lem:determinant bounds}

  Let $A,B \in \mathbb{R}^{m \times m}$ with $m \geq 2$ and $\det(A) > 0$. 
  If $0 \leq \epsilon \leq \frac{8}{34 m \rho(A^{-1} B)}$, then
  \begin{equation*}
    \det(A) \left(1 - \frac{17 \sqrt{e}}{8} m \epsilon  \rho(A^{-1}B) \right)
    \leq det(A + \epsilon B) \leq
    \det(A) \left(1 + \frac{17 \sqrt{e}}{8} m \epsilon \rho(A^{-1}B) \right),
  \end{equation*}
  where $\rho(X)$ is the spectral radius of the matrix $X$.
  In particular, for the stated choice of $\epsilon$, 
  \begin{equation*}
    \det(A) \left(1 - \frac{\sqrt{e}}{2}\right) 
    \leq \det(A + \epsilon B) \leq
    \det(A) \left(1 + \frac{\sqrt{e}}{2}\right).
  \end{equation*}
\end{lem}
As $\rho\rightarrow 0$, note that
$g(t_{i,\rho}-s_j) = g(t_i-s_j)+\rho \cdot g'(t_i-s_j)+\frac{\rho^2}{2} g''(\xi)$,
for some $\xi \in [t_i-s_j,t_i-s_j+\rho]$. After applying this expansion, 
we subtract the rows with $g(t_i-s_j)$ from the rows with $g(t_{i,\rho}-s_j)$,
take $\rho^k$ outside of the determinant and we can write $M^{\rho}$ as:
\begin{align*}
  M^{\rho} = \rho^k \det(M_N + \rho M_P),
\end{align*}
$M_N$ is a matrix with entries independent of $\rho$ and with determinant:
\begin{equation}
  N = \det(M_N) = 
  \left|\begin{array}{cccc}
    f_{0} & g\left(s_{1}\right) & \cdots & g\left(s_{m}\right)\\
    0 & g\left(t_1-s_{1}\right) & \cdots & g\left(t_1-s_{m}\right)\\ 
    0 &  g'\left(t_{1}-s_{1}\right) & \cdots &  g'\left(t_{1}-s_{m}\right)\\
    \vdots  \\
    F\l(\tau_{\underline{l}} \r) & g\left(\tau_{\underline{l}}-s_{1}\right) & \cdots & g\left(\tau_{\underline{l}}-s_{m}\right)\\
    \vdots \\
    0 & g\left(t_k-s_{1}\right) & \cdots & g\left(\tau-s_{m}\right)\\
    0 &  g'\left(t_{k}-s_{1}\right) & \cdots &  g'\left(t_{k}-s_{m}\right)\\
    f_{1} & g\left(1-s_{1}\right) & \cdots & g\left(1-s_{m}\right)
  \end{array}\right|.
  \label{eq:gauss ex}
\end{equation}
Moreover, while the entries of $M_P$ depend on $\rho$, the magnitude of each entry can
be bounded from above independently of $\rho$. Consequently, $\norm{M_P}_F$ is bounded from above independently
of $\rho$. Let us assume for the moment that $N > 0$ and so $M_N$ is invertible.
Since $\rho(M_N^{-1} M_P) \leq \norm{M_N^{-1}}_2 \norm{M_P}_F$, this implies that
$\rho(M_N^{-1} M_P)$ is bounded from above independently of $\rho$.
We can then apply the stronger result of
Lemma \ref{lem:determinant bounds} to $M^{\rho}$ and obtain:
\begin{equation}
  0 < (1 - \rho C_N) \rho^k N \leq M^{\rho}
  \leq (1 + \rho C_N) \rho^k N,
  \label{eq:mrho lemma}
\end{equation}
where $C_N > 0$ is a constant that does not depend on $\rho$.
Note that we do not need write the condition on $\rho$ required by Lemma \ref{lem:determinant bounds}
explicitly because $\rho \to 0$ and also that \eqref{eq:mrho lemma} applies to the 
minors of $M^{\rho}$ and $N$ along the row containing $\tau_{\underline{l}}$.
That $N$ is indeed positive (and therefore we can apply Lemma \ref{lem:determinant bounds}) 
is established below
By its definition in \eqref{eq:f example}, $F(\tau_{\underline{l}})$ can take two values above. Either 
\begin{itemize}
\item $F(\tau_{\underline{l}})=0$,  which happens when there exists $i_0\in [k]$ such 
  that $\tau_{\underline{l}}\in T_{i_0,\epsilon}$, namely when $\tau_{\underline{l}}$ 
  is close to the support $T$.  In this case, by applying the Laplace expansion to $N$, 
  we find that  
  \begin{align}
    N = 
    f_{0}
    \cdot N_{1,1}
    + f_{1}  \cdot N_{m+1,1},
    \label{eq: special case G taul zero}
  \end{align} 
  where $N_{1,1}$ are $N_{m+1,1}$ are the corresponding minors in \eqref{eq:gauss ex}. 
  Note that  both $ N_{1,1}$ and $N_{m+1,1}$ are positive because the  Gaussian
  window is \emph{extended totally positive}, see Example 5 in \cite{karlin1966tchebycheff}. 
  Recalling that $f_{0},f_{1}>0$, we conclude that $N$ is positive. Therefore, when $\rho$ 
  is sufficiently small, \eqref{eq:gauss ex} implies that $M^\rho$ is non-negative 
  when $F(\tau_{\underline{l}})=0$. Or 
\item $F(\tau_{\underline{l}})=\bar{f}$, which happens when $\tau_{\underline{l}}\in T_\epsilon^C$, 
  namely when $\tau_{\underline{l}}$ is away from the support $T$.  Suppose that $f_{0} \gg \bar{f}$ 
  so that  $N$ is dominated by its first minor, namely  $N_{1,1}$. More precisely, by applying 
  the Laplace expansion to $N$ in \eqref{eq:gauss ex}, we find that 
  \begin{align}
    \label{eq:special case N 2}
    N = f_0\cdot N_{1,1}-\bar{f} \cdot N_{\underline{l},1} +f_1\cdot N_{m+1,1},
  \end{align} 
  in which all three minors are positive because the Gaussian window is extended totally positive. 
  Also, note that $N_{\underline{l},1}$ does \emph{not} depend on $\tau_{\underline{l}}$ and 
  recall also that $f_0,\bar{f},f_1$ are all positive.  Therefore $N$ in \eqref{eq:special case N 2} 
  is positive if 
  \begin{align}
    \frac{f_0}{\bar{f}}  > \frac{N_{\underline{l},1}}{\min_{\tau_{\underline{l}}} N_{1,1}},
    \label{eq:dominated}
  \end{align}
  where the minimum is over $\tau_{\underline{l}} \in T_{\epsilon}^C$.  The right-hand side 
  above is well-defined because $N_{1,1}= N_{1,1}(\tau_{\underline{l}})$ is positive for 
  every $\tau_{\underline{l}}\in I$, $N_{1,1}(\tau_{\underline{l}})$ is a continuous function of 
	$\tau_{\underline{l}}$, and $I$ is compact. Indeed, $N_{1,1}(\tau_{\underline{l}})$ 
  is positive because the Gaussian window is extended totally positive. 
  As before, $N$ being positive implies that  $ M^\rho  $ is non-negative when $\rho$ is sufficiently 
  small, see \eqref{eq:gauss ex}.  
\end{itemize}
By combining both cases above, we conclude that $M^{\rho}$ is non-negative for sufficiently 
small $\rho$  provided that \eqref{eq:dominated} holds,  
thereby verifying  part 1  of Definition \ref{def:(T-systems,-modified)-For}.  
To verify part 2 of that definition, consider the minors along 
the row containing  $\tau_{\underline{l}}$ in $M^\rho$, see \eqref{eq:gauss ex}. 
Starting with the first minor along this row
and applying the same arguments as before for $M^{\rho}$,
we observe, after applying Lemma \ref{lem:determinant bounds}, that
\begin{equation}
  M^\rho_{\underline{l},1}  
  \geq
  (1 - \rho C_{\underline{l},1}) \rho^k
  \l|\begin{array}{ccc}
    g\l(s_1 \r) & \cdots & g\l(s_m \r) \\
    g\l(t_1-s_1\r) & \cdots & g\l(t_1-s_m\r) \\
    g'\l(t_{1}-s_1\r) & \cdots &  g'\l(t_{1}-s_m\r)\\
    \vdots & & \vdots\\
    g\l(t_k-s_1\r) & \cdots & g\l(t_1-s_m\r) \\
    g'\l(t_{k}-s_1\r) & \cdots &  g'\l(t_{k}-s_m\r)\\
    g\l(1-s_1 \r) & \cdots & g\l(1-s_m \r) 
  \end{array}\r| 
  =: (1 - \rho C_{\underline{l},1}) \rho^k\cdot N_{\underline{l},1},
  \label{eq:1st minor}
\end{equation}
and also $M^\rho_{\underline{l},1} 
\leq (1 + \rho C_{\underline{l},1}) \rho^k\cdot N_{\underline{l},1}$ as $\rho\rightarrow0$. 
Here $C_{\underline{l},1} > 0$ is a constant that does not depend on $\rho$. Moreover, 
$N_{\underline{l},1}$ does not depend on $\rho$ and is positive because the Gaussian 
window is extended totally positive.
Therefore $M^{\rho}_{\underline{l},1}$ in \eqref{eq:1st minor}  
approaches zero at the rate $\rho^k$. 

Consider next the $(j+1)$th minor along the row 
containing $\tau_{\underline{l}}$ of $M^\rho$ in \eqref{eq:gauss ex},
namely $M_{\underline{l},j+1}^\rho$ with $j=1,\ldots,m$. Using the 
same arguments as before, we obtain after applying 
Lemma \ref{lem:determinant bounds} that
\begin{align}
  M_{\underline{l},j+1}^\rho &
  \geq (1 - \rho C_{\underline{l},j+1}) \rho^k
  \l|\begin{array}{ccccccc}
    f_{0} & g\l(s_1\r) & \cdots & g(s_{j-1}) & g(s_{j+1}) & \cdots & g\l(s_m \r) \\
    0 & g\l(t_1-s_1\r) & \cdots & g(t_1-s_{j-1}) & g(t_1-s_{j+1}) & \cdots & g\l(t_1-s_m \r) \\
    0 & g'\l(t_1-s_1\r) & \cdots & g'(t_1-s_{j-1}) & g'(t_1-s_{j+1}) & \cdots & g'\l(t_1-s_m \r) \\
    \vdots \\
    0 & g\l(t_k-s_1\r) & \cdots & g(t_k-s_{j-1}) & g(t_k-s_{j+1}) & \cdots & g\l(t_k-s_m \r) \\
    0 & g'\l(t_k-s_1\r) & \cdots & g'(t_k-s_{j-1}) & g'(t_k-s_{j+1}) & \cdots & g'\l(t_k-s_m \r) \\
    f_{1} & g\l(1-s_1\r) & \cdots & g(1-s_{j-1}) & g(1-s_{j+1}) & \cdots & g\l(1-s_m \r)
  \end{array}\r| 
  \nonumber\\
  & = 
  (1 - \rho C_{\underline{l},j+1}) \rho^k f_{0}
  \l| \begin{array}{cccccc}
    g\l(t_{1}-s_1\r) & \cdots & g\l(t_{1}-s_{j-1}\r) & g\l(t_{1}-s_{j+1}\r) & \cdots  & g\l(t_{1}-s_m\r) \\
    g'\l(t_{1}-s_1\r) & \cdots & g'\l(t_{1}-s_{j-1}\r) & g'\l(t_{1}-s_{j+1}\r) & \cdots  & g'\l(t_{1}-s_m\r) \\
    \vdots\\
    g\l(t_{k}-s_1\r) & \cdots & g\l(t_{k}-s_{j-1}\r) & g\l(t_{k}-s_{j+1}\r) & \cdots & g\l(t_{k}-s_m\r) \\
    g'\l(t_{k}-s_1\r) & \cdots & g'\l(t_{k}-s_{j-1}\r) & g'\l(t_{k}-s_{j+1}\r) & \cdots &  g'\l(t_{k}-s_m\r)\\
    g\l(1-s_1\r) & \cdots & g\l(1-s_{j-1}\r) & g\l(1-s_{j+1}\r) & \cdots & g\l(1-s_m\r) \\
  \end{array} \r| \nonumber\\
  & \qquad -(1 - \rho C_{\underline{l},j+1})\rho^k f_{1}
  \l| \begin{array}{cccccc}
    g\l(s_1\r) & \cdots & g\l(s_{j-1}\r) & g\l(s_{j+1}\r) & \cdots &  g\l(s_m\r) \\
    g\l(t_{1}-s_1\r) & \cdots & g\l(t_{1}-s_{j-1}\r) & g\l(t_{1}-s_{j+1}\r) & \cdots &  g\l(t_{1}-s_m\r) \\
    g'\l(t_{1}-s_1\r) & \cdots & g'\l(t_{1}-s_{j-1}\r) & g'\l(t_{1}-s_{j+1}\r) & \cdots &  g'\l(t_{1}-s_m\r) \\
    \vdots\\
    g\l(t_{k}-s_1\r) & \cdots & g\l(t_{k}-s_{j-1}\r) & g\l(t_{k}-s_{j+1}\r) & \cdots &  g\l(t_{k}-s_m\r) \\
    g'\l(t_{k}-s_1\r) & \cdots & g'\l(t_{k}-s_{j-1}\r) & g'\l(t_{k}-s_{j+1}\r) & \cdots  & g'\l(t_{k}-s_m\r)
  \end{array}\r|
  \nonumber\\
  &
  =: (1- \rho C_{\underline{l},j+1}) 
    \rho^k \l(f_0 \cdot N_{\underline{l},j+1,0} - f_1 \cdot N_{\underline{l},j+1,1}\r)
  =: (1 - \rho C_{\underline{l},j+1} ) \rho^k \cdot N_{\underline{l},j+1},
  \label{eq:minors}
\end{align}
and also $M_{\underline{l},j+1}^\rho \leq
(1 + \rho C_{\underline{l},j+1} ) \rho^k \cdot N_{\underline{l},j+1}$ 
as $\rho \rightarrow 0$, provided $N_{\underline{l},j+1} > 0$.
Here, $N_{\underline{l},j+1}$ is the determinant on the first line of \eqref{eq:minors}
and $C_{\underline{l},j+1} > 0$ is a constant independent of $\rho$.
Note that $N_{\underline{l},j+1,0}$ and $N_{\underline{l},j+1,1}$ are both 
positive because the Gaussian window is extended totally positive. 
To ensure $N_{\underline{l},j+1} > 0$, we require $f_{0}\gg f_{1}$, or more precisely, the 
following to hold.
\begin{equation}
  \frac{f_0}{f_1} > \frac{N_{\underline{l},j+1,1}}{N_{\underline{l},j+1,0}}.
  \label{eq:f0 n f1}
\end{equation}
It then follows that $M_{\underline{l},j+1}^\rho$ approaches zero at the 
rate $\rho^k$ for every $j$, thereby verifying part 2 in 
Definition \ref{def:(T-systems,-modified)-For} for $\tau_{\underline{l}} \in int(I)$. 
In conclusion, we find that $\{F\}\cup\{\phi_j\}_{j=1}^m$ form a T*-system on $I$ 
with $\phi_j(t) = g(t-s_j) = e^{-\frac{(t - s_j)^2}{\sigma^2}} $ as in \eqref{eq:gaussian window}, 
provided that  \eqref{eq:dominated} and \eqref{eq:f0 n f1} hold. 

To establish that $\{F^\pi\}\cup \{\phi_j\}_{j=1}^m$ form a T*-system on $I$ with the Gaussian window,
we note that $F(\tau_{\underline{l}})$ is replaced by by $F^{\pi}(\tau_{\underline{l}})$, which takes
values $\pm 1$ when $\tau_{\underline{l}} \in T_{i,\epsilon}$ for some $i$, 
as indicated in \eqref{eq:g pi}.
The previous argument, thus, goes through similarly, showing that $\{F^\pi\}\cup \{\phi_j\}_{j=1}^m$ is  a T*-system on $I$ for arbitrary sign pattern $\pi$ when  \eqref{eq:dominated} and \eqref{eq:f0 n f1} hold 
and $f_{0} \gg 1$. 
The extra condition $f_{0} \gg 1$ comes from the only difference between the two proofs, 
namely that, instead of \eqref{eq: special case G taul zero}, we have
\begin{align}
  N = 
  f_{0} \cdot N_{1,1}
  \pm 1 \cdot N_{\underline{l},1}
  + f_{1}  \cdot N_{m+1,1}.
  \label{eq:special case N 2 pi}
\end{align}
Therefore to ensure $N > 0$, we now require the following to hold.
\begin{align} \label{eq:tstar_cond_f0}
  f_0 > \frac{N_{\underline{l},1}}{\min_{\tau_{\underline{l}}} N_{1,1}}.
\end{align}
We leave out the mostly repetitive details. 
Hence if \eqref{eq:dominated}, \eqref{eq:f0 n f1} and \eqref{eq:tstar_cond_f0} hold, 
then $\{F^\pi\}\cup \{\phi_j\}_{j=1}^m$ form a T*-system on $I$ with the Gaussian window.
This completes the proof of Lemma \ref{lem:Gaussian is Tstar sys}. 

\textit{Remark 1.} While we do not require that $f_{1} \gg \bar{f}$, if we impose
that both $f_0 \gg \bar{f}$ and $f_{1} \gg \bar{f}$, then \eqref{eq:special case N 2} 
holds for smaller $\frac{f_{0}}{\bar{f}}$, so it is useful in practice. 
Similarly, from \eqref{eq:special case N 2 pi} we want that also $f_1 \gg 1$.

\textit{Remark 2.} From this proof and in light of the first remark, we see that we only 
need to specify one end point rather than both, so we only need $m=2k+1$. 
However, the dual polynomial is better in practice if we have conditions at both 
end points (if we specify both $f_0$ and $f_1$).

\section{Proof of Lemma \ref{lem:bounds on b for Gaussian} \label{sec:proof of lemma bounds on b for Gaussian}}

Recall our assumptions that 
\begin{equation}
  m=2k+2,\qquad s_1 = 0, \quad s_m = s_{2k+2}=1,
\label{eq:begin n end points}
\end{equation}
and that 
\begin{equation}
  | s_{2i} -t_i | \le \eta, \qquad s_{2i+1}-s_{2i} =\eta,
  \qquad \forall i\in [k],
  \label{eq:assump}
\end{equation} 
That is, we collect two 
samples near  each impulse in $x$, supported on $T$. 
In addition, we make the following assumptions on $\eta$ and $\sigma$:
\begin{equation}
  \sigma \leq \sqrt{2}, \quad 
  \Delta > \sigma \sqrt{\log\left(3 + \frac{4}{\sigma^2} \right)}, 
  \quad \eta \leq \sigma^2.
	\label{eq:simplifying assumptions}
\end{equation}
After studying Appendix \ref{sec:Proof-of-Proposition dual construction}, it becomes clear that the entries of $b\in\mathbb{R}^m$ are specified as 
\begin{equation}
  b_j = \lim_{\rho\rightarrow 0}
  (-1)^{j+1}\frac{M^\rho_{\underline{l},j+1}}{M^{\rho}_{\underline{l},1}}
	= (-1)^{j+1} \lim_{\rho\rightarrow 0} 
    \frac{M^\rho_{\underline{l},j+1}}{M^{\rho}_{\underline{l},1}},
  \qquad j\in[m],
  \label{eq:reminder of coeffs}
\end{equation}
where the numerator and the denominator are the minors $\{M^\rho_{\underline{l},j}\}_{j=1}^{m+1}$ 
of $M^{\rho}$ in \eqref{eq:m rho} along the row containing $\tau_{\underline{l}}$. 
Using the upper and lower bounds on these quantities derived earlier, we obtain:
\begin{equation}
  \frac{
    (1 - \rho C_{\underline{l},j+1}) N_{\underline{l},j+1}
  }{
   (1 + \rho C_{\underline{l},1}) N_{\underline{l},1}
  }
  \leq
  \frac{M^\rho_{\underline{l},j+1}}{M^{\rho}_{\underline{l},1}} 
  \leq
  \frac{
    (1 + \rho C_{\underline{l},j+1}) N_{\underline{l},j+1}
  }{
   (1 - \rho C_{\underline{l},1}) N_{\underline{l},1}
  }, 
  \qquad j\in[m],
  \label{eq:ratio bounds rho}
\end{equation}
which in turn implies the following expression for $b_j$:
\begin{equation}
  b_j =  (-1)^{j+1} \frac{
    N_{\underline{l},j+1}
  }{
    N_{\underline{l},1}
  }, \qquad j \in [m].
  \label{eq:ratio N}
\end{equation}
Recall that $N_{\underline{l},1}, N_{\underline{l},j+1} > 0$ and so,
$|b_j| = \frac{N_{\underline{l},j+1}}{N_{\underline{l},1}}$ for $j \in [m]$.
Therefore, in order to upper bound each $|b_j|$, we will respectively lower
bound the denominator $N_{\underline{l},1}$ and upper bound the numerators
$N_{\underline{l},j+1}$.

\subsection{Bound on the Denominator of \eqref{eq:ratio N}\label{sec:bnd on denom}}

We now find a lower bound for the first minor, 
namely $N_{\underline{l},1}$ in \eqref{eq:ratio N}. 
Let us conveniently assume that the spike locations $T=\{t_i\}_{i=1}^k$ and the sampling 
points $S=\{s_j\}_{j=2}^{m-1}$ are away from the boundary of interval $I=[0,1]$, namely 
\begin{equation*}
  \sigma\sqrt{ \log(1/\eta^3)} \le t_i \le 1-\sigma\sqrt{\log(1/\eta^3)},
  \qquad \forall i\in [k],
\end{equation*}
\begin{equation}
  \sigma\sqrt{ \log(1/\eta^3)} \le s_j \le 1-\sigma\sqrt{ \log(1/\eta^3)},
  \qquad \forall j\in [2:m-1].
  \label{eq:away frm bndr}
\end{equation}
In particular, \eqref{eq:away frm bndr} implies that 
\begin{equation*}
  g(t_i) \leq \eta^3, \quad g(1-t_i) \leq \eta^3, \qquad i\in [k], 
\end{equation*}
\begin{equation}
  g(s_j) \leq \eta^3, \quad g(1-s_j) \leq \eta^3, \qquad j\in [2:m-1].
  \label{eq:small bndr}
\end{equation}
For the derivatives, we have that:
\begin{equation*}
  |g'(t_i)| = \frac{2t_i}{\sigma^2} g(t_i) 
  \leq \frac{2\eta^3}{\sigma^2},
\end{equation*}
where we used the fact that $0 \leq t_i \leq 1$ and \eqref{eq:small bndr}.
Similarly, for the $1-t_i$, $s_j$ and $1-s_j$:
\begin{equation*}
  |g'(t_i)| \leq \frac{2\eta^3}{\sigma^2}, \quad 
  |g'(1-t_i)| \leq \frac{2\eta^3}{\sigma^2}, \qquad i\in [k], 
\end{equation*}
\begin{equation}
  |g'(s_j)| \leq \frac{2\eta^3}{\sigma^2}, \quad 
  |g'(1-s_j)| \leq \frac{2\eta^3}{\sigma^2}, \qquad j\in [2:m-1].
  \label{eq:small bndr deriv}
\end{equation}
With the assumptions in \eqref{eq:begin n end points}, \eqref{eq:assump}, \eqref{eq:small bndr}, \eqref{eq:small bndr deriv} 
and the fact that $g(0)=1$,
we have that the determinant $N_{\underline{l},1}$ in \eqref{eq:1st minor} is equal to
\begin{equation}
  N_{\underline{l},1} = 
  \l|\begin{array}{cccccc}
    1 & \cdots & O(\eta^3) & O(\eta^3) & \cdots & O(\eta^3) \\
    \vdots & & \vdots & \vdots & & \vdots \\
    O(\eta^3) & \cdots & g(t_i - s_{2u}) &  g(t_i - s_{2u+1}) & \cdots & O(\eta^3) \\
    O(\eta^3/\sigma^2)  & \cdots & g'(t_{i}-s_{2u}) & g'(t_i - s_{2u+1}) & \cdots & O(\eta^3/\sigma^2) \\
    \vdots & & \vdots & \vdots & & \vdots \\
    O(\eta^3) &  \cdots & O(\eta^3) & O(\eta^3) & \cdots & 1
  \end{array}\r|,
\end{equation}
where we wrote
\begin{equation}
  g(t_i) = O(\eta^3)
  \text{ and } 
  g'(t_i) = O(\frac{\eta^3}{\sigma^2})
  \iff
  |g(t_i)| \leq M_1 \eta^3
  \text{ and }
  |g'(t_i)| \leq M_2 \frac{\eta^3}{\sigma^2},
  \label{eq:big o notation}
\end{equation}
for $t_i$ within the bounds defined in \eqref{eq:away frm bndr}
and some $M_1,M_2 > 0$. Here, we can take $M_1=M_2=2$
and we use the same notation 
for $1-t_i$, $s_j$ and $1-s_j$ with the same constants $M_1=M_2=2$.
We then take the Taylor expansion of $g(t_i-s_{2u+1})$ and $g'(t_i-s_{2u+1})$
around $t_i - s_{2u}$, subtract the columns with $t_i - s_{2u}$ from
the columns where we performed the expansion, take $\eta^k$ outside
of the determinant and we obtain\footnote{
  Note that the equality in \eqref{eq:ncc} is due to the
  fact that $s_{2u+1}-s_{2u}=\eta$. If we relax this 
  condition to \eqref{eq:C1etaC2eta}, where $C_1\eta$ and $C_2\eta$ take 
  the roles of lower and maximum upper bounds of $\eta$, we obtain:
  \begin{equation*}
    C_1^k \eta^k | C + \bar{C'} |
    \leq N_{\underline{l},1} \leq
    C_2^k \eta^k | C + \bar{C'} |,
  \end{equation*}
  where $\bar{C'}$ is the same as $C'$
  in \eqref{eq:C' def} but every entry 
  multiplied by $\eta$ or 
  $s_{2u+1}-s_{2u} \in [C_1\eta,C_2\eta]$,
  so the order of its entries is the same
  as in $\eta C'$ in \eqref{eq:ncc}.
  The computations in this subsection will then follow 
  in a similar way 
  except that we will work with the lower bound
  involving $C_1^k \eta^k$ instead of the term 
  involving $\eta^k$.
  \label{fn:C1etaC2_1}  
}:
\begin{equation}
  N_{\underline{l},1} = 
  \eta^k | C + \eta C' |,
  \label{eq:ncc}
\end{equation}
where
\begin{equation}
  C = 
  \begin{bmatrix}
    1  & \cdots & 0 & 0 & \cdots & 0 \\
    \vdots &    & \vdots & \vdots & & \vdots \\
    0  & \cdots & g(t_i-s_{2u}) & - g'(t_i-s_{2u}) & \cdots & 0 \\
    0  & \cdots & g'(t_i-s_{2u}) & - g''(t_i-s_{2u}) & \cdots & 0 \\
    \vdots &    & \vdots & \vdots & & \vdots \\
    0  & \cdots & 0 & 0 & \cdots & 1
  \end{bmatrix}
\end{equation}
\begin{equation}
  C' = 
  \begin{bmatrix}
    0  & \cdots & O(\eta^2) & O(2\eta) & \cdots & O(\eta^2) \\
    \vdots &    & \vdots & \vdots & & \vdots \\
    O(\eta^2)  & \cdots & 0 & \frac12 g''(\xi_{i,u}) & \cdots & O(\eta^2) \\
    O(\eta^2/\sigma^2)  & \cdots & 0 & \frac12 g'''(\xi'_{i,u}) & \cdots & O(\eta^2/\sigma^2) \\
    \vdots &    & \vdots & \vdots & & \vdots \\
    O(\eta^2)  & \cdots & O(\eta^2) & O(2\eta) & \cdots & 0
  \end{bmatrix},
  \label{eq:C' def}
\end{equation}
for some $\xi_{i,u}, \xi'_{i,u} \in [t_i-s_{2u}-\eta, t_i-s_{2u}]$
for all $i,u = 1,\ldots,k$.
Note that, using the notation in \eqref{eq:big o notation}, if we subtract a function
that is $O(\eta)$ with constant $M_1>0$ from another function that
is $O(\eta)$ with constant $M_2>0$, we obtain a function that
is $O(\eta)$ with constant $M_1 + M_2$, which is why we wrote $O(2\eta)$ on the
first and last columns where we subtracted two functions $O(\eta)$ with the same 
constant $M_1$ (so $O(2\eta)$ implies $\leq 2M_1 \eta$).
Next, we apply Taylor expansion around $t_i - t_u$ in the terms with $t_i - s_{2u}$
in $C$ as follows:
\begin{equation}
 g(t_i - s_{2u}) = g(t_i - t_u + t_u - s_{2u}) 
 = g(t_i - t_u) + (t_u - s_{2u}) g'(\xi^*_{i,u}),
\end{equation}
for some $\xi^*_{i,u} \in [t_i-t_u - |t_u-s_{2u}|, t_i-t_u + |t_u-s_{2u}|]$
for all $i,u = 1,\ldots,k$.
Note that $|t_u - s_{2u}| \leq \eta$ according to \eqref{eq:assump}. 
By applying a similar Taylor expansion to $g'$ and $g''$, we can write
\begin{equation}
  C = A + \eta A',
  \label{eq:caa}
\end{equation}
where
\begin{equation}
  A = \begin{bmatrix}
    1  & \cdots & 0 & 0 & \cdots & 0 \\
    \vdots &    & \vdots & \vdots & & \vdots \\
    0  & \cdots & g(t_i-t_u) & - g'(t_i-t_u) & \cdots & 0 \\
    0  & \cdots & g'(t_i-t_u) & - g''(t_i-t_u) & \cdots & 0 \\
    \vdots &    & \vdots & \vdots & & \vdots \\
    0  & \cdots & 0 & 0 & \cdots & 1
  \end{bmatrix}
\end{equation}
and
\begin{equation}
  A' = \begin{bmatrix}
    0  & \cdots & 0 & 0 & \cdots & 0 \\
    \vdots &    & \vdots & \vdots & & \vdots \\
    0  & \cdots & \frac{t_u-s_{2u}}{\eta} g'(\xi^*_{i,u}) 
        & - \frac{t_u-s_{2u}}{\eta} g''(\xi^{*\prime}_{i,u}) & \cdots & 0 \\
    0  & \cdots &  \frac{t_u-s_{2u}}{\eta}  g''(\xi^{*\prime}_{i,u}) 
        & - \frac{t_u-s_{2u}}{\eta} g'''(\xi^{*\prime\prime}_{i,u}) & \cdots & 0 \\
    \vdots &    & \vdots & \vdots & & \vdots \\
    0  & \cdots & 0 & 0 & \cdots & 0
  \end{bmatrix}
\end{equation}
for some $\xi^*_{i,u},\xi^{*\prime}_{i,u},\xi^{*\prime\prime}_{i,u} \in [t_i-t_u - |t_u-s_{2u}|, t_i-t_u + |t_u-s_{2u}|]$
for all $i,u = 1,\ldots,k$. We now substitute \eqref{eq:caa} into \eqref{eq:ncc}
and we obtain:
\begin{equation}
  N_{\underline{l},1} = \eta^k|A + \eta(A'+C')|,
\end{equation}
Assuming for the moment that $|A| > 0$ holds, we obtain via
Lemma \ref{lem:determinant bounds} the following bound:
\begin{equation}
  \eta^k \left(1 - \frac{\sqrt{e}}{2} \right) \det(A)
  \leq N_{\underline{l},1} \leq
  \eta^k \left(1 + \frac{\sqrt{e}}{2} \right) \det(A)
  \label{eq:bound denom}
\end{equation}
if
\begin{equation}
  \eta \leq \frac{8}{34(2k+2)\rho(A^{-1}(A' + C'))}.
  \label{eq:cond eta unfinished}
\end{equation}
We look closer at the condition \eqref{eq:cond eta unfinished} on $\eta$
in Section \ref{sec:cond on eta}.
That $|A| > 0$ (and therefore our application of Lemma \ref{lem:determinant bounds} 
above is valid) is established below. 
We can in fact write $A$ more compactly as follows. For scalar $t$, let 
\begin{equation}
  H(t) := 
  \l[
  \begin{array}{cc}
    g(t) & -g'(t)\\
    g'(t) & -g''(t)
  \end{array}
  \r]\in\mathbb{R}^{2\times 2},
\end{equation}
which allows us to rewrite $A$ as 
\begin{equation}
  \label{eq:def of A}
  A = \l[
  \begin{array}{ccc}
    1 & 0_{1\times 2k} & 0\\
    0_{2k \times 1} & B & 0_{2k \times 1} \\
    0 & 0_{1\times 2k} & 1
  \end{array}
  \r] \in\mathbb{R}^{m \times m},
\end{equation}
\begin{equation}
  B := \l[
  \begin{array}{ccccc}
    H(0) & H(t_1-t_2) & H(t_1-t_3) & \cdots & H(t_1-t_k) \\
    H(t_2-t_1) & H(0) & H(t_2-t_3) & \cdots & H(t_2-t_k) \\
    \vdots \\
    H(t_k-t_1) & H(t_k-t_2) & H(t_k-t_3) & \cdots & H(0)
  \end{array}
  \r] \in \mathbb{R}^{2k\times 2k},
  \label{eq:def of B}
\end{equation}
where $0_{a\times b}$ is the matrix of zeros of size $a\times b$.
It follows from \eqref{eq:def of A} 
that $|A| = |B|$ by Laplace expansion of $|A|$. 
In particular, the eigenvalues of $A$ are: $1,1$ and the eigenvalues of $B$.
Let us note that $B$ 
is a symmetric matrix\footnote{
  Indeed, $g, g''$  are even functions and $g'$ is an odd function.
}, since $H(-t) = H(t)^T$; hence, $A$ is also symmetric.
We now proceed to lower bound $|B|$, the details of which are given 
in Appendix \ref{sec:proof of lwr bnd on B}. The main observation is 
that $H(0)$ is a diagonal matrix while the entries of $H(t_i-t_j)$ for $i \neq j$
decay with the separation of sources $\Delta$.
\begin{lem}\label{lem:lwr bnd on B}
 Let $\sigma \leq \sqrt{2}$,
 $\Delta > \sigma \sqrt{\log{\left( 3 + \frac{4}{\sigma^2} \right)}}$
 and
  \begin{align*}
    0 < \quad
    F_{\min}\left(\Delta,\frac{1}{\sigma}\right) = 
      1 - \left(1 + \frac{2}{\sigma^2} \right)
      \frac{2 e^{-\frac{\Delta^2}{\sigma^2}}}{1 - e^{-\frac{\Delta^2}{\sigma^2}}}
    \quad < 1.
  \end{align*}
  Then, for each $i = 1,\dots, 2k$, it holds that
  \begin{align*}
    \lambda_i(B) \geq 
    F_{\min}\left(\Delta,\frac{1}{\sigma}\right).
  \end{align*}
  \end{lem}
Since $|A| = |B|$, we obtain via Lemma \ref{lem:lwr bnd on B} that
$|A| \geq F_{\min}\left(\Delta,\frac{1}{\sigma}\right)^{2k} > 0$. Using
this in \eqref{eq:bound denom}, leads to the following bound:
\begin{align}
  N_{\underline{l},1} 
  & \geq \eta^k \left(1 - \frac{\sqrt{e}}{2}\right)
    F_{\min}\left(\Delta,\frac{1}{\sigma}\right)^{2k}.
    \label{eq:final bound denominator}
\end{align}

\subsection{Bound on the Numerator of \eqref{eq:ratio N} 
  \label{sec:bnd on num}}

Since \eqref{eq:mrho lemma} holds for all the minors $M^{\rho}_{\underline{l},j}$ 
and $N_{\underline{l},j}$, let us now upper bound the $N_{\underline{l},j}$ 
for $j = 2,\ldots,m+1$ in the numerator of \eqref{eq:ratio N}. Note that
we distinguish two cases: $j \in \{3,\ldots,m\}$ and $j \in \{2,m+1\}$.
To simplify the presentation for the first case, suppose, for example, that  $j=3$.  
Using the assumptions in \eqref{eq:begin n end points}, \eqref{eq:assump}, 
\eqref{eq:small bndr}, \eqref{eq:small bndr deriv} and the fact that $g(0)=1$, 
\begin{align}
  N_{\underline{l},3} 
  &=
  \begin{vmatrix}
    f_{0} & 1  & O(\eta^3) & \cdots & O(\eta^3) & O(\eta^3) & \cdots & O(\eta^3) \\
    \vdots \\
    0 & O(\eta^3)  & g(t_i-s_{3}) & \cdots & g(t_i-s_{2u}) & g(t_i-s_{2u+1}) & \cdots & O(\eta^3) \\
    0 & O(\eta^3/\sigma^2)  & g'(t_i-s_{3}) & \cdots & g'(t_i-s_{2u}) & g'(t_i-s_{2u+1}) & \cdots & O(\eta^3/\sigma^2) \\
    \vdots \\
    f_{1} & O(\eta^3)  & O(\eta^3) & \cdots & O(\eta^3) & O(\eta^3) & \cdots & 1
  \end{vmatrix}.
  \label{eq:det n3}
\end{align}
We now expand $g(t_i - s_{2u+1})$ around $g(t_i-s_{2u})$, subtract the columns and take
$\eta$ out of the determinant as before\footnote{
  Similarly to the issue addressed in footnote \ref{fn:C1etaC2_1},
  if instead of $s_{2u+1}-s_{2u}=\eta$ 
  we have \eqref{eq:C1etaC2eta}, then:
  \begin{equation*}
    C_1^{k-1} \eta^{k-1} \det(\tilde{N}_3)
    \leq N_{\underline{l},3} \leq
    C_2^{k-1} \eta^{k-1} \det(\tilde{N}_3),
  \end{equation*}
  and then the proof will continue in a similar way except
  that we work with the upper bound 
  of $N_{\underline{l},3}$ above and some of the calculations
  will involve $C_1$ and $C_2$ as well.
  \label{fn:C1etaC2_2}  
}:
\begin{align}
  N_{\underline{l},3} 
  &= \eta^{k-1}
  \begin{vmatrix}
    f_{0} & 1  & O(\eta^3) & \cdots & O(\eta^3) & O(2\eta^2) & \cdots & O(\eta^3) \\
    \vdots \\
    0 & O(\eta^3)  & g(t_i-s_{3}) & \cdots & g(t_i-s_{2u}) 
      & -g'(\xi_{i,u}) & \cdots & O(\eta^3) \\
    0 & O(\eta^3/\sigma^2)  & g'(t_i-s_{3}) & \cdots & g'(t_i-s_{2u}) 
      & -g''(\xi'_{i,u}) & \cdots & O(\eta^3/\sigma^2) \\
    \vdots \\
    f_{1} & O(\eta^3)  & O(\eta^3) & \cdots & O(\eta^3) & O(2\eta^2) & \cdots & 1
  \end{vmatrix}
  =: \eta^{k-1} \det(\tilde{N}_3),
  \label{eq:det tilde n}
\end{align}
where $\xi_{i,u},\xi'_{i,u} \in [t_i - s_{2u} - \eta, t_i - s_{2u}]$
and denote the matrix in $\eqref{eq:det tilde n}$ by $\tilde{N}_3$. 
Note the $\eta^{k-1}$ since we only perform column operations on $k-1$ columns 
and $\det(\tilde{N}_3) > 0$ (due to the choice of $f_0$ and $f_1$).
Next let $C\in\mathbb{R}^{m\times 3}$ consist of the first, second, 
and last columns of $\tilde{N}_3$ and  $D\in\mathbb{R}^{m\times (m-3)}$ consist 
of the rest of the columns of $\tilde{N}_3$. Then, we may write that 
\begin{equation}
  \det([C D])^{2} = 
  \left(
    \begin{bmatrix}
      C^* \\ D^*
    \end{bmatrix}
    [C D]
  \right)
  = \det \left( \begin{bmatrix}
    C^* C & C^* D \\
    D^*C  & D^* D
  \end{bmatrix} \right)
  \leq
  \det(C^* C) \det(D^* D),
\end{equation}
where we note that swapping columns only changes the sign in a determinant 
(here, $\det([CD]) = -\det(\tilde{N}_3)$) and in the last 
inequality we applied Fischer's inequality (see, for example, Theorem 7.8.3 in \cite{Horn1990}), 
which works because the matrix $[C D]^*[CD]$
is Hermitian positive definite.
Therefore, we have that
\begin{equation}
  \det(\tilde{N}_3) = |\det([CD])| \leq \det(C^*C)^\frac12 \det(D^*D)^\frac12,
  \label{eq:C and D}
\end{equation} 
and it suffices to bound the determinants on the right-hand side above.

\subsubsection{Bounding $\det(C^*C)$}

We now write $C$ as follows:
\begin{equation}
  C = \begin{bmatrix}
    f_0 & 1 & 0 \\
    \vdots & \vdots & \vdots \\
    0 & 0 & 0 \\
    0 & 0 & 0 \\
    \vdots & \vdots & \vdots \\
    f_1 & 0 & 1
  \end{bmatrix} + \begin{bmatrix}
    0 & 0 & O(\eta^3) \\
    \vdots & \vdots & \vdots \\
    0 & O(\eta^3) & O(\eta^3) \\
    0 & O(\eta^3/\sigma^2) & O(\eta^3/\sigma^2) \\
    \vdots & \vdots & \vdots \\
    0 & O(\eta^3) & 0
  \end{bmatrix}
  =: X + \pert,
\end{equation}
where we denote the first matrix by $X$ and the 
second matrix by $\pert$. We have that
\begin{equation}
  \det(C^*C) = \det((X + \pert)^* (X+\pert)) 
  = \det(X^*X + \pert'), 
\end{equation}
with $\pert' = X^*\pert + \pert^*X + \pert^*\pert$ 
and we apply Weyl's inequality to $C^*C$ to obtain:
\begin{equation}
  \det(C^*C) \leq 
    (\lambda_1(X^*X) + \lambda_{\max}(\pert'))
    (\lambda_2(X^*X) + \lambda_{\max}(\pert'))
    (\lambda_3(X^*X) + \lambda_{\max}(\pert')).
  \label{eq:weyl cc}
\end{equation}
Then
\begin{equation}
  X^*X = \begin{bmatrix}
    f_0^2 + f_1^2 & f_0 & f_1 \\
    f_0 & 1 & 0 \\
    f_1 & 0 & 1
  \end{bmatrix}
  \quad \text{ with } \quad
  \lambda_3(X^*X) = 0,
  \label{eq:xx eigs}
\end{equation}
and
\begin{equation}
  \|X^*X\|_F \leq \bar{C}(f_0,f_1), 
  \quad \text{ where } \quad
  \bar{C}(f_0,f_1) = f_0^2 + f_1^2 + 2f_0 + 2f_1 + 2,
\end{equation}
so
\begin{equation}
  \| X \|_2^2 = \| X^* X \|_2 
  \leq \|X^*X\|_F 
  \leq \bar{C}(f_0,f_1)
  \quad \text{ so } \quad
  \|X\|_2 = \sqrt{\bar{C}(f_0,f_1)}
  \label{eq:x norm}
\end{equation}
and 
\begin{align}
  \| \pert'\|_2 
  &= \|X^*\pert + \pert^*X + \pert^*\pert\|_2 \nonumber \\
  &\leq 2 \|X\|_2 \|\pert\|_2 + \|\pert\|_2^2 
    = \left(2 \sqrt{\bar{C}(f_0,f_1)} + \|\pert\|_2\right) \|\pert\|_2 
  \nonumber \\
  &\leq 3\sqrt{\bar{C}(f_0,f_1)} \|\pert\|_2,
  \label{eq:bound delta prime}
\end{align}
where the last inequality holds if
\begin{equation}
  \|\pert\|_2 \leq \sqrt{\bar{C}(f_0,f_1)}. 
  \label{eq:cond delta norm}
\end{equation}
Noting that $\pert'$ is symmetric 
and $\lambda_{\max}(\pert') \leq \|\pert'\|_2$,
we substitute \eqref{eq:xx eigs} and \eqref{eq:bound delta prime} 
into \eqref{eq:weyl cc} and obtain:
\begin{equation}
  \det(C^*C) \leq 
  \left(\bar{C}(f_0,f_1)+ 3\sqrt{\bar{C}(f_0,f_1)} \|\pert\|_2 \right)^2
  3\sqrt{\bar{C}(f_0,f_1)} \|\pert\|_2,
\end{equation}
if \eqref{eq:cond delta norm} holds. Further applying \eqref{eq:cond delta norm}
in the parentheses, we obtain:
\begin{equation}
  \det(C^*C) \leq
  48 \bar{C}(f_0,f_1)^\frac52
  \|\pert\|_2.
  \label{eq:det cc bound unfinished}
\end{equation}
Now, using \eqref{eq:big o notation} with $M_1=M_2=2$, we are able to 
upper bound $\|\pert\|_F$,
which is also an upper bound for $\|\pert\|_2$:
\begin{align}
  \|\pert\|_2 \leq \|\pert\|_F 
  &\leq \sqrt{(2k+2)(2\eta^3)^2 + 2k \left(\frac{2\eta^3}{\sigma^2}\right)^2}
  \nonumber \\
  &\leq (2k+2)2\eta^3 + 2k \frac{2\eta^3}{\sigma^2}
  = \eta^3 \left(4k+4+\frac{4k}{\sigma^2}\right).
  \label{eq:bound delta 2}
\end{align}
Therefore, to satisfy \eqref{eq:cond delta norm}, it is sufficient
to find $\eta$ such that:
\begin{equation}
  \eta^3 \left(4k+4 + \frac{4k}{\sigma^2} \right) 
  \leq \sqrt{\bar{C}(f_0,f_1)}.
  \label{eq:cond eta f}
\end{equation}
With this choice of $\eta$, by substituting \eqref{eq:bound delta 2}
into \eqref{eq:det cc bound unfinished}, we obtain:
\begin{equation}
  \det(C^*C) \leq
  48 \bar{C}(f_0,f_1)^\frac52
  \left(4k + 4 + \frac{4k}{\sigma^2}\right) 
  \eta^3.
  \label{eq:bound det cc}
\end{equation}

Note that all our calculations so far will be used for bounding $\|b\|_2$. In case
of $\|b^{\pi}\|_2$, everything is the same except that we have:
\begin{equation}
  X_{\pi}^* X_{\pi} = 
  \begin{bmatrix}
    f_0^2 +f_1^2 + 2k & f_0 & f_1 \\
    f_0 & 1 & 0 \\
    f_1 & 0 & 1
  \end{bmatrix},
\end{equation}
which does not have a zero eigenvalue, so we will not obtain the $\eta^3$ factor. 
We omit the calculations, but we obtain:
\begin{equation}
  \det(C_{\pi}^* C_{\pi}) \leq
  64 \left(\bar{C}(f_0,f_1) + 2k\right)^3
  \label{eq:detcc pi}
\end{equation}
with a condition on $\eta$ that is weaker than \eqref{eq:cond eta f}.

\subsubsection{Bounding $\det(D^*D)$}

Then, since $D^*D$ is Hermitian positive definite, 
we can apply Hadamard's inequality (Theorem 7.8.1 in \cite{Horn1990}) 
to bound its determinant by the product of its 
main diagonal entries (i.e. the squared 2-norms of the columns of $D$), so we obtain, 
after we use \eqref{eq:big o notation} with $M_1=M_2=2$:
\begin{align}
  \det(D^* D) & \leq
  \left( 8\eta^6
    + \sum_{i=1}^{k} g(t_i - s_3)^2 
    + \sum_{i=1}^{k} g'(t_i - s_3)^2 
  \right) \nonumber \\
  &\quad \cdot \prod_{u=2}^{k} \left( 8\eta^6
    + \sum_{i=1}^{k} g(t_i - s_{2u})^2 
    + \sum_{i=1}^{k} g'(t_i - s_{2u})^2 
  \right) \nonumber \\
  &\quad \cdot 
  \prod_{u=2}^{k} \left( 32\eta^4
    + \sum_{i=1}^{k} g'(\xi_{i,u})^2 
    + \sum_{i=1}^{k} g''(\xi'_{i,u})^2 
  \right),
  \label{eq:num 0}
\end{align}
where $\xi_{i,u},\xi'_{i,u} \in [t_i - s_{2u}-\eta, t_i-s_{tu}]$.

For fixed $u \in \{1,\ldots,k\}$, we may write that 
\begin{align}
  \sum_{i=1}^k g(t_i-s_{2u})^2 
  &\leq 2 \sum_{i=0}^{\infty} g(i\Delta)^2
  \leq 2 \sum_{i=0}^{\infty} g(i\Delta)
  = 2 \sum_{i=0}^{\infty} e^{-\frac{i^2 \Delta^2}{\sigma^2}}
  \nonumber \\
  &\leq 2 \sum_{i=0}^{\infty} \left(e^{-\frac{\Delta^2}{\sigma^2}}\right)^i
  = \frac{2}{1-e^{-\frac{\Delta^2}{\sigma^2}}}.
  \label{eq:two sums}
\end{align}
To see the first inequality above, 
if we note that $s_{2u} \leq t_u \leq s_{2u+1}$, we have that

\begin{align*}
  g(t_{u-1}-s_{2u}) \leq g(0),
    &\quad \quad g(t_u-s_{2u}) \leq g(0),
    \\
  g(t_{u-2}-s_{2u}) \leq g(\Delta),
    &\quad \quad g(t_{u+1}-s_{2u}) \leq g(\Delta),
    \\
  \vdots \quad&\quad\quad\quad\quad \vdots
\end{align*}
and by adding the inequalities we obtain \eqref{eq:two sums}.
Likewise, it holds that 
\begin{equation*}
  \sum_{i=1}^{k} g'(t_i-s_{2u})^2 
  \leq \frac{4}{\sigma^4} \cdot \frac{2}{1-e^{-\frac{\Delta^2}{\sigma^2}}},
\end{equation*}
\begin{equation}
  \sum_{i=1}^{k} g''(t_i-s_{2u})^2 
  \leq \left(\frac{2}{\sigma^2} + \frac{4}{\sigma^4}  \right)^2
  \cdot \frac{2}{1-e^{-\frac{\Delta^2}{\sigma^2}}}.
  \label{eq:two sums prime}
\end{equation}
and, for $\xi_{i,u}, \xi'_{u,u} \in [t_i-s_{2u}-\eta,t_i-s_{2u}]$:

\begin{equation*}
  \sum_{i=1}^{k} g'(\xi_{i,u})^2 
  \leq \frac{4}{\sigma^4} \cdot \frac{2}{1-e^{-\frac{\Delta^2}{\sigma^2}}},
\end{equation*}
\begin{equation}
  \sum_{i=1}^{k} g''(\xi'_{i,u})^2 
  \leq \left(\frac{2}{\sigma^2} + \frac{4}{\sigma^4}  \right)^2
  \cdot \frac{2}{1-e^{-\frac{\Delta^2}{\sigma^2}}}.
  \label{eq:two sums xi}
\end{equation}
Note that we obtained \eqref{eq:two sums xi} in the same way as \eqref{eq:two sums}:
\begin{align*}
  g(\xi_{u,u}) \leq g(0)
    &\quad \quad g(\xi_{u-1,u}) \leq g(0) 
    \\
  g(\xi_{u+1,u}) \leq g(\Delta)
    &\quad \quad g(\xi_{u-2,u}) \leq g(\Delta) 
    \\
  g(\xi_{u+2,u}) \leq g(2\Delta)
    &\quad \quad g(\xi_{u-3,u}) \leq g(2\Delta) 
    \\
  \vdots \quad&\quad\quad\quad\quad \vdots
\end{align*}
Lastly, the above bounds also hold if we 
have $g(t_i-s_{2u+1})$ instead of $g(t_i-s_{2u})$.
Substituting these bounds back into \eqref{eq:num 0} 
and using \eqref{eq:big o notation} with $M_1=M_2=2$, we obtain:
\begin{align}
  \det(D^*D)
  &\leq \left[   
    8 \eta^6 
    + \left(1 + \frac{4}{\sigma^4} \right)
      \frac{2}{1-e^{-\frac{\Delta^2}{\sigma^2}}} 
  \right] \nonumber \\
  &\quad \cdot \prod_{u=2}^{k} \left[
    8 \eta^6 
    + \left(1 + \frac{4}{\sigma^4} \right)
      \frac{2}{1-e^{-\frac{\Delta^2}{\sigma^2}}} 
  \right] \nonumber \\
  &\quad \cdot \prod_{u=2}^{k} \left[
    32 \eta^4
      + \left( 
          \frac{4}{\sigma^4} 
          + \left(\frac{2}{\sigma^2} + \frac{4}{\sigma^4}  \right)^2
        \right)
        \frac{2}{1-e^{-\frac{\Delta^2}{\sigma^2}}} 
  \right] \nonumber \\
  &\leq \left[   
        8 + \left(1 + \frac{4}{\sigma^4} \right)
          \frac{2}{1-e^{-\frac{\Delta^2}{\sigma^2}}} 
    \right]^k
    \left[
      32 + \left( 
            \frac{4}{\sigma^4} 
          + \left(\frac{2}{\sigma^2} + \frac{4}{\sigma^4}  \right)^2
        \right)
      \frac{2}{1-e^{-\frac{\Delta^2}{\sigma^2}}} 
    \right]^{k-1} \nonumber \\
  &= F_1\left(\Delta,\frac{1}{\sigma}\right)^k \cdot
    F_2\left(\Delta,\frac{1}{\sigma}\right)^{k-1},
  \label{eq:final bnd on D2}
\end{align} 
where 
\begin{align}
  F_1\left(\Delta,\frac{1}{\sigma}\right) 
  &=
    8 + \left(1 + \frac{4}{\sigma^4} \right)
      \frac{2}{1-e^{-\frac{\Delta^2}{\sigma^2}}}, 
  \nonumber \\ 
  F_2\left(\Delta,\frac{1}{\sigma}\right) 
  &=
    32 + \left( 
      \frac{1}{\sigma^4} 
      +\frac{2}{\sigma^6} + \frac{2}{\sigma^8}  
    \right)
    \frac{16}{1-e^{-\frac{\Delta^2}{\sigma^2}}}.
  \label{eq:f1 and f2}
\end{align}
Note that the bound \eqref{eq:final bnd on D2} on $\det(D^*D)$
is the same for all $j=3,\ldots,m$. Combining \eqref{eq:bound det cc}
with \eqref{eq:final bnd on D2} in \eqref{eq:det tilde n}, we obtain:
\begin{equation}
  N_{\underline{l},j} \leq 
  \eta^{k+\frac12}
  C(f_0,f_1)^{\frac12} 
  \left(4k + 4 + \frac{4k}{\sigma^2}\right)^{\frac12}
  F_1\left(\Delta,\frac{1}{\sigma}\right)^{\frac{k}{2}}
  F_2\left(\Delta,\frac{1}{\sigma}\right)^{\frac{k-1}{2}}
  \label{eq:bound nj}
\end{equation}
for $j = 3,\ldots,m$ if \eqref{eq:cond eta f} holds.

Finally, we need to upper bound $N_{\underline{l},j}$ 
for $j=2$ and $j=m+1$. For simplicity, consider $j=2$. Applying
the same assumptions and operations as in \eqref{eq:det tilde n},
we have
\begin{align}
  N_{\underline{l},2} 
  &= \eta^k
  \begin{vmatrix}
    f_{0} & O(\eta^3) & \cdots & O(\eta^3) & O(2\eta^2) & \cdots & O(\eta^3) \\
    \vdots \\
    0 & g(t_i-s_{2}) & \cdots & g(t_i-s_{2u}) 
      & -g'(\xi_{i,u}) & \cdots & O(\eta^3) \\
    0 & g'(t_i-s_{2}) & \cdots & g'(t_i-s_{2u}) 
      & -g''(\xi'_{i,u}) & \cdots & O(\eta^3/\sigma^2) \\
    \vdots \\
    f_{1} & O(\eta^3) & \cdots & O(\eta^3) & O(2\eta^2) & \cdots & 1
  \end{vmatrix},
  \label{eq:det tilde n j2}
\end{align}
where $\xi_{i,u}, \xi'_{i,u} \in [t_i-s_{2u}-\eta,t_i-s_{2u}]$.
We bound $N_{\underline{l},2}$ by using Hadamard's inequality 
(the more general version, see \cite{Gradshteyn2000}) and, 
after we use \eqref{eq:big o notation} with $M_1=M_2=2$, 
we obtain:
\begin{align}
  N_{\underline{l},2} 
  &\leq \eta^k
    \sqrt{f_0^2 + f_1^2} 
    \left[
      1 + 4(k+1)\eta^6 + 4k \frac{\eta^6}{\sigma^4}
    \right]^{\frac12}
  \nonumber \\
  &\quad \cdot
    \prod_{u=1}^k \left[
      8\eta^6 + \sum_{i=1}^k g(t_i-s_{2u})^2
        + \sum_{i=1}^k g'(t_i-s_{2u})^2
    \right]^{\frac12}
  \nonumber \\
  &\quad \cdot
    \prod_{u=1}^k \left[
      32\eta^4
      + \sum_{i=1}^k g'(\xi_{i,u})^2
      + \sum_{i=1}^k g''(\xi'_{i,u})^2
    \right]^{\frac12},
\end{align}
and, by applying the bounds on the sums 
and $\eta \leq 1$, we obtain:
\begin{equation}
  N_{\underline{l},2} 
  \leq \eta^k
    \sqrt{f_0^2 + f_1^2} 
    \left(
      4k + 5 + \frac{4k}{\sigma^4}
    \right)^{\frac12}
    F_1\left(\Delta,\frac{1}{\sigma}\right)^{\frac{k}{2}}
    F_2\left(\Delta,\frac{1}{\sigma}\right)^{\frac{k}{2}},
  \label{eq:bound n2}
\end{equation}
for $F_1$ and $F_2$ defined as in \eqref{eq:f1 and f2}, 
and note that the same bound also holds for $N_{\underline{l},m+1}$.

To conclude, from \eqref{eq:bound nj} and \eqref{eq:bound n2}, we can derive
a general bound valid for all $j$:
\begin{equation}
  N_{\underline{l},j} 
  \leq \eta^k
  C(f_0,f_1)^{\frac12}
  \left(4k + 5 + \frac{4k}{\sigma^4}\right)^{\frac12}
  F_{\max}\left(\Delta,\frac{1}{\sigma}\right)^k
  \label{eq:final bound numerator}
\end{equation}
for all $j=2,\ldots,m+1$ if \eqref{eq:cond eta f} holds,
where
\begin{equation}
  F_{\max}\left(\Delta,\frac{1}{\sigma}\right) 
  = 
  \left(
    8 + \left(1 + \frac{4}{\sigma^4} \right)
      \frac{2}{1-e^{-\frac{\Delta^2}{\sigma^2}}}
  \right)^{\frac12}
  \left(
    32 + \left( 
      \frac{1}{\sigma^4} 
      +\frac{2}{\sigma^6} + \frac{2}{\sigma^8}  
    \right)
    \frac{16}{1-e^{-\frac{\Delta^2}{\sigma^2}}}
  \right)^{\frac12}.
\end{equation}

\subsection{Condition on $\eta$
    \label{sec:cond on eta}}

We now return to the condition \eqref{eq:cond eta unfinished} that $\eta$ 
must satisfy so that our application of Lemma \ref{lem:determinant bounds} is valid.
Since A is positive definite, we have $\norm{A^{-1}}_2 \leq 1/\lambda_{\min}(A)$. Using this,
we obtain:
\begin{equation}
  \rho(A^{-1}(A'+C')) 
  \leq \| A^{-1} (A'+C')\|_2 
  \leq \|A^{-1}\|_2 \|A' + C'\|_2
  \leq \frac{1}{\lambda_{\min}(A)} \cdot
    \|A'+C'\|_F,
\end{equation}
and 
\begin{align}
  A'+C' = \begin{bmatrix}
    0  & \cdots & O(\eta^2) & O(2\eta) & \cdots & O(\eta^2) \\
    \vdots &    & \vdots & \vdots & & \vdots \\
    O(\eta^2)   & \cdots & \frac{t_u-s_{2u}}{\eta} g'(\xi^*_{i,u}) 
        & - \frac{t_u-s_{2u}}{\eta} g''(\xi^{*\prime}_{i,u}) + \frac12 g''(\xi_{i,u}) 
        & \cdots & O(\eta^2) \\
    O(\eta^2/\sigma^2) & \cdots &  \frac{t_u-s_{2u}}{\eta}  g''(\xi^{*\prime}_{i,u}) 
        & - \frac{t_u-s_{2u}}{\eta} g'''(\xi^{*\prime\prime}_{i,u}) + \frac12 g'''(\xi'_{i,u}) 
        & \cdots & O(\eta^2/\sigma^2) \\
    \vdots &    & \vdots & \vdots & & \vdots \\
    O(\eta^2)  & \cdots & O(\eta^2) & O(2\eta) & \cdots & 0
  \end{bmatrix},
\end{align}
where
\begin{align}
  &\xi_{i,u}, \xi'_{i,u} \in [t_i - s_{2u} - \eta, t_i - s_{2u}],
  \nonumber \\
  &\xi^{*}_{i,u}, \xi^{*\prime}_{i,u},\xi^{*\prime\prime}_{i,u} 
    \in [t_i - t_u - |t_u - s_{2u}|, t_i - t_u + |t_u - s_{2u}|].
  \label{eq:xis}
\end{align}
for all $i,u = 1,\ldots,k$. 

Because the eigenvalues of $A$ are $1, 1$ and the eigenvalues of $B$ and
$\lambda_{\min}(B) \geq F_{\min}(\Delta, \frac{1}{\sigma})$ 
with $0 < F_{\min}(\Delta, \frac{1}{\sigma}) < 1$, then we have that 
that $\lambda_{\min}(A) \geq F_{\min}\left(\Delta,\frac{1}{\sigma}\right)$.
Moreover, after applying \eqref{eq:big o notation} with $M_1=M_2=2$, we have that:
\begin{align}
  \|A'+C'\|_F^2
  &\leq 8(k+1)\eta^4 + 8k\frac{\eta^4}{\sigma^4}
    \qquad \text{ (from the first and last columns) } \nonumber \\
  &+ 8k\eta^4 + 32 k\eta^2 +
    \qquad \text{ (from the first and last rows) } \nonumber \\
  &+ \sum_{i,u=1}^k g'(\xi^*_{i,u})^2 
    + \sum_{i,u=1}^k g''(\xi^{*\prime}_{i,u})^2
    \nonumber \\
  &+ \sum_{i,u=1}^k \left(- \frac{t_u - s_{2u}}{\eta} g''(\xi^{*\prime}_{i,u}) 
      + \frac12 g''(\xi_{i,u}) \right)^2
    + \sum_{i,u=1}^k \left(- \frac{t_u - s_{2u}}{\eta} g'''(\xi^{*\prime\prime}_{i,u}) 
      + \frac12 g'''(\xi'_{i,u}) \right)^2.
\end{align}
We upper bound this using the inequalities 
in \eqref{eq:two sums xi} and, similarly, for $g'''$:
\begin{equation}
  \sum_{i=1}^k g'''(\xi_{i,u})^2 \leq 
  \left(
    \frac{12}{\sigma^4} + \frac{8}{\sigma^6}
  \right)^2 \cdot
  \frac{2}{
    1 - e^{-\frac{\Delta^2}{\sigma^2}}
  },
\end{equation}
and note that these inequalities hold for all numbers $\xi$ 
in \eqref{eq:xis}.
By expanding the parentheses and applying Cauchy-Schwartz to the products,
and using that $|t_u - s_{2u}| \leq \eta$, we obtain:
\begin{align}
  \|A' + C'\|_F^2 
  &\leq
  8(2k+1)\eta^4 + 32 k\eta^2 + 8k\frac{\eta^4}{\sigma^4}
  \nonumber \\ &\quad
  + k \cdot \frac{4}{\sigma^4} \cdot \frac{2}{1-e^{-\frac{\Delta^2}{\sigma^2}}}
  + k \cdot \left(\frac{2}{\sigma^2}+\frac{4}{\sigma^4}\right)^2
        \cdot \frac{2}{1-e^{-\frac{\Delta^2}{\sigma^2}}}
  \nonumber \\ &\quad
  + \frac{9k}{4} \cdot \left(\frac{2}{\sigma^2}+\frac{4}{\sigma^4}\right)^2
        \cdot \frac{2}{1-e^{-\frac{\Delta^2}{\sigma^2}}}
  \nonumber \\ &\quad
  + \frac{9k}{4} \cdot \left(
    \frac{12}{\sigma^4} + \frac{8}{\sigma^6}
  \right)^2 \cdot
  \frac{2}{
    1 - e^{-\frac{\Delta^2}{\sigma^2}}
  }
\end{align}
Using the assumption that $\eta \leq \sigma^2$ to 
write $\frac{\eta^4}{\sigma^4} \leq \sigma^4$, and then 
by applying
$\eta \leq 1$ and $\sigma \leq \sqrt{2}$,
we can write
\begin{align}
  \|A' + C'\|_F^2
  \leq
  80k + 8 + k P\left(\frac{1}{\sigma}\right)
  \frac{2}{
    1 - e^{-\frac{\Delta^2}{\sigma^2}} 
  },
\end{align}
where $P\left(\frac{1}{\sigma}\right)$ is a polynomial in $\frac{1}{\sigma}$
defined as follows:
\begin{align}
  P\left(\frac{1}{\sigma}\right) = 
  \frac{4}{\sigma^4} 
  + \frac{13}{4} \left(
      \frac{2}{\sigma^2} + \frac{4}{\sigma^4}
    \right)^2
  + \frac{9}{4} \left(
      \frac{12}{\sigma^4} + \frac{8}{\sigma^6}
    \right)^2.
\end{align}
Inserting the above observations in the condition \eqref{eq:cond eta unfinished}, we finally obtain
\begin{equation}
  \eta \leq 
  \frac{
    8 F_{\min}(\Delta,\frac{1}{\sigma})
  }{
    34(2k+2)
    \left(
      80k + 8 + k P\left(\frac{1}{\sigma}\right)
      \frac{2}{
        1 - e^{-\frac{\Delta^2}{\sigma^2}}
      }
    \right)^{\frac12}
  }.
  \label{eq:cond eta 2}
\end{equation}
With this choice of $\eta$, the condition \eqref{eq:cond eta unfinished} 
is satisfied.

\subsection{Bound for \eqref{eq:reminder of coeffs} }

Combining the results from \ref{sec:bnd on denom} 
and \ref{sec:bnd on num}, we arrive at
\begin{align}
  |b_j| & = 
    \frac{N_{\underline{l},j+1}}{N_{\underline{l},1}}
    \qquad \mbox{(see \eqref{eq:ratio N})} \nonumber\\
  & \le 
  \frac{
    \bar{C}(f_0,f_1)^{\frac54}
    \left(4k + 5 + \frac{4k}{\sigma^4}\right)^{\frac12}
  }{
    1 - \frac{\sqrt{e}}{2}
  }
  \left[
    \frac{
      F_{\max}\left(\Delta,\frac{1}{\sigma}\right)
    }{
      F_{\min}\left(\Delta,\frac{1}{\sigma}\right)^2
    }
  \right]^k,
  \qquad \mbox{(see \eqref{eq:final bound numerator} and \eqref{eq:final bound denominator})}
  \\
  \text{and}\quad |b_j^{\pi}|
  &\leq
  \frac{
    \left(\bar{C}(f_0,f_1) + 2k\right)^{\frac32}
  }{
    \eta \left(1 - \frac{\sqrt{e}}{2}\right)
  }
  \left[
    \frac{
      F_{\max}\left(\Delta,\frac{1}{\sigma}\right)
    }{
      F_{\min}\left(\Delta,\frac{1}{\sigma}\right)^2
    }
  \right]^k,
  \qquad \mbox{(see \eqref{eq:detcc pi})}
\end{align}
for every $j\in[m]$, provided that the 
conditions \eqref{eq:cond eta f} and \eqref{eq:cond eta 2} hold. 
Consequently, 
\begin{align}
  \|b\|_2 & = \sqrt{\sum_{j=1}^m b_j^2} 
  \leq
  \frac{
    \sqrt{
      (2k+2)
      \left(
        4k + 5 + \frac{4k}{\sigma^4}
      \right)
    }
  }{
    1 - \frac{\sqrt{e}}{2}
  }
  \bar{C}(f_0,f_1)^{\frac54}
  F\left(\Delta,\frac{1}{\sigma}\right)^k,
  \\
  \|b^{\pi}\|_2 & \leq
  \frac{
    \sqrt{
      2k+2
    }
  }{
    \eta \left(1 - \frac{\sqrt{e}}{2}\right)
  }
  \left(\bar{C}(f_0,f_1) + 2k\right)^{\frac32}
  F\left(\Delta,\frac{1}{\sigma}\right)^k,
\end{align}
where we write 
$
  F\left(\Delta,\frac{1}{\sigma}\right) =   
  \frac{
    F_{\max}\left(\Delta,\frac{1}{\sigma}\right)
  }{
    F_{\min}\left(\Delta,\frac{1}{\sigma}\right)^2
  }
$
if, $\sigma \leq \sqrt{2}$ 
and $\Delta > \sigma \sqrt{ \log{\left(3 + \frac{4}{\sigma^2} \right)}}$ and the
conditions \eqref{eq:cond eta f} and \eqref{eq:cond eta 2} hold.
This completes the proof of Lemma \ref{lem:bounds on b for Gaussian} since $m=2k+2$ by assumption.

\section{Proof of Lemma \ref{lem:determinant bounds}
  \label{sec:proof of determinant bounds lem}}

To begin with, let us note that
\begin{equation*}
  \det(A + \epsilon B) = \det(A) \det(I + \epsilon A^{-1} B)
\end{equation*}
We will now upper and lower bound the second term in this equality.
Denoting $C = A^{-1}B$, consider
\begin{equation*}
  \log \det (I + \epsilon C) =
  \sum_{i \in R} \log(1 + \epsilon \lambda_i(C)) +
  \sum_{i \in I} \log(1 + 2 \epsilon \operatorname{Re}(\lambda_i(C)) + \epsilon^2 | \lambda_i(C) |^2),
\end{equation*}
where 
\begin{align*}
  &R = \{i \in [m]: \operatorname{Im}(\lambda_i(C)) = 0 \} \text{ and } \\
  &I = \{i \in [m]: \lambda_{i_1}(C), \lambda_{i_2}(C) 
    \text{ are complex conjugate and } 
    \operatorname{Im}(\lambda_{i_1}(C)) \ne 0\}.
\end{align*}

\begin{enumerate}
  \item For $i \in R$, if $\epsilon < \frac{1}{|\lambda_i(C)|}$, 
    use apply Taylor expansion for $\log(1+x)$ and obtain:
    \begin{equation}
      \log(1 + \epsilon \lambda_i(C)) = \epsilon \lambda_i(C) - 
        \frac{\epsilon^2 \lambda_i^2(C)}{2 \xi^2_{i,\epsilon}}, 
        \quad \text{ where }
        \xi_{i,\epsilon} \in 
        \left[ 1 - \epsilon |\lambda_i(C)|, 1 +\epsilon |\lambda_i(C)| \right].
      \label{eq:i in r}
    \end{equation}
  
  \item For $i \in C$, we apply the same Taylor expansion and, writing 
      $y = 2 \epsilon \operatorname{Re}(\lambda_i(C)) + \epsilon^2|\lambda_i(C)|^2$,
    for $|y| < 1$ we obtain:
    \begin{equation*}
      \log(1 + y) = \frac{y}{\xi_{i,\epsilon}}, 
      \quad \text{ where }
      \xi_{i,\epsilon} \in \left[1-|y|, 1+|y|\right].
    \end{equation*}

    Then, we have that
    \begin{align*}
      2 \epsilon |\lambda_i(C)| + \epsilon^2|\lambda_i(C)|^2
      \leq 4\epsilon |\lambda_i(C)| < 1,
    \end{align*}
    where the first inequality is true 
    if $\epsilon \leq \frac{2}{|\lambda_i(C)|}$ and the second 
    inequality is true if $\epsilon < \frac{1}{4|\lambda_i(C)|}$. 
    From the condition on $\xi_{i,\epsilon}$ and 
    noting that 
      $|\operatorname{Re}(\lambda_i(C))| \le |\lambda_i(C)|$, 
    we have that
    \begin{equation}
      \label{eq:i in i}
      \xi_{i,\epsilon} \leq 1 + |y| 
      \leq 1 + 4 \epsilon |\lambda_i(C)|
      \quad\quad \text{ and } \quad\quad
      \xi_{i,\epsilon} \geq 1 - |y| 
      \geq 1 - 4 \epsilon |\lambda_i(C)| 
      \geq \frac12,
    \end{equation}
    where the last inequality holds 
    if $\epsilon \leq \frac{1}{8|\lambda_i(C)|}$. 
    Therefore, $\frac{1}{\xi_{i,\epsilon}} \leq 2$
    if $\epsilon \leq \frac{1}{8|\lambda_i(C)|}$.
\end{enumerate}
We now use \eqref{eq:i in r}, \eqref{eq:i in i}
and $\epsilon \leq \frac{1}{8|\lambda_i(C)|}$.
Let $\epsilon \leq \frac{1}{8|\lambda_i(C)|}$ 
for all eigenvalues $\lambda_i(C)$, then:
\begin{equation}
  \log \det(I + \epsilon C) = 
  \sum_{i \in R} \epsilon \lambda_i(C) 
  - \frac{\epsilon^2}{2} \sum_{i \in R} \frac{\lambda^2_i(C)}{\xi^2_{i,\epsilon}}
  + \sum_{i \in I} \frac{1}{\xi_{i,\epsilon}}
  \left(
    2 \epsilon \operatorname{Re}(\lambda_i(C)) 
    + \epsilon^2 | \lambda_i(C) |^2
  \right),
\end{equation}
and, by using that
  $\xi_{i,\epsilon} \geq 1 - \epsilon|\lambda_i(C)| \geq \frac78 $ 
for $i \in R$ and $\frac{1}{\xi_{i,\epsilon}} \leq 2 $ for $i \in C$
and the fact that the index $i \in I$ accounts for two eigenvalues, 
we obtain:
\begin{align}
  |\log \det(I + \epsilon C)| &\leq
    \epsilon \sum_{i \in R} |\lambda_i(C)| + 
    \frac{32 \epsilon^2}{49} \sum_{i \in R} |\lambda_i(C)|^2 +
    4 \epsilon \sum_{i \in I} |\lambda_i(C)| + 
    2 \epsilon^2 \sum_{i \in I} |\lambda_i(C)|^2
  \nonumber \\ 
  &\leq
  2 \epsilon \sum_{i=1}^m |\lambda_i(C)| 
  + \epsilon^2 \sum_{i=1}^m |\lambda_i(C)|^2 
  \nonumber \\ 
  &\leq
  m \epsilon \rho(C) (2 + \epsilon \rho(C))
  \nonumber \\ &\leq
  m \epsilon \rho(C) (2 + \frac18) 
  = \frac{17}{8} m \epsilon \rho(C)
  \nonumber \\ 
  &\leq
  \frac12
\end{align}
where the second last inequality holds if
$ \epsilon \leq \frac{1}{8 \rho(C)} $
and the last inequality holds if 
$ \epsilon \leq \frac{8}{34 m \rho(C)} $, 
which is also the dominating condition for $\epsilon$ 
in the proof if $m \geq 2$.

Now, note that for $|x| \leq \frac12$, we have that 
$e^x = 1 + x e^{\xi}$ for 
some $\xi \in [-|x|, |x|] \subseteq [-\frac12,\frac12]$
by Taylor expansion, and by taking $x = \log\det(I + \epsilon C)$, 
we obtain:
\begin{equation*}
  \det(I + \epsilon C) = e^x
  \leq 1 + |x| e^{\frac12}
  \leq 1 + \frac{17 \sqrt{e}}{8} m \epsilon \rho(C) 
  \leq 1 + \frac{\sqrt{e}}{2},
\end{equation*}
and similarly
\begin{equation*}
  \det(I + \epsilon C) = e^x
  \geq 1 - |x| e^{\frac12}
  \geq 1 - \frac{17 \sqrt{e}}{8} m \epsilon \rho(C) 
  \geq 1 - \frac{\sqrt{e}}{2}.
\end{equation*}
From the last two inequalities, by multiplying
by $\det(A)$, we obtain the result of our lemma.

\section{Proof of Lemma \ref{lem:lwr bnd on B} \label{sec:proof of lwr bnd on B}}

Recalling the definition of $B$ in \eqref{eq:def of B}, 
we apply Gershgorin disc theorem 
to find the discs $D(a_{ii}, \sum_{i \ne j} a_{ij})$ 
which contain the eigenvalues of $B$. 
Due to the structure of $H$, we consider two cases:
\begin{enumerate}
  \item On odd rows, the centre of the disc is $a_{\text{odd}_{ii}} = g(0) = 1$ and the radius is
    \begin{subequations}
      \begin{align}
        R_{\text{odd}_i} 
        &= \sum_{\substack{j=1\\j\ne i}}^{k} |g(t_i - t_j)| + |-g'(t_i - t_j)| \\
        &= \sum_{\substack{j=1\\j\ne i}}^{k} g(t_i-t_j)
            \left(1 + \frac{2|t_i-t_j|}{\sigma^2}\right), 
        \quad\quad \text{for} \quad i=1,\ldots,k.
      \end{align}
    \end{subequations}
  
  \item On even rows, the centre of the disc is $a_{\text{even}_{ii}} = -g''(0) = \frac{2}{\sigma^2}$
    and the radius is
    \begin{subequations}
      \begin{align}
        R_{\text{even}_i} 
        &= \sum_{\substack{j=1\\j \ne i}}^{k} |g'(t_i - t_j)| + |-g''(t_i - t_j)| \\
        &= \sum_{\substack{j=1\\j\ne i}}^{k} g(t_i-t_j)
          \left(\frac{2|t_i-t_j|}{\sigma^2} + 
            \left|-\frac{2}{\sigma^2} + \frac{4(t_i-t_j)^2}{\sigma^4}\right|
          \right),
        \quad\quad \text{for} \quad i=1,\ldots,k.
      \end{align}
    \end{subequations}
\end{enumerate}
Since the eigenvalues of $B$ are real, they are lower bounded by 
\begin{equation}
  \min_{i = 1,\ldots,k} 1 - R_{\text{odd}_{i}}
  \quad\text{ or }\quad
  \min_{i = 1,\ldots,k} \frac{2}{\sigma^2} - R_{\text{even}_{i}}.
\end{equation}
Because $|t_i - t_j| \leq 1$, we have that
\begin{equation}
  1 - R_{\text{odd}_i} \geq 
  1 - \left(1 + \frac{2}{\sigma^2} \right)
  \sum_{\substack{j=1\\j\ne i}}^{k} g(t_i - t_j)
  \label{eq:r odd tmp}
\end{equation}
Since $ \frac{4(t_i-t_j)^2}{\sigma^4} -\frac{2}{\sigma^2}
 > \frac{4 \Delta^2}{\sigma^4} -\frac{2}{\sigma^2} > 0 $
due to our assumptions on $\sigma$ (see \eqref{eq:simplifying assumptions}), we obtain
\begin{equation}
  \frac{2}{\sigma^2} - R_{\text{even}_{i}} \geq
  \frac{2}{\sigma^2} - \frac{4}{\sigma^4} 
  \sum_{\substack{j=1\\j\ne i}}^{k} g(t_i - t_j),
  \label{eq:r even tmp}
\end{equation}
for all $i = 1,\ldots,k$.
Using the fact that $g$ is decreasing and $|t_i - t_j| \geq |i-j| \Delta$, 
we obtain
\begin{subequations}
  \begin{align}
    \sum_{\substack{j=1\\j\ne i}}^{k} g(t_i - t_j)
    &\leq
    \sum_{\substack{j=1\\j\ne i}}^{k} g(|i-j|\Delta) \\
    &\leq
    2 \sum_{j=1}^{\infty} g(j \Delta) 
    = 2 \sum_{j=1}^{\infty} e^{-\frac{j^2 \Delta^2}{\sigma^2}}
    \leq 2 \sum_{j=1}^{\infty} 
      \left(e^{-\frac{ \Delta^2}{\sigma^2}}\right)^j \\
    &=
    \frac{2 e^{-\frac{\Delta^2}{\sigma^2}}}{1 - e^{-\frac{\Delta^2}{\sigma^2}}},
    \quad \quad \text{ for } i = 1,\ldots,k.
  \end{align}
\end{subequations}
The sum of of the series is valid because $e^{-\frac{\Delta^2}{\sigma^2}} < 1$. 
Combining this with \eqref{eq:r odd tmp} and \eqref{eq:r even tmp}, we obtain:
\begin{subequations}
  \begin{align}
    &1 - R_{\text{odd}_i} 
    \geq 
    1 - \left(1 + \frac{2}{\sigma^2} \right)
    \frac{2 e^{-\frac{\Delta^2}{\sigma^2}}}{1 - e^{-\frac{\Delta^2}{\sigma^2}}},
    \label{eq:bound r odd}
    \\
    &\frac{2}{\sigma^2} - R_{\text{even}_{i}} 
    \geq
    \frac{2}{\sigma^2} - \frac{4}{\sigma^4} 
    \frac{
      2 e^{-\frac{\Delta^2}{\sigma^2}}
    }{
      1 - e^{-\frac{\Delta^2}{\sigma^2}}
    },
    \label{eq:bound r even}
  \end{align}
  \label{eq:bounds r}
\end{subequations}
for all $i = 1,\ldots,k$.
By using the assumption that
$\Delta > \sigma \sqrt{ \log{\left(3 + \frac{4}{\sigma^2} \right)}}$, 
we can check that the lower bound 
in \eqref{eq:bound r odd} is  greater than zero, 
and by also using that $\sigma \leq \sqrt{2}$,
we can check that it is smaller than
than the lower bound in \eqref{eq:bound r even}.

To conclude, since all the eigenvalues of the matrix $B$ are real 
and in the union of the discs $D(1, R_{\text{odd}_{i}})$ 
and $D(\frac{2}{\sigma^2}, R_{\text{even}_{i}})$ for $i=1,\ldots,k$, then, 
by using the above observation and the lower bound 
in \eqref{eq:bound r odd}, we obtain a lower bound of all the eigenvalues of $B$:
\begin{align}
  \lambda_j(B) 
  \geq
    1 - \left(1 + \frac{2}{\sigma^2} \right)
    \frac{
      2 e^{-\frac{\Delta^2}{\sigma^2}}
    }{
      1 - e^{-\frac{\Delta^2}{\sigma^2}}
    },
    \quad\quad \forall j=1,\ldots,2k.
\end{align}
Note that we may be able to obtain better bounds given by \eqref{eq:bound r even} 
for $k$ eigenvalues if we scale the Gershgorin discs so that 
$D(1,R_{\text{odd}_{i}})$ and $D(1,R_{\text{even}_{i}})$ become disjoint.

\section{Proof of Lemma \ref{lem:in particular}}
\label{sec:in particular}

Because $\Delta > \sigma\sqrt{\log{\frac{5}{\sigma^2}}}$, we have that 
\begin{equation}
  e^{-\frac{\Delta^2}{\sigma^2}} < \frac{\sigma^2}{5},
  \label{eq:in partic 1}
\end{equation}
which implies that 
\begin{equation}
  \frac{1}{1 - e^{-\frac{\Delta^2}{\sigma^2}}} < \frac{5}{5 - \sigma^2}.
  \label{eq:in partic 2}
\end{equation}
Then, 
$
  - \frac{2e^{-\frac{\Delta^2}{\sigma^2}}}{1 - e^{-\frac{\Delta^2}{\sigma^2}}} 
  > - \frac{2 \sigma^2}{5 - \sigma^2}
$,
so
\begin{equation}
  1 - (1 + \frac{2}{\sigma^2})
  \frac{2e^{-\frac{\Delta^2}{\sigma^2}}}{1 - e^{-\frac{\Delta^2}{\sigma^2}}} 
  >
  \frac{1 - 3\sigma^2}{5 - \sigma^2}, 
  \quad \text{ i.e. } \quad
  F_{\min}\left(\Delta,\frac{1}{\sigma}\right) 
  >
  \frac{1 - 3\sigma^2}{5 - \sigma^2}.
  \label{eq:in partic 3}
\end{equation}
Similarly, we have that:
\begin{equation}
  8 + \left(1 + \frac{4}{\sigma^4}\right)
  \frac{2}{1 - e^{-\frac{\Delta^2}{\sigma^2}}}
  <
  \frac{
    -8 \sigma^6 + 50 \sigma^4 + 40
  }{
    \sigma^4 (5 - \sigma^2)
  }
  < 
  \frac{98}{
    \sigma^4 (5 - \sigma^2)
  }
  \label{eq:in partic 4}
\end{equation}
and
\begin{equation}
  32 + \left(\frac{1}{\sigma^4} 
        + \frac{2}{\sigma^6}
        + \frac{2}{\sigma^8}
      \right)
  \frac{16}{1 - e^{-\frac{\Delta^2}{\sigma^2}}}
  <
  \frac{
    -32\sigma^{10} + 160\sigma^8 + 80\sigma^4 + 160\sigma^2 + 160
  }{
    \sigma^8 (5 - \sigma^2)
  }
  <
  \frac{592}
  {
    \sigma^8 (5 - \sigma^2)
  },
  \label{eq:in partic 5}
\end{equation}
where in the last two inequalities we used $\sigma < 1$.
Combining \eqref{eq:in partic 3},\eqref{eq:in partic 4} 
and \eqref{eq:in partic 5}, we obtain \eqref{eq:in particular 0}:
\begin{equation}
  \frac{
    F_{\max}\left(\Delta,\frac{1}{\sigma}\right)
  }{
    F_{\min}\left(\Delta,\frac{1}{\sigma}\right)^2   
  }
  <
  \frac{c_1}
  {
    \sigma^6(1-3\sigma^2)^2 
  }.
  \label{eq:in partic 6}
\end{equation}
Then, using that $\bar{f} < 1$ and that 
  $ 
    \frac{\bar{C}^{\frac54}}{\bar{f}} 
    \le 
    \frac{(\bar{C}+2k)^{\frac32}}{\bar{f}}, 
    (\bar{C}+2k)^{\frac32}
    \le 
    \frac{(\bar{C}+2k)^{\frac32}}{\bar{f}},    
  $
we have that:
\begin{align}
  \left(
    (6 + \frac{2}{\bar{f}})
    \sqrt{4k + 5 + \frac{4k}{\sigma^4}}
    \bar{C}^\frac54
    +
    \frac{6}{\eta}(\bar{C}+2k)^\frac32
  \right)
  \frac{\sqrt{2k+2}}{1-\frac{\sqrt{e}}{2}}
  &<
  c_2 \frac{(\bar{C}+2k)^\frac32}{\bar{f}} \sqrt{2k+2}
  \frac{\eta\sqrt{4k \sigma^4 + 5\sigma^4 + 4k}+\sigma^2}{\eta\sigma^2}
  \nonumber \\
  &<
  c_2 \frac{(\bar{C}+2k)^\frac32}{\bar{f}} \sqrt{2k+2}
  \cdot 
  \frac{\sqrt{8k+5} + 1}{\eta\sigma^2}
  \quad\quad (\sigma < 1)
  \nonumber \\
  &<
  c_2 \frac{(\bar{C}+2k)^\frac32}{\bar{f}}
  \frac{k}{\eta\sigma^2},
  \label{eq:in partic 7}
\end{align}
for a large enough constant $c_2$ 
and, from \eqref{eq:in partic 6} and \eqref{eq:in partic 7}, 
we obtain \eqref{eq:in particular 1}. Similarly we show 
that \eqref{eq:in particular 12} holds:
\begin{align}
  \frac{
    \sqrt{
      (2k+2)
      \left(4k + 5 + \frac{4k}{\sigma^4}\right)
    }
  }{
    1 - \frac{\sqrt{e}}{2}
  }
  \frac{
    \bar{C}(f_0,f_1)^\frac54}{
    \bar{f}
  }
  &<
  c_5 \frac{\sqrt{(k+1)(4k\sigma^4+5\sigma^4+4k)}}{\sigma^2} 
  \frac{\bar{C}^\frac54}{\bar{f}}
  \nonumber \\
  &< 
  c_5 \cdot \frac{\sqrt{(k+1)(8k+5)}}{\sigma^2} 
  \frac{\bar{C}^\frac54}{\bar{f}}
  \quad\quad (\sigma < 1)
  \nonumber \\
  &<
  c_5 \frac{k}{\sigma^2} \frac{\bar{C}^{\frac54}}{\bar{f}}.
\end{align}
To show that \eqref{eq:main cond eta} is satisfied 
if \eqref{eq:cond eta partic} holds,
from \eqref{eq:in partic 2} and \eqref{eq:in partic 3}, we have that:
\begin{align}
  \frac{
    8 F_{\min}(\Delta,\frac{1}{\sigma})
  }{
    34(2k+2)
    \left(
      80k + 8 + k P\left(\frac{1}{\sigma}\right)
      \frac{2}{
        1 - e^{-\frac{\Delta^2}{\sigma^2}}
      }
    \right)^{\frac12}
  }
  &>
  \frac{
    8(1 - 3\sigma^2)
  }{
    34(5-\sigma^2)(2k+2)
    \left(
      80k + 8 + 
      \frac{
        10kP\left(\frac{1}{\sigma}\right)
      }{
        5 - \sigma^2
      }
    \right)^\frac12
  }
  \nonumber \\
  &=
  \frac{
    8(1 - 3\sigma^2)
  }{
    34(2k+2)
    \left(
      (80k + 8)(5-\sigma^2)^2 + 
      10(5-\sigma^2)kP\left(\frac{1}{\sigma}\right)
    \right)^\frac12
  }
  \nonumber \\
  &>
  \frac{
    8(1 - 3\sigma^2)
  }{
    34(2k+2)
    \left(
      2000k + 200 + 
      50kP\left(\frac{1}{\sigma}\right)
    \right)^\frac12
  } 
  \quad\quad (5-\sigma^2 < 5)
  \nonumber \\
  &>
  \frac{
    8(1-3\sigma^2)\sigma^6
  }{
    34(2k+2)
    \left(
      2000k\sigma^{12} + 200\sigma^{12} + c_3 k
    \right)^\frac12
  }
  \quad\quad 
  \left(
    P\left(\frac{1}{\sigma}\right) < \frac{c_3}{\sigma^{12}}
  \right)
  \nonumber \\
  &>
  \frac{
    8(1-3\sigma^2)\sigma^6
  }{
    34(2k+2)
    \left(
      c_3 k + 200 
    \right)^\frac12
  }
  \quad\quad (\sigma < 1)
  \nonumber \\
  &>
  c_3 \cdot
  \frac{
    \sigma^6 (1-3\sigma^2)
  }{
    (k+1)^{\frac32}
  },
\end{align}
and similarly:
\begin{align}
  \frac{
    \bar{C}^{\frac16}
  }{
    \left(4k+4 + \frac{4k}{\sigma^2}\right)^\frac13  
  }
  &=
  \frac{
    \bar{C}^\frac16 \sigma^\frac23
  }{
    \left(4k\sigma^2 + 4\sigma^2 + 4k\right)^\frac13
  }
  >
  \frac{
    \bar{C}^\frac16 \sigma^\frac23
  }{
    (8k + 4)^\frac13
  } 
  >
  c_4 \cdot
  \frac{
    \bar{C}^\frac16 \sigma^\frac23
  }{
    (k + 1)^\frac13
  },
\end{align}
for some constants $c_3,c_4$.
Finally, from \eqref{eq:in partic 1} and \eqref{eq:in partic 2}, we also obtain:
\begin{equation}
  \frac{
    e^{-\frac{\Delta^2}{\sigma^2}} + e^{-\frac{2\Delta^2}{\sigma^2}}
  }{
    1 - e^{-\frac{\Delta^2}{\sigma^2}}
  }
  <
  \left(
    \sigma^2 + \frac{\sigma^4}{5}
  \right)
  \frac{1}{5 - \sigma^2}
  < \frac{16}{210},
\end{equation}
where, in the last inequality, we used $\sigma < \frac{1}{\sqrt{3}}$.
Furthermore, if $\lambda < \frac25$, then $1-2\lambda > \frac15$ and
\begin{equation}
  e^{-\frac{\Delta^2(1-2\lambda)}{\sigma^2}} 
  <
  \left( \frac{\sigma^2}{5} \right)^{1-2\lambda}
  <
  \frac{1}{
    15^{1-2\lambda}
  }
  <
  \frac{1}{15^{1/5}},
\end{equation}
so we can combine the last two inequalities to obtain:
\begin{equation}
  1 - 
  \frac{
    e^{-\frac{\Delta^2}{\sigma^2}} + e^{-\frac{2\Delta^2}{\sigma^2}}
  }{
    1 - e^{-\frac{\Delta^2}{\sigma^2}}
  }
  - e^{-\frac{\Delta^2(1-2\lambda)}{\sigma^2}} 
  >
  1 - \frac{16}{210} - \frac{1}{15^{1/5}}
  > \frac13.
\end{equation}
Finally, to bound $e^{\frac{\Delta^2\lambda^2}{\sigma^2}}$, 
note from the definition of $\lambda$ and $\eta$ 
that $\lambda \Delta < \frac{\eta}{2}$, and from the assumption
that $\eta \leq \sigma^2$, we obtain that:
\begin{align}
  e^{\frac{\Delta^2\lambda^2}{\sigma^2}}
  <
  e^{\frac{\eta^2}{4\sigma^2}}
  \leq
  e^{\sigma^2/4}
  < c_5,
\end{align}
where we used that $\sigma < 1$.
Combining the last two inequalities, we obtain 
that \eqref{eq:in particular 2} holds for some constant $c_5 > 0$.

\section{Proof of Lemma \ref{lem:small epsilon}}
  \label{sec:proof of small epsilon}

In the proof of Lemma \ref{lem:Gaussian is Tstar sys} in 
Appendix \ref{sec:proof of lemma Gaussian is Tstar sys} we require 
that $f_0 \gg \bar{f}$, $f_0 \gg f_1$ and $f_0 \gg 1$. These conditions
come from equations \eqref{eq:special case N 2}, \eqref{eq:f0 n f1} 
and \eqref{eq:tstar_cond_f0} respectively. We can, therefore, 
fix $\bar{f}$ and $f_1$ such that $\bar{f} < 1$ and $1 < f_1 < f_0$, 
and give an expression of $f_0$
as a function of $\bar{f}$ as $\epsilon \to 0$. 
From equation \eqref{eq:special case N 2}, we have that
\begin{equation}
  \frac{f_0}{\bar{f}} 
  >
  \frac{N_{\underline{l},1}}{
    \min_{\tau_{\underline{l}} \in T_{\epsilon}^C} N_{1,1}(\tau_{\underline{l}})
  }.
  \label{eq:f0 fbar 2}
\end{equation}
which is required by the condition that $\det(M_N)$ 
in \eqref{eq:gauss ex} is positive when $\tau_{\underline{l}} \in T_{\epsilon}^C$.
While $N_{\underline{l},1}$ does not depend on $\epsilon$, we
can argue that, for $\epsilon \to 0$,  
the minimum in the denominator is lower bounded by $N_{1,1}(\tau_{\underline{l}})$, 
where
\begin{equation}
  \tau_{\underline{l}} \in \bar{T}_{\epsilon} =
  \{t_1 - \epsilon,t_1+\epsilon,\ldots,t_k-\epsilon,t_k+\epsilon\}.
  \label{eq:bar t def}
\end{equation}
Therefore, a sufficient condition to ensure \eqref{eq:f0 fbar 2} is:
\begin{equation}
  \frac{f_0}{\bar{f}} 
  >
  \frac{N_{\underline{l},1}}{
    \min_{\tau_{\underline{l}} \in \bar{T}_{\epsilon}} N_{1,1}(\tau_{\underline{l}})
  }.
\end{equation}
Note that  
$\min_{\tau_{\underline{l}} \in \bar{T}_{\epsilon}} 
  N_{1,1}(\tau_{\underline{l}}) \to 0$ as $\epsilon \to 0$,
since two rows of the determinant become equal.
More explicitly, let us assume, for simplicity, that the minimum over the finite set
of points $\bar{T}_{\epsilon}$ is attained $t_1 + \epsilon$.
Then, we subtract the row with $\tau_{\underline{l}}$ from
the first row of $N_{1,1}(\tau_{\underline{l}})$, and then expand the determinant along this
row. By taking $\tau_{\underline{l}}=t_1+\epsilon$, we obtain
\begin{equation*}
  N_{1,1}(t_1 + \epsilon) = \sum_{j=1}^{m} (-1)^{j+1}  N_{1,1,j} 
  \left[g(t_1-s_j) - g(t_1-s_j + \epsilon)\right],
\end{equation*}
where the minors $N_{1,1,j}$ are fixed (i.e. independent of $\epsilon$).
Therefore, as $\epsilon \to 0$, we have that $N_{1,1}(t_1 + \epsilon) \to 0$, 
with $N_{1,1}(t_1 + \epsilon) > 0, \forall \epsilon$. Then, everything else
in $N_{1,1}(t_1 + \epsilon)$ being fixed, there exists $\epsilon_0 >0$
such that 
\footnote{
    To see this, take $f(\epsilon) = N_{1,1}(t_1 + \epsilon)$
    and we know that $f$ is continuous, $f(0)=0$ and $f(\epsilon)>0,\forall \epsilon$.
    If $f(\epsilon) = C \epsilon^2$  on $[0,\epsilon']$ for
    some $\epsilon' > 0$ and $C > 0$, which we show later in the proof (see \eqref{eq:mn eps}),
    then take $B = \min_{\epsilon \ge \epsilon'} f(\epsilon)$, and then
    there exists $\epsilon_0 \leq \epsilon'$ such that
    $f(\epsilon) \leq B, \forall \epsilon < \epsilon_0$.
    So we have that 
    $f(\epsilon) \leq \min_{\tau \geq \epsilon} f(\tau),
    \forall \epsilon < \epsilon_0$, 
    which implies that 
    $\min_{\tau \ge \epsilon} f(\tau) = f(\epsilon),
    \forall \epsilon < \epsilon_0$.
}
\begin{equation}
  \label{eq:min n11 at eps}
  \min_{\tau_{\underline{l}} \in T_{\epsilon}^C} 
  N_{1,1}(\tau_{\underline{l}})
  = N_{1,1}(t_1 + \epsilon),
  \quad \forall \epsilon < \epsilon_0.
\end{equation}
We will now find the exact rate at which $N_{1,1}(t_1+\epsilon) \to 0$
for $\epsilon < \epsilon_0$.
In the row with $\tau_{\underline{l}}$ in $N_{1,1}(t_1+\epsilon)$,
we Taylor expand the entries in the columns $j=1,\ldots,m$ as follows:
\begin{equation}
  g(\tau_{\underline{l}}-s_j) = g(t_1-s_j+\epsilon)
  = g(t_1-s_j) + \epsilon g'(t_1-s_j) 
  + \frac{\epsilon^2}{2} g''(\xi_j),
\end{equation}
for some $\xi_j \in [t_1-s_j, t_1-s_j+\epsilon]$, 
and note that $\xi_j \to t_1-s_j$ as $\epsilon \to 0$.
Then 
\begin{align}
  N_{1,1}(t_1+\epsilon)
  &= \epsilon \left|\begin{array}{ccc}
    g\left(t_1-s_{1}\right) & \cdots & g\left(t_1-s_{m}\right)\\ 
    g'\left(t_{1}-s_{1}\right) & \cdots &  g'\left(t_{1}-s_{m}\right)\\
    \vdots & & \vdots \\
    g'(t_1-s_1)+\frac{\epsilon}{2}g''(\xi_1) & 
          \cdots & g'(t_1-s_m)+\frac{\epsilon}{2}g''(\xi_m)   \\
    \vdots & & \vdots \\
    g\left(t_k-s_{1}\right) & \cdots & g\left(\tau-s_{m}\right)\\
    g'\left(t_{k}-s_{1}\right) & \cdots &  g'\left(t_{k}-s_{m}\right)\\
    g\left(1-s_{1}\right) & \cdots & g\left(1-s_{m}\right)
  \end{array}\right|
  \qquad \text{(subtract 1st row from } \tau_{\underline{l}} \text{ row)}
  \nonumber \\
  &=  \frac{\epsilon^2}{2} \left|\begin{array}{ccc}
    g\left(t_1-s_{1}\right) & \cdots & g\left(t_1-s_{m}\right)\\ 
    g'\left(t_{1}-s_{1}\right) & \cdots &  g'\left(t_{1}-s_{m}\right)\\
    \vdots & & \vdots \\
    g''(\xi_1) & \cdots & g''(\xi_m)   \\
    \vdots & & \vdots \\
    g\left(t_k-s_{1}\right) & \cdots & g\left(\tau-s_{m}\right)\\
    g'\left(t_{k}-s_{1}\right) & \cdots &  g'\left(t_{k}-s_{m}\right)\\
    g\left(1-s_{1}\right) & \cdots & g\left(1-s_{m}\right)
  \end{array}\right|  
  =: \frac{\epsilon^2}{2} N_{1,1}^{\epsilon},
  \qquad \text{(subtract 2nd row from } \tau_{\underline{l}} \text{ row)}
  \label{eq:mn eps}
\end{align}
for $\epsilon < \epsilon_0$
and note that swapping the $\tau_{\underline{l}}$ row with the third row involves
an even number of adjacent row swaps, so the sign of the determinant remains the same.
Also, for $\epsilon < \epsilon_0$, we have that:
\begin{equation}
  N_{1,1}^{\epsilon} \to 
  \left|\begin{array}{ccc}
    g\left(t_1-s_{1}\right) & \cdots & g\left(t_1-s_{m}\right)\\ 
    g'\left(t_{1}-s_{1}\right) & \cdots &  g'\left(t_{1}-s_{m}\right)\\
    g''(t_1-s_1) & \cdots & g''(t_1-s_m)   \\
    \vdots & & \vdots \\
    g\left(t_k-s_{1}\right) & \cdots & g\left(\tau-s_{m}\right)\\
    g'\left(t_{k}-s_{1}\right) & \cdots &  g'\left(t_{k}-s_{m}\right)\\
    g\left(1-s_{1}\right) & \cdots & g\left(1-s_{m}\right)
  \end{array}\right|
  =: N_{1,1}' > 0,
\end{equation}
where the last inequality is true because Gaussians form an 
extended T-system (see \cite{karlin1966tchebycheff}) and the determinant 
in the limit does not depend on $\epsilon$.
\footnote{
  Note that in the previous footnote we assumed 
  that $f(\epsilon) = C\epsilon^2$, where $C$ is independent of $\epsilon$, 
  but actually from \eqref{eq:mn eps} we have that
  $f(\epsilon) = N_{1,1}(t_1+\epsilon) = \frac{\epsilon^2}{2} N_{1,1}^{\epsilon}$,
  where $N_{1,1}^{\epsilon} \to N_{1,1}' > 0$ as $\epsilon \to 0$. Our conclusion that 
  the $N_{1,1}$ goes to zero at the rate $\epsilon^2$ does not change because,
  for small enough $E > 0$, there exists $\epsilon' > 0$ such that:
  \begin{equation*}
    0 < N_{1,1}' - E < N_{1,1}^{\epsilon} < N_{1,1}' + E,
    \quad \forall \epsilon < \epsilon',
  \end{equation*}
  so $f(\epsilon) \geq C_2 \epsilon^2$ for all $\epsilon < \epsilon'$, 
  where $C_2 = (N_{1,1}' - E)/2 > 0$ is independent of $\epsilon$. Then, using the argument in the
  previous footnote, there exists $\epsilon_0 > 0$ such that
  $\min_{\tau \geq \epsilon} f(\tau) \geq C_2 \epsilon^2$ 
  for all $\epsilon < \epsilon_0$, and therefore the obtain a version
  of \eqref{eq:eps zero final cond} where the factor in front of $1/\epsilon^2$ is 
  independent of $\epsilon$.
}

Substituting \eqref{eq:mn eps} and \eqref{eq:min n11 at eps} 
into \eqref{eq:f0 fbar 2}, and noting the the minimum in \eqref{eq:min n11 at eps}
can be attained at any $\tau_{\underline{l}} \in \bar{T}_{\epsilon}$ 
defined in \eqref{eq:bar t def} (not necessarily at $t_1+\epsilon$, which we 
assumed above for simplicity), we obtain:
\begin{equation}
  \frac{f_0}{\bar{f}} 
  >
  \frac{2 N_{\underline{l},1}}{
    \epsilon ^2 
    \min_{\tau_{\underline{l}} \in \bar{T_{\epsilon}}} N_{1,1}^{\epsilon}
  },
  \quad \forall \epsilon < \epsilon_0,
  \label{eq:eps zero final cond}
\end{equation}
which is the condition we must impose on $f_0 / \bar{f}$
so that $\det(M_N)>0$ for $\epsilon < \epsilon_0$
instead of \eqref{eq:dominated}.
Therefore, for $\epsilon < \epsilon_0$, we set
\begin{equation}
  f_0 = C_{\epsilon} \cdot \frac{\bar{f}}{\epsilon^2},
  \quad \text{ where } \quad
  C_{\epsilon} = \frac{3 N_{\underline{l},1}}{
    \min_{\tau_{\underline{l}} \in \bar{T}_{\epsilon}} N_{1,1}^{\epsilon}
  } 
  \quad \text{ and } \quad
  \lim_{\epsilon \to 0} C_{\epsilon} \in (0,+\infty).
  \label{eq:f0 fbar eps}
\end{equation}
Using that $f_1 < f_0$ and $1 < f_0$, we 
bound $\bar{C}(f_0,f_1)$ in \eqref{eq:c bar def}:
\begin{equation}
  \bar{C}(f_0,f_1) < \bar{c}_1 f_0^2,
  \label{eq:cbar bound}
\end{equation}
where $\bar{c}_1$ is a universal constant.
Finally, we insert \eqref{eq:f0 fbar eps} and \eqref{eq:cbar bound} 
into \eqref{eq:c1 eps def} and obtain:
\begin{align}
  C_1\left(\frac{1}{\epsilon}\right) 
  &< \frac{(\bar{c}_1 f_0^2 +2k)^\frac32}{\bar{f}}
  = \frac{(\bar{c}_1 C_{\epsilon}^2 \bar{f}^2/\epsilon^4 + 2k)^\frac32}{\bar{f}}
  \nonumber \\
  &< \frac{(\bar{c}_1 C_{\epsilon}^2 + 2k \epsilon^4)^\frac32}{\bar{f}} \cdot \frac{1}{\epsilon^6},
  \quad \forall \epsilon < \epsilon_0,
\end{align}
where we have used that $\bar{f} < 1$. Similarly, from \eqref{eq:c2 eps def}
we obtain
\begin{align}
  C_2\left(\frac{1}{\epsilon}\right) 
  &< \frac{\bar{c}_2 f_0^\frac52}{\bar{f}} 
  = \bar{c}_2 C_{\epsilon}^\frac52 \bar{f}^\frac32 
    \cdot \frac{1}{\epsilon^5}
  \nonumber \\
  &< \bar{c}_2 C_{\epsilon}^\frac52
    \cdot \frac{1}{\epsilon^5},
  \quad \forall \epsilon < \epsilon_0,
\end{align}
where $c_2$ is a universal constant, and this concludes the proof.

\end{document}